\documentclass[12pt]{article}
\usepackage{amsmath,amssymb}
\usepackage{amscd}
\usepackage{bm}
\usepackage[usenames]{color}
\newtheorem{tm}{Theorem}[section]
\newtheorem{lm}[tm]{Lemma}
\newtheorem{co}[tm]{Corollary}
\newtheorem{re}[tm]{Remark}

\newtheorem{exm}[tm]{Example}
\newtheorem{pr}[tm]{Proposition}
 
 \newenvironment{demo}[1]{\par\smallskip\par\begin{trivlist}
\item[]{\bf #1}\ }{\end{trivlist}\par\smallskip\par}
\newcommand{\Proof}{\begin{demo}{{\it Proof.\ }}}
\newcommand{\qed}{\end{demo}}
\newcommand{\toy}{\ \rule[0em]{0.5ex}{1.8ex}}
\newcommand{\QED}{\toy\end{demo}}

\newcommand{\la}{\langle}
\newcommand{\ra}{\rangle}
\newcommand{\nn}{\nonumber}
\newcommand{\III}{{\vert \kern-.10em \vert \kern-.10em \vert}}
\newcommand{\ve}{\varepsilon}

\newcommand{\al}{\alpha}

\makeatletter
 
 \@addtoreset{equation}{section}
\makeatother
 
\begin{document}
\setlength{\baselineskip}{15pt} 
%
\bibliographystyle{plain}
\title{
Short time full asymptotic expansion of hypoelliptic
heat kernel at the cut locus
}
\author{  Yuzuru \textsc{Inahama} and Setsuo \textsc{Taniguchi}
}
\date{   }
%
%
\maketitle


\begin{abstract}
In this paper we prove a short time asymptotic expansion of 
a hypoelliptic heat kernel 
on a Euclidean space and a compact manifold.
We study the ``cut locus" case, namely, the case where
energy-minimizing  paths
which  join the two points under consideration
form not a finite set, but a compact manifold.
Under mild assumptions we obtain an asymptotic expansion
of the heat kernel up to any order.
Our approach is probabilistic
and the heat kernel is regarded as
 the density of the law of a hypoelliptic diffusion process,
 which is realized as a unique solution of the corresponding 
 stochastic differential equation.
Our main tools are S. Watanabe's distributional 
Malliavin calculus and T. Lyons' rough path theory.
\end{abstract}

\section{Introduction}

Short time asymptotics of heat kernels 
is a huge topic in analysis, probability and geometry
and have been studied intensively and extensively.
In this paper we focus on the short time off-diagonal 
asymptotic expansion of heat kernels.  
We do not discuss the on-diagonal asymptotics, 
the estimates 
or the logarithmic asymptotics of heat kernels.
Even if we restrict our attention to the off-diagonal expansion,
so many papers have been published 
and our references are probably far from complete.

Known proofs of off-diagonal asymptotics are either analytic or probabilistic
and our proof is the latter.
The first probabilistic proof was given by Molchanov \cite{mol}
in the case of the standard heat kernel on a Riemannian manifold.
In it  pinned diffusion processes are used.
Bismut \cite{bi} was first to use Malliavin calculus, 
then followed by many others.
The Feynman-Kac representation is the key in these probabilistic approaches.


On a Riemannian manifold ${\cal N}$, 
if $x$ and $x'$ are outside the cut locus, 
the following asymptotic expansion of the standard heat kernel  
is well-known:
\begin{equation}\label{intro1.eq}
p_t (x,x') \sim e^{- d(x,x')^2 /2t} t^{-n/2}
(c_0 + c_1 t + c_2 t^2 +\cdots)
\qquad
\mbox{
as $t \searrow 0$.}
\end{equation}
Here, $n = \dim {\cal N}$ and $d(x,x')$ is the Riemannian distance between 
$x, x' \in {\cal N}$.
Ben Arous \cite{bena2} proved that the heat kernel of a
hypoelliptic diffusion process 
of a wide class on a Euclidean space
admits essentially the same expansion as above 
if $x$ and $x'$ are outside the cut locus.
In this case, however, the distance should be replaced by
 the sub-Riemannian distance.

When the two points are not near, 
this asymptotics becomes quite difficult.
Molchanov \cite{mol} already mentioned that 
there is an example of a Riemannian manifold ${\cal N}$
and  $x, x' \in {\cal N}$ in the cut locus
such that the factor
$t^{-n/2}$ in (\ref{intro1.eq}) must be replaced by 
$t^{-(n+n')/2}$  with some $n' \neq 0$.
Therefore, 
the off-diagonal asymptotics 
in the cut locus case is qualitatively different and 
much less understood than the case of ``near points."
(Note that we use the term ``cut locus" in a loose sense in this section.
What we precisely mean by the term will be given in 
Assumption {\bf (B2)} in the next section.)


Our present work
aims to study the short time  off-diagonal asymptotic
expansion for hypoelliptic heat kernels in the cut locus case
in a satisfactory general way.
It has the following three features.
To our knowledge,  
no other works
have all of these features simultaneously:
\begin{enumerate}
\item
The manifold and the hypoelliptic diffusion process on it are rather general.
In other words, this is not a study of special examples.
\item
The ``cut locus" case is studied.
More precisely, we mean by this that
the set of energy-minimizing paths (or controls)
which connect the two points under consideration is assumed to be
a compact manifold of finite dimension.
\item
The asymptotic expansion is full, that is, the polynomial part 
of the asymptotics is up to any order. 
\end{enumerate}

A recent paper by Barilari, Boscain and Neel  \cite{bbn}
satisfy (a) and (b), but 
only the leading term of the asymptotic expansion is obtained.
Chang and Li \cite{cl}  satisfies (b) and (c) for special examples.
(Many papers obtain the leading term in the  ``cut locus" case
for special examples. See \cite{bj, bw, bgg, bm, ga, sm, ue, uewa} among many others.)
Note that the result in \cite{bena2} is formulated on a Euclidean space,
but can be modified 
to the manifold case (if the manifold 
and the diffusion process are not too wild).
Hence, it satisfies (a) and (c).


On a Euclidean space, however, 
there are two famous works which satisfy (b), (c) and the latter half of (a). 
Both of them are probabilistic and use generalized versions of Malliavin calculus.
One is Takanobu and Watanabe \cite{tw}.
They use Watanabe's distributional Malliavin calculus 
developed in  \cite{wa_tata, wa}.
The other is Kusuoka and Stroock \cite{ks5}.
They use their version of generalized Malliavin calculus 
developed in  \cite{ks4}.
In this paper we use the former.

Though we basically follow Takanobu-Watanabe's argument,
the main difference is that
we use T. Lyons' rough path theory together, 
which is something like a deterministic version of the SDE theory.
The main advantage of using rough path theory is 
that 
while the usual It\^o map i.e. 
the solution map of an SDE  is discontinuous, 
the Lyons-It\^o map i.e.
the solution map of a rough differential equation (RDE)  is continuous.

This fact enables us to do  ``local analysis"  of the Lyons-It\^o map 
(for instance, 
restricting  the map on a neighborhood of its critical point and doing
Taylor-like expansion)
in a somewhat similar way we do in the Fr\'echet calculus.
Recall that in the standard SDE theory,
this type of local operation is often very hard and sometimes impossible,
due to the discontinuity of  the It\^o map.
For this reason, 
the localization procedure in \cite{tw}  looks so complicated
that it might be difficult to generalize their method
if rough path theory did not exist.


In this paper we first reprove and generalize 
the main result in \cite{tw} in the Euclidean setting
by using rough path theory.
Malliavin calculus for Brownian rough path was already 
studied in \cite{in1, in3}
in the proof of large deviation principle (LDP for short) for conditional measures.
It turned out that these two theories fit very well.

The strategy of our proof is as follows:
First,  the LDP above implies that 
contributions from
the complement  of the set of energy-minimizers (action-minimizers)
are negligibly small.
In other words, only a neighborhood of the set of energy-minimizers matters.
Next, around each energy-minimizer in the Cameron-Martin space, 
we take a small neighborhood with respect to the rough path topology
with certain nice properties. 
Since the set of energy-minimizers is compact, 
we need to compute contributions from finitely many of them only.
Then, on each of such small neighborhoods,
we do local analysis.
To prove the asymptotic expansion up to any order, 
we use
a modified version of Watanabe's asymptotic theory given in \cite{tw},
which can be regarded as a  ``localized" version 
of his asymptotic theory developed in \cite{wa, iwbk}.

Even in this Euclidean setting
many parts of the proof are technically 
improved thanks to rough path theory.
We believe that the following are worth mentioning:
{\rm (i)}~Large deviation upper bound 
(Theorem \ref{tm.ub_ldp} and Lemma \ref{lm.ldp_up}).
{\rm (ii)}~Asymptotic partition of unity (Section 6).
{\rm (iii)}~A Taylor-like expansion of the Lyons-It\^o map 
(Subsection 3.3)
and the uniform
exponential integrability lemma for the ordinary and the remainder terms
of the expansion
(Lemmas \ref{lm.tw64} and \ref{lm.tw65}).
{\rm (iv)}~Quasi-sure analysis for solutions of SDEs, 
although this is implicit in this paper. (See \cite{in1, in3} for details.)


Then, we study the manifold case. 
Recall that Malliavin calculus for a manifold-valued SDE was 
studied by Taniguchi \cite{ta}.
Thanks to rough path theory 
and the localized version of Watanabe's asymptotic theory, 
it is enough to localized in a small neighborhood 
of energy-minimizers in the geometric rough path space.
Restricted on such a  small neighborhood, 
the (Lyons-)It\^o map  corresponding to the given SDE on the manifold
can easily be transferred to one on a linear space.
For these reasons, without much technical effort 
the problem reduces to the one on a Euclidean space.
As a result, 
the proof of the asymptotic expansion
in the manifold setting is not so different from the one in the Euclidean setting.

Of course, there is a possibility that 
our main result can be proved without rough path theory,
but we believe that the theory is quite suitable for this problem
and gives us a very clear view (in particular, in the manifold case).

%
%
%
%
%


The structure of the paper is as follows:
In Section 2,
we provide setting, assumptions and our main results.
First, we discuss the Euclidean case.
Our assumptions are weaker that those in \cite{tw}
and hence our main result (Theorem \ref{tm.main}) in the Euclidean case 
generalizes the main result in \cite{tw}.
Next, we formulate a parallel problem on a compact manifold.
As mentioned before,
our main result (Theorem \ref{tm.main_mfd}) in the manifold case
satisfies all of (a), (b) and (c) above.
At the end of the section, a few examples are given.

Sections 3--7 are devoted to proving 
our main theorem in the Euclidean setting (Theorem \ref{tm.main}).
In Section 3, we gather various facts from relevant research areas 
such as Malliavin calculus,  rough path theory
and differential geometry.
All of them are already known basically, but some of them are not elementary.
They will be used in the proof of the asymptotic expansion.

In Section 4, we calculate the skeleton ODE 
driven by an energy-minimizing Cameron-Martin path
and show that
the solution satisfies  a naturally defined Hamiltonian ODE.
We look at this well-known argument from a viewpoint of rough path theory.
In Section 5, we prove the uniform exponential integrability of 
Wiener functionals which appear in the proof of the asymptotic expansion.
Since we use rough path theory, 
the Taylor-like  expansion of the Lyons-It\^o map is 
deterministic and the remainder terms satisfy  a simple and reasonable estimate.
This simplifies our proof.

In Section 6, we study an asymptotic partition of unity
which was first introduced in \cite{tw}.
Though it plays an important role in the proof of the asymptotics in  \cite{tw},
it is of very complicated form.
We construct a rough path version of it,
which is written in terms of Besov norms of Brownian rough path  
and looks simple and natural.

Section 7 is the main part of our proof of Theorem \ref{tm.main}.
Our main tool is a modified version of Watanabe's asymptotic theorem.
In Section 8, we prove our main result in the manifold case (Theorem \ref{tm.main_mfd}).


\section{Setting, assumptions and main results}

\subsection{Setting: The Euclidean case}
Let ${\cal W} =C_0 ([0,1], {\mathbb R}^r)$ be the set
of the continuous functions from $[0,1]$ to ${\mathbb R}^r$
which start at $0$. This is equipped with the usual sup-norm.
The Wiener measure on ${\cal W}$ is denoted by $\mu$.
We denote by 
\[
{\cal H}
=
\Bigl\{  h \in {\cal W} 
~\Big\vert~
 \mbox{absolutely continuous and } 
\|h\|^2_{{\cal H}} :=\int_0^1 |h^{\prime}_s|^2 ds <\infty
\Bigr\}
\]
the Cameron-Martin subspace of ${\cal W}$.
In some (non-probabilistic) literatures, the derivative of 
Cameron-Martin paths are called controls.
The triple $({\cal W}, {\cal H}, \mu)$ is called the classical Wiener space.
The canonical realization on ${\cal W}$
of $r$-dimensional
Brownian motion is denoted by 
$(w_t)_{0 \le t \le 1}=(w_t^1, \ldots, w_t^r)_{0 \le t \le 1}$.

Let $V_{i} \colon {\mathbb R}^d \to {\mathbb R}^d$ be a vector field 
with sufficient regularity ($0 \le i \le r$).
Precise conditions on $V_{i}$ will be specified later.
For  a small parameter $\ve \in (0,1]$,
we consider
the following SDE of Stratonovich type:
\begin{equation}\label{sdeX.def}
dX^{\ve}_t =\ve \sum_{i=1}^r  V_i ( X^{\ve}_t) \circ  dw_t^i  + \ve^2  V_0 (X^{\ve}_t)   dt
\qquad
\qquad
\mbox{with \quad $X^{\ve}_0 =x \in {\mathbb R}^d$.}
\end{equation}
When necessary, we will write $X^{\ve}_t = X^{\ve}(t, x, w)$ or $X^{\ve}(t, x)$
and sometimes write $\lambda^{\ve}_t =\ve^2 t$.
When $\ve =1$, we simply write $X_t =X^{1}_t$.
By the scaling property of Brownian motion, 
the laws of the processes
$(X^{\ve}_t )_{t \ge 0}$ and $(X_{\ve^2 t})_{t \ge 0}$ are the same.

Now we introduce the skeleton ODE which corresponds to SDE (\ref{sdeX.def}).
For a Cameron-Martin path $h \in {\cal H}$,
we consider the following controlled ODE:
\begin{equation}\label{ode.def}
d\phi_t = \sum_{i=1}^r  V_i ( \phi_t) dh_t^i 
\qquad
\qquad
\mbox{with \quad $\phi_0 =x \in {\mathbb R}^d$.}
\end{equation}
The solution will often be denoted by $\phi_t(h),~ \phi(t, x, h)$, etc.
Note the absence of the drift term in (\ref{ode.def}).

Let ${\cal V}$ be an $n$-dimensional linear subspace of ${\mathbb R}^d$ ($1 \le n \le d$)
and $\Pi_{{\cal V}} \colon {\mathbb R}^d \to {\cal V}$ be the orthogonal projection.
(For our purpose,  we may and sometimes will 
assume without loss of generality that ${\cal V} = {\mathbb R}^n \times \{ {\bf 0}_{d-n}\}$, 
where ${\bf 0}_{d-n}$ is the zero vector of ${\mathbb R}^{d-n}$.)
Set 
$Y^{\ve}_t = \Pi_{{\cal V}} (X^{\ve}_t)$ and $\psi_t =  \Pi_{{\cal V}} (\phi_t)$, 
which will often be denoted by $Y^{\ve}(t, x, w)$,
and
$\psi(t, x, h)$, respectively. 

For $a \in {\cal V}$,
define 
${\cal K}_a = \{h \in {\cal H}  \mid  \psi (1,x,h) =a \}$.
We set 
$d_a :=\min \{  \|h\|_{{\cal H}}  \mid  h \in  {\cal K}_a\}$
and 
${\cal K}_a^{min} =\{ h \in  {\cal K}_a \mid \|h\|_{{\cal H}} =d_a\}$.
It is known that
if ${\cal K}_a \neq \emptyset$, 
then $\inf\{  \|h\|_{{\cal H}} \mid h \in  {\cal K}_a\}$ actually attains a minimum and ${\cal K}^{min}_a \neq \emptyset$.
(We can see this from  the goodness of the rate function in a Schilder-type large 
deviation principle on the rough path space.)
We will basically assume $x \notin \Pi_{{\cal V}}^{-1} (a)$
and 
${\cal K}_a \neq \emptyset$, which imply $0 <d_a <\infty$.



\subsection{Assumptions: The Euclidean case}

In this subsection we introduce assumptions.
First, we impose two conditions on the coefficient vector fields.
The first one is on regularity of $V_i$.
\\
\\
{\bf (A1)}:  $V_{i} \colon {\mathbb R}^d \to {\mathbb R}^d$
is of $C^{\infty}$ with bounded derivatives of all order $\ge 1$ $(0 \le i \le r)$.
\\
\\
Note that $V_{i}$ itself may have linear growth.
When $V_{i}$ is also bounded,  $V_{i}$ is said to be of $C^{\infty}_b$.
($C^{k}_b$ is similarly defined for $1 \le k <\infty$.)
This assumption is very standard in Malliavin calculus
and most of the results in Malliavin calculus for SDEs are proved under this assumption.
Under {\bf (A1)}, a global solution $X^{\ve}$ of (\ref{sdeX.def}) uniquely exists
and $X^{\ve}_t$ is smooth in the sense of Malliavin calculus for any $(t, \ve) \in [0,1]^2$
and their Sobolev norm of any index
is bounded in $(t, \ve) \in [0,1]^2$.
Similarly, {\bf (A1)} is  a standard assumption for the skeleton ODE (\ref{ode.def}).
A unique global solution $\phi$ for any $h \in {\cal H}$ exists 
and,  moreover,  $\phi(\cdot\,  ,x, h)$ 
is   Fr\'echet smooth in $h \in {\cal H}$ for any $x$
under {\bf (A1)}.


Next we introduce  
a kind of H\"ormander's bracket generating 
condition not only on the  ``upper" space ${\mathbb R}^d$, 
but also on the "lower" space ${\cal V}$. 
We denote by $\sigma[ Y^{\ve}_1]$ and $\sigma[\psi_1]$
the Malliavin covariance matrix
of $Y^{\ve}_1$ and 
the deterministic Malliavin covariance matrix
of
$\psi_1$, respectively.
(Precise definitions will be given in Section 3.)

Set  
$$
\Sigma_1 =\{ V_i \mid 1 \le i \le r\} 
\qquad
\mbox{and} 
\qquad
\Sigma_k =\{ [V_i, W] \mid 0 \le i \le r, W \in \Sigma_{k-1}\} 
$$
for $k \ge 2$ recursively.
Note that the drift vector field is not involved in the definition of $\Sigma_1$. 
\\
\\
{\bf (A2)}:
We say that (A2) holds if either of the following two conditions holds:
\\
{\rm (i)}~
At the starting point $x \in {\mathbb R}^d$,  
$\cup_{k=1}^{\infty}  \{  W(x) \mid  W \in \Sigma_{k}\}$ spans 
${\mathbb R}^d$ in the sense of linear algebra.
\\
{\rm (ii)}~ There exists $L \in {\mathbb N}$ such that 
\begin{equation}\label{upHc}
\inf_{ x \in  {\mathbb R}^d }  \inf_{\eta \in {\cal V}, |\eta| =1} 
\sum_{k=1}^L \sum_{W \in \Sigma_k}
\la \Pi_{{\cal V}}  W (x), \eta \ra^2 >0.
\end{equation}
\\
\\
The first one is quite standard.
The latter one is called the uniform partial H\"ormander condition.
This condition {\bf (A2)} implies that $Y^{\ve}_t =  Y^{\ve}(t,x,w)$ is  non-degenerate in 
the sense of Malliavin calculus for any $t>0$ and $\ve \in (0,1]$.
Moreover, 
there exist positive constants $\nu$ and $C_p~(1 <p <\infty)$
such that
\begin{equation}\label{ks_bnd.ineq}
\| \det \sigma[ Y^{\ve}_1]^{-1} \|_{L^p} \le C_p \ve^{-\nu}
\qquad
\qquad
(1 <p <\infty,~0< \ve \le 1).
\end{equation}
(See Theorem 2.17 and Lemma 5.1
in Kusuoka and Stroock \cite{ks2}.)
Note that $\nu$ is independent of $p$.
Non-degeneracy implies that, if $\ve >0$ and $t>0$,
 $\delta_a (Y^{\ve}_t )$ is 
well-defined as a positive Watanabe distribution (i.e.
generalized Wiener functional),
where $\delta_a$ stands for the delta function at $a \in {\cal V}$,
and its generalized expectation ${\mathbb E}[\delta_a (Y^{\ve}_t )]$
equals the smooth density of the 
of the law of $Y^{\ve}_t$ 
with respect to the Lesbegue measure on ${\cal V}$
for every $a$.
(By the way, the partial partial H\"ormander condition only at the starting point 
does not imply the non-degeneracy of $Y^{\ve}_t$.
There is a simple counterexample.)


From here
we introduce assumptions on the subset ${\cal K}^{min}_a$ of energy minimizers 
($a \neq \Pi_{{\cal V}} x$).
\\
\\
{\bf (B1)}:  Assume that ${\cal K}_a \neq \emptyset$
and the deterministic Malliavin covariance matrix
$\sigma [\psi_1] (h)$ is non-degenerate for any $h \in {\cal K}^{min}_a$.
\\
\\
Note that 
non-degeneracy of $\sigma [\psi_1] (h)$ is equivalent to surjectivity of  
the tangent map $D \psi_1 (h) \colon {\cal H} \to {\cal V}$.
If we use terminology of sub-Riemannian geometry,
this condition loosely means that any element of ${\cal K}^{min}_a$
is a non-abnormal geodesics.
(Keep in mind that we are not in the framework of sub-Riemannian geometry, however.)

As we wrote, when ${\cal K}^{min}_a$  is not a finite set
 things usually get quite complicated.
To exclude the cases where ${\cal K}^{min}_a$  is not  nice, 
we assume the following condition:
\\
\\
{\bf (B2)}: 
${\cal K}^{min}_a$ is a smooth, compact  manifold of finite dimension $n'$
regularly embedded in ${\cal H}$.
\\
\\
In the above condition the case $n' =0$ is also allowed,
 which means that ${\cal K}^{min}_a$ is a finite set.
Note that compactness of  ${\cal K}^{min}_a$ is essential and will be used throughout this paper.
%
%
It will turn out in Remark \ref{re.top=} that {\bf (A1)} and {\bf (B1)} imply 
compactness of ${\cal K}^{min}_a$ with respect to ${\cal H}$-topology. 
From this viewpoint, the compactness assumption in {\bf (B2)} seems  natural.


Finally, we assume positivity of the Hessian of 
$I |_{{\cal K}_a}$ on ${\cal K}^{min}_a$,
where
 $I (h) := \| h \|^2_{{\cal H}} /2$, in an appropriate sense.
 Without the positivity our proof breaks down.
Since $I |_{{\cal K}_a}$ attains minimum on ${\cal K}^{min}_a$,
its Hessian 
is naturally defined on ${\cal K}^{min}_a$
without assuming additional structures.
The precise definition of the Hessian at $h \in {\cal K}^{min}_a$
is as follows:
For any $k \in T_h {\cal K}_a$
(which is equivalent to $k \in \ker D\psi_1 (h)$),
take a smooth curve $(-\tau_0, \tau_0) \ni \tau \mapsto c(\tau) \in {\cal K}_a$ 
such that 
$c(0) =h$ and $c^{\prime} (0) =k$. (Here, $\tau_0 >0$.
Because of {\bf (A1)} and the implicit function theorem, 
${\cal H}$-neighborhood of $h$ in ${\cal K}_a$ has a manifold structure.
Hence, the notion of  ``smooth" curve in the neighborhood is well-defined.)
Then, we set 
\begin{eqnarray*}
I^{ \prime\prime} (h) \la k,k \ra
&:=&
\frac{d^2}{d\tau^2} \Big|_{\tau =0}   I (c(\tau)).
\end{eqnarray*}
As is well-known, this  is independent of the choice of the curve $c$
and $I^{ \prime\prime} (h)$ becomes 
a bounded bilinear form on $T_h {\cal K}_a \times T_h {\cal K}_a$.
It is obvious that $I^{ \prime\prime} (h) \la k,k \ra =0$ if $k \in T_h {\cal K}^{min}_a$. 

We assume that the elements of
$T_h  {\cal K}^{min}_a$ are the  ``only directions" in $T_h {\cal K}_a$ 
such that $I^{ \prime\prime} (h)$ vanishes.
For simplicity, we will set ${\cal H}_0(h) := T_h {\cal K}_a  \cap (T_h {\cal K}^{min}_a )^{\bot}$,
which is a closed subspace of codimension $n+ n'$.
\\
\\
{\bf (B3)}: 
For any $h \in {\cal K}^{min}_a$ and 
any $k \in {\cal H}_0(h) \setminus \{0\}$,
$I^{ \prime\prime} (h) \la k,k \ra >0$.
\\
\\
Some equivalent conditions to {\bf (B3)} are known.
(See e.g. Bismut \cite{bi}
or Takanobu and Watanabe \cite{tw}.)
One obvious example is 
that {\bf (B3)} remains the same if we replace 
``$k \in {\cal H}_0(h) \setminus \{0\}$" 
by  ``$k \in T_h {\cal K}_a  \setminus T_h {\cal K}^{min}_a$".
Note also  that 
{\bf (B3)}  is also equivalent to the exponential integrability of 
a certain quadratic Wiener functional that appears in the proof of our main theorem.
(See Section 5.)

Another possible choice of sufficient condition in place of {\bf (B3)}  is the
``non-focality condition."
In Section 3, \cite{dfjv}, 
equivalence of  
{\bf (B3)}  and 
the non-focality condition is shown under slightly different assumptions from ours.
Although 
{\bf (B3)} is quite typical and does not look bad, 
some people may prefer the non-focality condition 
since it is written in terms of
finite dimensional mathematics, in particular, Hamiltonian ODEs.
However, we do not use the non-focality condition in this paper.


We introduce a family $G^{\ve} (w) = G(\ve, w)$
of scalar-valued ${\mathbb D}_{\infty}$-functionals ($0 <\ve \le 1$).
Main examples we have in mind are 
multiplicative functionals that appear in 
the Feynman-Kac(-It\^o) formulae.
Here, ${\mathbb D}_{\infty}$ stands for the space of test Wiener functionals.
\\
\\
{\bf (C1)}: For any $\ve \in (0,1]$, 
$G^{\ve}$ is a real-valued smooth Wiener functional, that is,
$G^{\ve}  \in {\mathbb D}_{\infty}$ and 
 the following asymptotic expansion holds for any $h \in {\cal H}$:
\[
G(\ve, w+ \frac{h}{\ve} ) 
\sim 
\Gamma_0 (h) 
+ 
\ve \Gamma_1 (w; h) + \ve^2 \Gamma_2 (w; h) +\cdots
\qquad
\mbox{in ${\mathbb D}_{\infty}$  as $\ve \searrow 0$}
\]
uniformly on $\{ h  \in {\cal H}\mid \|h\|_{{\cal H}}  \le \rho \}$ for any $\rho >0$.
Here, $\Gamma_0 (h)$ is a smooth function on ${\cal H}$
and $\Gamma_j (\,\cdot\,; h) \in {\mathbb D}_{\infty}$ depends smoothly on $h \in {\cal H}$ ($j \ge 1$).
\\
\\
In the asymptotic expansion of  heat kernels,
the odd order terms of $\ve$ vanish.
%
(As usual we use a scale change $t= \ve^2$.)
To formulate and prove this phenomenon, 
we set the following assumption:
\\
\\
{\bf (C2)}: For any $h \in {\cal K}_a^{min}$,
$\Gamma_j (\,\cdot\,; h)\in {\mathbb D}_{\infty}$ 
is an even (resp. odd) Wiener functional in $w$-argument
if $j$ is even (resp. odd).
\\
\\
We say a (generalized) Wiener functional $F$ is even (resp. odd)
if $F(-w) =F(w)$ (resp. $F(-w) =-F(w)$).
These are well-defined since the map $w \mapsto -w$ leaves
$\mu$ invariant.
If $G^{\ve} \equiv 1$, then {\bf (C2)} is clearly satisfied.

\begin{re}\normalfont\label{rem_ee}
The reader may wonder why we introduce the projection $\Pi_{{\cal V}}$
in our formulation.
In the early days of Malliavin calculus, 
this type of assumption was common in the study of 
``partially hypoelliptic" problems (See \cite{ta, ks2} and references within.)
This formulation probably has its origin 
in Eells-Elworthy's construction of Brownian motion on a Riemannian 
manifold (Example \ref{exm.riem}).
Recently, this kind of formulation attracted attention from mathematical finance,
too.
(See \cite{dfjv} for instance.)
\end{re}



\subsection{Main result: The Euclidean case}

In this subsection we present our main result in the Euclidean setting.
The proof of this theorem will be given in Sections \ref{sec.apu}--\ref{sec_pf}.
An explicit expression of the leading term $c_0$ will be given in (\ref{tp_trm.eq}).
If $\Gamma_0 $ is non-negative and not identically zero
on ${\cal K}_a^{min}$, then $c_0 >0$.
Note that the left hand side of the asymptotic expansion (\ref{eq.tusika}) below
is the generalized expectation.
\begin{tm}\label{tm.main}
Let $x \in {\mathbb R}^d$ and $a \in {\cal V}$ such that $\Pi_{{\cal V}} (x) \neq a$.
Assume {\bf (A1)},  {\bf (A2)}, {\bf (B1)}, {\bf (B2)}, {\bf (B3)} and {\bf (C1)}.
Then, we have the following asymptotic expansion:
\begin{equation}\label{eq.tusika}
{\mathbb E} [G(\ve, w) \delta_a (Y_1^{\ve}) ]
 \sim
e^{ - d_a^2 /2 \ve^2}
 \ve^{-(n+ n')}
 (c_0 + c_1 \ve + c_2 \ve^2 + \cdots)
 \qquad
 \mbox{ as $\ve\searrow 0$}
\end{equation}
for certain constants $c_j \in {\mathbb R}~(j \ge 0)$.
If we assume {\bf (C2)} in addition, then $c_{2j +1} =0~(j \ge 0)$. 
\end{tm}

\begin{re}\normalfont
The only difference between our Thorem \ref{tm.main}
and Theorem 5.1 in \cite{tw}
is the 
H\"ormander's  bracket generating condition {\bf (A2)}.
In \cite{tw} 
the stronger version of
bracket generating condition is assumed at every $x \in {\mathbb R}^d$
on the upper space.
(By  ``stronger" we mean the drift vector field $V_0$ is not involved in the condition.)

On the other hand, 
we assume the weaker version of 
bracket generating condition. 
Moreover, the condition is either  ``on the lower space" 
or  ``at the starting point only on the upper space" as in {\bf (A2)}.
Therefore, we believe we made sizable relaxation of 
the bracket generating condition.

Put in a broader context, however, 
whether this improvement is large or not is somewhat unclear.
In other words, we do not know how many of the newly allowed 
examples by this relaxation 
satisfy the other assumptions.
(For instance, under the weaker version of H\"ormander's  condition,
{\bf (B1)} may not hold very often.)
It needs further investigations. 

On the other hand, the relaxation of the H\"ormander-type  condition 
to the  ``partial" version
(i.e. the bracket generating condition on the lower space)
is crucial for our geometric purpose.
Without it 
we could not even treat the case of 
Brownian motion on a Riemannian manifold.
(See Remark \ref{rem_ee} and Examples \ref{exm.riem}--\ref{exm.subR}.)
\end{re}


The formulation of Theorem \ref{tm.main},
(in particular, the generalized expectation) 
may look too general 
for non-experts of Watanabe's distributional  Malliavin calculus.
Therefore, we give basic examples below.

\begin{exm}\normalfont
\label{exm.fk}
Consider the case ${\cal V} ={\mathbb R}^d$ and write $a = x' (\neq x)$.
As usual we set $t =\ve^2$.
Then, the heat kernel $p_t (x,x')$ associated with the differential operator $L :=(1/2) \sum_{i=1}^r V_i^2 + V_0$
is expressed as
$$
p_{\ve^2} (x,x') = {\mathbb E} [\delta_{x'} (X^{\ve} (1,x,w) ) ].
$$
In this example $G(\ve, w) \equiv 1$.
(In our convention, 
$(e^{t L} f) (x) =\int_{{\mathbb R}^d} p_t (x,x') f(x' ) dx'$,
where $e^{t L}$ denotes the heat semigroup associated with $L$.)

Let $C$ be a tempered smooth function on ${\mathbb R}^d$
which is bounded from below.
Then, the heat kernel $p_t^C(x,x')$ associated with 
$L^C :=(1/2) \sum_{i=1}^r V_i^2 + V_0 - C$
is expressed as
\[
p_{\ve^2}^C  (x,x') 
= {\mathbb E} 
\Bigl[  
\exp \Bigl( - \ve^2 \int_0^1 C( X^{\ve} (s,x,w)) ds  \Bigr)
 \, \delta_{x'} (X^{\ve} (1,x,w) ) 
 \Bigr].
\]
In this case $G(\ve, w) =  \exp ( - \ve^2 \int_0^1 C( X^{\ve} (s,x,w)) ds )$
is the Feynman-Kac multiplicative functional
and satisfies {\bf (C1)}.
See pp. 414--415, Ikeda and Watanabe \cite{iwbk}.

Due to the Feynman-Kac-It\^o formula,
the heat kernels of magnetic Schr\"odinger operators  
also admit a similar expression
(We omit details.
See Section 6 in Ikeda's survey \cite{ike}.)
\end{exm}


\subsection{The manifold case}

Let ${\cal M}$ and ${\cal N}$ be compact manifolds 
with dimension $d$ and $n$ ($1 \le n \le d$), respectively,
and let $\Pi \colon {\cal M} \to {\cal N}$ be a smooth submersion.
${\cal N}$ is equipped with a smooth volume  ${\rm vol}$.
(A measure on ${\cal N}$ is said to be a smooth volume
if it is expressed on each coordinate chart as a strictly positive
smooth density function times the Lebesgue measure.)

For vector fields $V_i~(0 \le i \le r)$ on ${\cal M}$,
we study the following scaled SDE and its corresponding skeleton ODE:
\begin{eqnarray}
dX^{\ve}_t &=&
\ve \sum_{i=1}^r  V_i ( X^{\ve}_t) \circ  dw_t^i  + \ve^2  V_0 (X^{\ve}_t)   dt
\qquad
\qquad
\mbox{with \quad $X^{\ve}_0 =x \in {\cal M}$,}
\label{sde.mfd}
\\
d\phi_t
&=&
\sum_{i=1}^r  V_i ( \phi_t) dh_t^i 
\qquad
\qquad
\qquad
\qquad
\mbox{with \quad $\phi_0 =x \in {\cal M}$.}
\label{ode.mfd}
\end{eqnarray}
Set 
$Y^{\ve}_t = \Pi (X^{\ve}_t)$
and
$\psi_t = \Pi (\phi_t )$.
When necessary we will write $X^{\ve}_t = X^{\ve}(t, x, w)$ etc.
as before.
For $x \in {\cal M}$ and $a \in {\cal N}$ such that $\Pi(x) \neq a$,
we define 
${\cal K}_a$,  ${\cal K}^{min}_a$, and $d_a$
in the same way as in the Euclidean case.

Malliavin calculus for SDEs on manifolds was 
studied by Taniguchi \cite{ta}.
Roughly speaking, most of important results 
in the flat space case still hold true in the manifold case.


We impose two assumptions on the coefficient vector fields.
They are  similar to 
the corresponding ones in the Euclidean case.
\\
\\
{\bf (A1)'}:  $V_{i}$
is a  smooth vector fields on ${\cal M}$ $(0 \le i \le r)$.
\\
\\
Under this condition, 
$X^{\ve}_t$ is smooth in the sense of Malliavin calculus for any $(t, \ve) \in [0,1]^2$.
(For the definition of manifold-valued smooth Wiener functionals, 
see \cite{ta}.)

Now we introduce a H\"ormander-type condition.
$\Sigma_{k}$ is defined as in the Euclidean case.
Since the manifold ${\cal M}$ is compact, the partial H\"ormander condition 
below is automatically uniform as in (\ref{upHc}).
Therefore, {\bf (A2)'} and {\bf (A2)} are parallel.
\\
\\
{\bf (A2)'}: 
We say that (A2)' holds if either of the following two conditions holds:
\\
{\rm (i)}~
At the starting point $x \in {\cal M}$,  
$\cup_{k=1}^{\infty}  \{  W(x) \mid  W \in \Sigma_{k}\}$ spans 
$T_{x} {\cal M}$ in the sense of linear algebra.
\\
{\rm (ii)}~
For every $x \in {\cal M}$,  
$\cup_{k=1}^{\infty}  \{ (\Pi_*)_{x}  W(x) \mid  W \in \Sigma_{k}\}$ spans 
$T_{\Pi (x)} {\cal N}$ in the sense of linear algebra.
\\
\\
Choose a Riemannian metric on ${\cal N}$
so that the determinant of the (determinisitic) Malliavin covariance 
of ${\cal N}$-valued functionals are well-defined.
Results in
Taniguchi \cite{ta} or 
Section 5 in Kusuoka-Stroock \cite{ks2} apply to this case.
$Y_t^{\ve}$ is non-degenerate in 
the sense of Malliavin calculus for any $t>0$ and $\ve \in (0,1]$.
Moreover, Kusuoka-Stroock's moment estimate (\ref{ks_bnd.ineq}) 
for $\det \sigma [Y_1^{\ve} ]^{-1}$
also holds in this case.


It is known that integration by parts formula also
holds for manifold-valued Wiener functional. (See Section 8 for a proof.)
%
Therefore, as in the Euclidean case, if an ${\cal N}$-valued 
smooth Wiener functional $F$ is non-degenerate 
in the sense of Malliavin, 
then the composition $T(F) =T\circ F \in \tilde{{\mathbb D}}_{-\infty}$
is well-defined 
as a Watanabe distribution
for any distribution $T$ on ${\cal N}$.
In particular, 
$\delta_a ( Y^{\ve}_1)$
a positive Watanabe distribution under {\bf (A2)'} for every $a \in {\cal N}$.


On the other hands,
Assumptions {\bf (B1)}, {\bf (B2)}, {\bf (B3)}, {\bf (C1)}  and {\bf (C2)} 
 need not be modified 
and will be imposed in the manifold setting again.
In this case, too, 
non-degeneracy of $\sigma [\psi_1] (h)$ is equivalent to surjectivity of  
the tangent map $D \psi_1 (h) \colon {\cal H} \to T_a{\cal N}$.


Without loss of generality, we may assume that 
a Riemannian metric is given on both ${\cal M}$ and ${\cal N}$.
The reason we introduce them is as follows:
One on ${\cal M}$ is needed when a
tubular neighborhood on ${\cal M}$ is used.
Since Malliavin covariance matrix
of ${\cal N}$-valued functional is actually a bilinear form 
on the cotangent space, 
a metric on the cotangent space is needed 
when 
the determinant and the eigenvalues of the 
Malliavin covariance matrix are considered.
For our purpose any Riemannian metric will do.
In particular, even if we change
the Riemannian metric on ${\cal N}$,
{\bf (A2)'} and {\bf (B1)} (and the other assumptions) remain equivalent.
Hence, $ \delta_a (Y_1^{\ve})$
does not depend on the choice of the Riemannian metrics on 
${\cal M}$ and ${\cal N}$.

On the other hand, $\delta_a$ and 
$ \delta_a (Y_1^{\ve})$ depend on the choice of  ${\rm vol}$.
Since any other smooth volume can be expressed as
$\widehat{\rm vol} (dy) =\rho (y) {\rm vol} (dy)$
for some strictly positive smooth function $\rho$ on ${\cal N}$,
the delta function with respect to 
$\widehat{\rm vol}$ is given by $\hat{\delta}_a =\rho (a)^{-1} \delta_a$.
Therefore, it is sufficient to prove 
Theorem \ref{tm.main_mfd} below for one particular smooth volume.
In the proof of Theorem \ref{tm.main_mfd} in Section \ref{sec.pr_mfd}, 
we will assume that ${\rm vol}$ is the Riemannian measure on ${\cal N}$.

\begin{tm}\label{tm.main_mfd}
Let $x \in {\cal M}$ and $a \in {\cal N}$ such that $\Pi (x) \neq a$.
Assume {\bf (A1)'},  {\bf (A2)'}, {\bf (B1)}, {\bf (B2)}, {\bf (B3)} and {\bf (C1)}.
Then, we have the following asymptotic expansion:
\[
{\mathbb E} [G(\ve, w)  \delta_a (Y_1^{\ve}) ]
 \sim
e^{ - d_a^2 /2\ve^2}
 \ve^{-(n+ n')}
 (c_0 + c_1 \ve + c_2 \ve^2 + \cdots)
 \qquad
 \mbox{ as $\ve\searrow 0$}
\]
for certain constants $c_j \in {\mathbb R}~(j \ge 0)$.
If we assume {\bf (C2)} in addition, then $c_{2j +1} =0~(j \ge 0)$. 
\end{tm}

\begin{re}\normalfont
Compactness of ${\cal M}$ and 
${\cal N}$ in Theorem  \ref{tm.main_mfd} is assumed for simplicity.
Since manifold-valued Malliavin calculus in \cite{ta}
is developed for a reasonable class of non-compact manifolds,
we believe that Theorem  \ref{tm.main_mfd} extends to the case of 
non-compact manifolds under fairly mild assumptions.
%
(In this paper, however, only the compact case is proved.
The non-compact case is just a conjecture.)
\end{re}

\begin{re}\normalfont\label{re.sos}
A  remark similar to Example \ref{exm.fk} holds 
in the manifold setting, too.
In particular, if ${\cal M}= {\cal N}$ and $\Pi$ is the identity map,
then ${\mathbb E} [\delta_{x'} (X^{\ve} (1,x,w) ) ]$
is the Feynman-Kac representation of 
the heat kernel $p_{\ve^2} (x,x')$
associated with $L :=(1/2) \sum_{i=1}^r V_i^2 + V_0$.
\end{re}


However, some important second order differential operators
cannot be written in a  ``sum of squares" form as in Remark \ref{re.sos} above,
which is the reason why two manifolds are considered.
Here are examples of manifolds ${\cal M}$ and ${\cal N}$
we have in mind.

The first example is the orthonormal frame bundle over 
a compact Riemannian manifold.
This is a familiar example in stochastic analysis on manifolds 
and called Eells-Elworthy's construction 
of It\^o's stochastic parallel transport.
(See 
Section V-4, \cite{iwbk} or Hsu \cite{hsu}, Stroock \cite{st} for example.)
\begin{exm}\normalfont\label{exm.riem}
Let ${\cal N}$ be a compact Riemannian manifold with $\dim {\cal N} = n$
and 
let ${\cal M} ={\cal O} ( {\cal N})$ be its orthonormal frame bundle. 
Hence, $d := \dim {\cal M} =n(n+1)/2$.
Let $\Pi \colon {\cal O} ( {\cal N}) \to {\cal N}$ be the canonical projection.
We take $V_i$ to be the $i$th canonical horizontal vector field on ${\cal O} ( {\cal N})$
for $1 \le i \le n$
and 
set $V_0 \equiv 0$.

Let $\hat{a} \in {\cal N}$
and consider SDE (\ref{sde.mfd})
with $r=n$, $x \in \Pi^{-1} (\hat{a}) \subset {\cal O} ( {\cal N})$  and $a \in {\cal N}$.
Then, regardless of the choice of $x$,
 the process $t \mapsto Y^{\ve} (t, x,w)$
is the Brownian motion on ${\cal N}$, that is, 
the diffusion process associated to the (minus one half of) Laplace-Beltrami 
operator $\triangle_{{\cal N}}$ starting at $\hat{a} \in {\cal N}$.

In this case {\bf (A1)'} and {\bf (A2)'} are clearly satisfied
and $\delta_a (Y^{\ve} (1, x,w))$ is well-defined.
Moreover, ${\mathbb E} [\delta_a (Y^{\ve} (1, x,w)) ]
= p_{\ve^2} (\hat{a},a)$,
the heat kernel associated to $- \triangle_{{\cal N}} /2$.
Note that the deterministic Malliavin covariance is non-degenerate 
at any $h \in {\cal H}$ since the coefficient vector fields are elliptic at $x$.
So, {\bf (B1)} is also satisfied.

We denote by ${\cal H}_{\hat{a}} ({\cal N})$
the set of absolutely continuous path $\xi$ on ${\cal N}$
starting at $\hat{a}$
with ${\cal E} (\xi) := \int_0^1 \| \xi^{\prime}_s \|^2_{ T_{\xi_s} {\cal N}}  ds <\infty$.
The Cartan development map
 $h \mapsto \psi(h)$ is a diffeomorphism (of Hilbert manifolds)
 between ${\cal H}$ and ${\cal H}_{\hat{a}} ({\cal N})$
 which preserves the energy, that is, ${\cal E} (\psi (h)) =\|h\|^2_{{\cal H}}$
 for any $h \in {\cal H}$.
  Therefore, 
 \[ d_a^2 =\min \{ {\cal E} (\xi) \mid \xi \in  {\cal H}_{\hat{a}} ({\cal N}), \xi_1 =a \}
\mbox{ and } 
  \psi ({\cal K}_a^{min}) 
  = \{ \xi  \in  {\cal H}_{\hat{a}} ({\cal N}) \mid \xi_1 =a,  {\cal E} (\xi)=d_a^2 \}
 \]
  which will be denoted by ${\cal K}_{\hat{a}, a}^{min} ({\cal N})$.
  In other words, the set of minimal energy paths 
  does not change through $\psi$.
  Moreover, $d_a$ is equal to the Riemannian distance $d_{{\cal N}} (\hat{a},a)$.

  Let us rewrite the other assumptions in terms 
  of ${\cal K}_{\hat{a}, a}^{min} ({\cal N})$.
First,  {\bf (B2)} is equivalent to   
 \\
 \\
 {\bf (R2)}: ${\cal K}_{\hat{a}, a}^{min} ({\cal N})$ is a smooth
 and  compact manifold of finite dimension $n'$
regularly embedded in ${\cal H}_{\hat{a}} ({\cal N})$.
 \\
 \\
Set ${\cal K}_{\hat{a}, a} ({\cal N})
:= \{  \xi \in  {\cal H}_{\hat{a}} ({\cal N}) \mid \xi_1 =a \}$
and restrict ${\cal E}$ to this set.
Then,  {\bf (B3)} is equivalent to
\\
\\
{\bf (R3)}: 
For any $\xi \in {\cal K}_{\hat{a}, a}^{min} ({\cal N})$ and 
any $k \in T_{\xi}  {\cal K}_{\hat{a}, a} ({\cal N}) \setminus 
T_{\xi} {\cal K}_{\hat{a}, a}^{min} ({\cal N})$,
${\cal E}^{ \prime\prime} (\xi) \la k,k \ra >0$.
\\
\\
Then, as a corollary of our main theorem, we obtain the following:
Assume $\hat{a} \neq a$, {\bf (R2)} and {\bf (R3)}.
Then,  as a special case ($G^{\ve} \equiv 1$) 
of Theorem \ref{tm.main_mfd},
we have the following asymptotic expansion:
\begin{equation}\label{tm.Rmain}
p_{\ve^2} (\hat{a}, a)
 \sim
\exp \Bigl( -  \frac{ d_{{\cal N}} (\hat{a},a)^2 }{2\ve^2 }  \Bigr)
 \ve^{-(n+ n')}
 (c_0 + c_2 \ve^2 + c_4 \ve^4 + \cdots)
 \qquad
 \mbox{ as $\ve\searrow 0$}
\end{equation}
for certain constants $c_{2j} \in {\mathbb R}~(j \ge 0)$.
 \end{exm}

Conditions {\bf (R2)} and {\bf (R3)} are natural and 
often appear in 
Morse theory and Riemannian geometry
(possibly in a slightly different form).
Note that the dimension of the null eigenspace of ${\cal E}^{ \prime\prime} (\xi)$ 
equals the dimension of the set of Jacobi fields along $\xi$ 
which vanish at both ends.

The simplest example which satisfies {\bf (R2)} and {\bf (R3)} 
is the sphere ${\mathbb S}^{n}$ and two antipodal points $\hat{a}$ and $a$ on it.
In that case, ${\cal K}_{\hat{a}, a}^{min}$ consists of great circles from $\hat{a}$ to $a$
and therefore diffeomorphic to ${\mathbb S}^{n-1}$.


\begin{re}\normalfont \label{re.ludw}
We are not aware of any preceding work that proved (\ref{tm.Rmain}) above.
However, a few months after we finished our present work, 
Ludewig \cite{lu} proved a similar theorem.
His conditions on 
${\cal K}_{\hat{a}, a}^{min}$, which is called Morse-Bott condition in  \cite{lu},
 are essentially equivalent to {\bf (R2)} and {\bf (R3)}.
His method is basically analytic.
Compered to  (\ref{tm.Rmain}), his result is stronger 
in the following senses.
First, he obtained the leading constant $c_0$
in two explicit ways.
Second, he studied heat kernels associated with 
Laplace-type operators 
acting on sections of a vector bundle in a systematic way.
On the other hand, no hypoelliptic case is studied in \cite{lu}.
\end{re}


The next example is 
a compact, strictly pseudo-convex CR manifold
and its unitary frame bundle in Kondo and Taniguchi \cite {kt}.
For fundamental facts on CR manifolds, 
the reader is referred to \cite{dt} among others.

\begin{exm}\normalfont\label{exm.CR}
A CR (Cauchy-Riemann) manifold ${\cal N}$ is a real smooth 
manifold equipped with a complex subbundle $T_{1,0}$ of the
complexified tangent bundle ${\mathbb C}T{\cal N}$ with properties 
that $T_{1,0}\cap T_{0,1}=\{0\}$, where 
$T_{0,1}=\overline{T_{1,0}}$ is the complex conjugate of $T_{1,0}$ 
and 
$[T_{1,0},T_{1,0}]\subset T_{1,0}$.
Assume that ${\cal N}$ is compact, orientable and of real 
dimension $2k+1$ and $T_{1,0}$ is of complex dimension $k$.

There exists a real non-vanishing $1$-form $\theta$ on ${\cal N}$  
annihilating $H=\text{Re}(T_{1,0}\oplus T_{0,1})$.
We assume that ${\cal N}$ is strictly pseudo-convex, i.e. 
the Levi form $L_\theta$ defined by
\[
    L_\theta(Z,W)=-\sqrt{-1}\,d\theta(Z,W)
    \quad\text{for }Z,W\in\Gamma^\infty(T_{1,0}\oplus T_{0,1}),
\]
where $\Gamma^\infty(V)$ is the space of smooth  sections
of vector bundle $V$, is positive definite.
Associated with $\theta$, the characteristic direction $T$, the unique
real vector field on ${\cal N}$ transverse to $H$ is given by
\[
    T \rfloor d\theta=0, \quad T\rfloor\theta=1,
\]
where $T\rfloor \omega$ is the interior product.

The Webster metric $g_\theta$ is defined by 
\[
    g_\theta(X,Y)=d\theta(X,JY),\quad
    g_\theta(X,T)=0,\quad
    g_\theta(T,T)=1
    \quad\text{for }X,Y\in \Gamma^\infty(H),
\]
where $J \colon T{\cal N}\to T{\cal N}$ is defined so that its complex
linear extension to ${\mathbb C}T{\cal N}$ is equal to the 
multiplication by $\sqrt{-1}$ on $T_{1,0}$, that by $-\sqrt{-1}$
on $T_{0,1}$, and $J(T)=0$.
The Tanaka-Webster connection is a unique linear connection 
$\nabla$ on ${\cal N}$ such that
\begin{align*}
   & \nabla_X Y\in\Gamma^{\infty}(H),
      \quad \nabla J=0, \quad \nabla g_{\theta}=0, \quad
      T_{\nabla}(Z,W)=0, 
   \\
   & T_{\nabla}(Z,W^\prime)
       =2\sqrt{-1}  L_{\theta}(Z,W^\prime)T,\quad
       T_{\nabla}(T,J(X))+J(T_{\nabla}(T,X))=0
\end{align*}
for $X\in\Gamma^{\infty}(T{\cal N})$,
$Y\in\Gamma^{\infty}(H)$, $Z,W\in \Gamma^{\infty}(T_{1,0})$,
and $W^\prime\in \Gamma^{\infty}(T_{0,1})$, 
where $\nabla_X$ is the covariant derivative in the direction of 
$X$ and $T_{\nabla}$ is the torsion tensor field of $\nabla$.

Define the unitary frame bundle over ${\cal N}$ by
\[
    U(T_{1,0}) 
    =\coprod_{x\in M}\{u \colon {\mathbb C}^k  \to(T_{1,0})_x; 
         \text{$u$ is a unitary isometry}\}.
\]
For $u\in U(T_{1,0})$ with $u \colon {\mathbb C}^k  \to(T_{1,0})_x$, 
let $\Pi(u)=x$. 
Every $v\in T_x{\cal N}$ admits a unique horizontal lift
$\eta_u(v)\in T_u U(T_{1,0})$, where $u\in\Pi^{-1}(x)$,
so that there exist smooth curves $p \colon [-\tau ,\tau ]\to {\cal N}$ and
$\widehat{p} \colon [-\tau ,\tau]\to U(T_{1,0})$ such that
$\Pi\circ\widehat{p}=p$, 
$\widehat{p}(0)=u$, $\widehat{p}^\prime(0)=\eta_u(v)$,
$(\Pi_*)_u\eta_u(v)=v$, and
the curve $[-\tau ,\tau]\ni t\mapsto \widehat{p}(t)\zeta\in T_{1,0}$ is
a parallel section along $p$ with respect to $\nabla$  for any 
$\zeta\in{\mathbb C}^k$.
Extending $\eta_u$ naturally to a complex linear mapping of
${\mathbb C}T_xM$ to $T_u U(T_{1,0})$,  set
$L(\zeta)_u=\eta_u(u\zeta)$ for $\zeta\in{\mathbb C}^k$.
Using the standard complex basis $\{e_1,\dots,e_k\}$ of 
${\mathbb C}^k$, define the canonical vector fields $L_1,\dots,L_k$
by
\[
    L_i=L(e_i),\quad i=1,\dots,k.
\]
%
%
%

Let $V_1, \ldots, V_{2k}$ be vector fields on $U(T_{1,0})$ defined by
\[
V_i = \frac{1}{\sqrt{2}} ( L_i + \overline{L_i})
\quad
\mbox{and}
\quad
V_{k+i} = \frac{1}{\sqrt{2}\sqrt{-1} } ( L_i - \overline{L_i}),
\quad
i=1,\dots,k.
\]
Consider the stochastic differential equation on $U(T_{1,0})$ given by
\[
    d X_t^{\ve} =\ve \sum_{i=1}^{2k} V_i( X_t^{\ve} )\circ dw^i(t)
\]
for $0 < \ve \le 1$.
Then $Y_t^{\ve}:= \Pi(X_t^{\ve})$ determines a diffusion process generated by 
$- (\ve^2/2) \triangle_b$.
Here the operator $\triangle_b$ is a hypoelliptic sub-Laplacian given by
\[
    \int_{\cal N}  (\triangle_b f) g  \, d{\rm vol}
    =\int_{\cal N} L_\theta^*(d_bf,d_b g) \,  d{\rm vol}
    \]
for any smooth $f,g$ on ${\cal N}$, where
${\rm vol}=\theta\wedge(d\theta)^k$, $L_\theta^*$ is the dual
metric on $H^*$ of $L_\theta$, and  $d_bf$ is the projection of $df$
onto $H^*$.
%

Assumption {\bf (A1)'} clearly holds. Moreover 
{\bf (A2)'} is also satisfied as was seen in \cite{kt}.
Indeed, if $\{Z_i\}_{i=1}^n$ is a local orthonormal frame of 
$T_{1,0}$, then 
\[ 
    (\Pi_*)_u(L_i)=\sum_{j=1}^k  e_j^i Z_j
    \quad\text{and}\quad
    (\Pi_*)_u[L_i,L_{\overline{i}}] =-2\sqrt{-1}T
    \quad\text{mod } \{Z_i,Z_{\overline{i}} \mid i=1,\dots,k \},
\]
where $(e_i^j)\in U(n)$.
See Kondo and Taniguchi \cite{kt} for details.

We say that an absolutely continuous path $\xi$ on ${\cal N}$
is {\it horizontal in the CR sense} if $\xi^{\prime}_t \in H_{\xi_t}$
for almost all $t$.
We denote by 
$\tilde{\cal H}_{\hat{a}} ({\cal N})$
the set of absolutely continuous path $\xi$ on ${\cal N}$
which is horizontal in the CR sense and start at $\hat{a}$.
The energy functional $\tilde{\cal E}$ on $\tilde{\cal H}_{\hat{a}} ({\cal N})$
is defined 
in the same way as in Example \ref{exm.riem}, this time with 
the Webster metric instead of the Riemannian metric. 
In this case again, 
the  development map
 $h \mapsto \psi(h)$ is an energy-preserving 
 diffeomorphism (of Hilbert manifolds)
 between ${\cal H}$ and $\tilde{\cal H}_{\hat{a}} ({\cal N})$.
 Therefore, 
 \[ 
 d_a^2 
 =\min \{ \tilde{\cal E} (\xi) \mid \xi \in  \tilde{\cal H}_{\hat{a}} ({\cal N}), \xi_1 =a \}
\mbox{ and } 
  \psi ({\cal K}_a^{min}) 
  = \{ \xi  \in  \tilde{\cal H}_{\hat{a}} ({\cal N}) \mid \xi_1 =a,  \tilde{\cal E} (\xi)=d_a^2 \}
 \]
  which will be denoted by $\tilde{\cal K}_{\hat{a}, a}^{min} ({\cal N})$.
  In other words, the set of minimal energy paths 
  does not change through $\psi$.
  Moreover, $d_a$ is equal to the sub-Riemannian
  (Carnot-Carath\'eodory) distance $\tilde{d}_{{\cal N}} (\hat{a},a)$.
In a strictly pseudo-convex CR manifold, 
no non-trivial energy-minimzing path (i.e. geodesics) is
abnormal in the sense of sub-Riemannian geometry.
(See pp. 24--25 \cite{ri} for example.)
Hence, {\bf (B1)} is always satisfied.
%

Let us rewrite the other assumptions in terms 
  of $\tilde{\cal K}_{\hat{a}, a}^{min} ({\cal N})$.
First,  {\bf (B2)} is equivalent to   
 \\
 \\
 {\bf (CR2)}: $\tilde{\cal K}_{\hat{a}, a}^{min} ({\cal N})$ is a smooth
 and  compact manifold of finite dimension $n'$
regularly embedded in $\tilde{\cal H}_{\hat{a}} ({\cal N})$.
 \\
 \\
Set $\tilde{\cal K}_{\hat{a}, a} ({\cal N})
:= \{  \xi \in  \tilde{\cal H}_{\hat{a}} ({\cal N}) \mid \xi_1 =a \}$ 
and restrict $\tilde{\cal E}$ to this set.
Then,  {\bf (B3)} is equivalent to
\\
\\
{\bf (CR3)}: 
For any $\xi \in \tilde{\cal K}_{\hat{a}, a}^{min} ({\cal N})$ and 
any $l \in T_{\xi}  \tilde{\cal K}_{\hat{a}, a} ({\cal N}) \setminus 
T_{\xi} \tilde{\cal K}_{\hat{a}, a}^{min} ({\cal N})$,
$\tilde{\cal E}^{ \prime\prime} (\xi) \la l,l \ra >0$.
\\
\\
Then, as a corollary of our main theorem, we obtain the following:
Assume $\hat{a} \neq a$, {\bf (CR2)} and {\bf (CR3)}.
Then,  as a special case ($G^{\ve} \equiv 1$) 
of Theorem \ref{tm.main_mfd},
we have the following asymptotic expansion:
\begin{equation}\label{tm.CRmain}
\tilde{p}_{\ve^2} (\hat{a}, a)
 \sim
\exp \Bigl( -  \frac{ \tilde{d}_{{\cal N}} (\hat{a},a)^2 }{2\ve^2 }  \Bigr)
 \ve^{-(2k+1+ n')}
 (c_0 + c_2 \ve^2 + c_4 \ve^4 + \cdots)
 \qquad
 \mbox{ as $\ve\searrow 0$}
\end{equation}
for certain constants $c_{2j} \in {\mathbb R}~(j \ge 0)$.
Here, $\tilde{p}$ stands for the heat kernel 
associated with the sub-Laplacian $\triangle_b/2$.
\end{exm}

The next example is a continuation of Example \ref{exm.CR} above.
We provide a very concrete example of ${\cal N}$ and $\hat{a},a \in {\cal N}$
which satisfy {\bf (CR2)} and {\bf (CR3)}.

\begin{exm}\normalfont\label{exm.CR}
In this example, we show that antipodal points on the standard
CR sphere  ${\mathbb S}^{2k+1}$  ($k \ge 1$) satisfy 
{\bf (CR2)} and {\bf (CR3)}.
For basics of  the standard CR sphere, we refer to \cite{bw, cmv, mm}.
To keep notations simple, 
we avoid the complex coordinates and
describe the standard CR sphere as a real manifold.

Let ${\mathbb S}^{2k+1}=\{  (x_0,y_0, \ldots, x_k, y_k) \in {\mathbb R}^{2k+2}
\mid  \sum_{i=0}^k (x_i^2 +y_i^2) =1\}$.
It has a natural strictly pseudo-convex CR structure inherited from 
${\mathbb C}^{k+1} \cong {\mathbb R}^{2k+2}$, which is given as follows:
\begin{eqnarray*}
(T_{1,0})_{ (x_0,y_0, \ldots, x_k, y_k)} 
&=&
\Bigl\{
\sum_{i=0}^k \Bigl[
\Bigl(  
p_i \frac{\partial}{\partial x_i}   + q_i  \frac{\partial}{\partial y_i}  
\Bigr)
+
\sqrt{-1}
\Bigl(
q_i \frac{\partial}{\partial x_i}   -p _i  \frac{\partial}{\partial y_i}  
\Bigr)
\Bigr]
\quad
\Big|
\\
&&
(p_0,q_0, \ldots, p_k, q_k) \in {\mathbb R}^{2k+2},
\,
\sum_{i=0}^k  (p_i x_i + q_i y_i)=0,
\,
\sum_{i=0}^k  (q_i x_i -p_i y_i)=0
\Bigr\}.
\end{eqnarray*}
The horizontal subbundle $H$ coincides with the kernel 
of the following contact one-form:
$
\omega = \sum_{i=0}^k  (-y_i dx_i + x_i dy_i).
$
The Webster metric on ${\mathbb S}^{2k+1}$
 is just restriction of the canonical one on ${\mathbb R}^{2k+2}$.
%

%
We choose $(\pm 1, 0,  \ldots, 0,0)$ 
as the starting point and the end point, respectively, 
and denote them by ${\bf e}$ and $-{\bf e}$. 
By the way the metric is defined, we have
 ${\cal E} (\xi)  = \tilde{\cal E} (\xi)$ for a horizontal path $\xi$ on 
 ${\mathbb S}^{2k+1}$.
 In particular, we have 
 \[
 \tilde{\cal K}_{ {\bf e}, - {\bf e} }^{min} ({\mathbb S}^{2k+1} ) 
 =
 \{  \xi \in {\cal K}_{ {\bf e}, -{\bf e}}^{min} ({\mathbb S}^{2k+1} )
 \mid  
  \mbox{$\xi$ is horizontal in the CR sense}\}.
   \]
As we will see, the right hand side above is not empty.
The sub-Riemannian 
distance between the two points is $\pi$.
%

Recall that ${\cal K}_{ {\bf e}, -{\bf e}}^{min} ({\mathbb S}^{2k+1} )$
is the set of all the great circles from $ {\bf e}_+$ to $ {\bf e}_-$ (in time $1$).
We will write them down explicitly.
Set 
\begin{eqnarray*}
{\mathbb S}^{2k} &=&  
\{(p_0,q_0, \ldots, p_k, q_k) \in {\mathbb S}^{2k+1} \mid
p_0 =0\},
\nn\\
{\mathbb S}^{2k-1} &=&  
\{(p_0,q_0, \ldots, p_k, q_k) \in {\mathbb S}^{2k+1} \mid
p_0=q_0 =0\}.
\end{eqnarray*} 
For ${\bf v} = (0, q_0, \ldots, p_k, q_k)\in {\mathbb S}^{2k}$, 
we write 
\[
\xi_{ {\bf v}} (t) = \cos (\pi t) {\bf e}   + \sin (\pi t) {\bf v}.
\]
Then, $\xi_{ {\bf v}} $ is a great circle and 
we have 
\[
{\cal K}_{ {\bf e}, -{\bf e}}^{min} ({\mathbb S}^{2k+1} )
=\{ \xi_{ {\bf v}} \mid {\bf v} \in   {\mathbb S}^{2k}\}.
\]
Moreover, by straightforward computation, we have 
$
\langle \xi_{ {\bf v}}^{\prime},  \omega_{\xi_{{\bf v}}}   \rangle \equiv \pi q_0
$.
Hence, $\xi_{ {\bf v}}$ is horizontal in the CR sense 
if and only if $q_0 =0$
and we have 
\[
\tilde{\cal K}_{ {\bf e}, -{\bf e}}^{min} ({\mathbb S}^{2k+1} )
=\{ \xi_{ {\bf v}} \mid {\bf v} \in   {\mathbb S}^{2k-1}\}.
\]
Thus, we have shown {\bf (CR2)} with $n' =2k-1$.
%

Now we turn to {\bf (CR3)}.  
A key fact is 
\begin{equation}\label{eq.kyex}
T_{\xi} \tilde{\cal K}_{ {\bf e}, -{\bf e}}^{min} ({\mathbb S}^{2k+1} )
=
T_{\xi}  {\cal K}_{ {\bf e}, -{\bf e}}^{min} ({\mathbb S}^{2k+1} )
\cap 
T_{\xi} \tilde{\cal K}_{ {\bf e}, -{\bf e}} ({\mathbb S}^{2k+1} )
\qquad
(\xi \in \tilde{\cal K}_{ {\bf e}, -{\bf e}}^{min} ({\mathbb S}^{2k+1} )).
\end{equation}
Once (\ref{eq.kyex}) is established, the problem reduces 
to the corresponding one for ${\mathbb S}^{2k+1}$
as a Riemannian manifold.
The reason is as follows.
If  $l \in T_{\xi} \tilde{\cal K}_{ {\bf e}, -{\bf e}} ({\mathbb S}^{2k+1} )
 \setminus 
T_{\xi} \tilde{\cal K}_{ {\bf e}, -{\bf e}}^{min} ({\mathbb S}^{2k+1} )$,
then 
$l \in T_{\xi} {\cal K}_{ {\bf e}, -{\bf e}} ({\mathbb S}^{2k+1} )
 \setminus 
T_{\xi} {\cal K}_{ {\bf e}, -{\bf e}}^{min} ({\mathbb S}^{2k+1} )$ by  (\ref{eq.kyex}).
We have already seen 
in the paragraph just above Remark \ref{re.ludw} 
that
$\tilde{\cal E}^{ \prime\prime} (\xi) \la l,l \ra
={\cal E}^{ \prime\prime} (\xi) \la l,l \ra >0$ holds
for such $l$. 
This proves  {\bf (CR3)}.

It is obvious that the left hand side of (\ref{eq.kyex})
is included in the right hand side.
To see the converse inclusion, we use  
$F \colon {\cal K}_{ {\bf e}, -{\bf e}} ({\mathbb S}^{2k+1} ) \to {\mathbb R}$ 
defined by
$F(\xi) = \int_0^1 
 \langle \xi^{\prime} (t),  \omega_{\xi (t)}   \rangle  dt$.
Let $(-1,1) \ni \tau \mapsto  c(\tau) \in 
{\cal K}_{ {\bf e}, -{\bf e}} ({\mathbb S}^{2k+1} )$  be a smooth curve
such that $c(0) =\xi \in  \tilde{\cal K}_{ {\bf e}, -{\bf e}}^{min} ({\mathbb S}^{2k+1} )$.
Then, $\tau \mapsto  c(\tau)^{\prime}$ is a smooth curve in 
$L^2 ([0,1], {\mathbb R}^{2k+2})$
since the differentiation in $t$ (denoted by ``prime")
is a unitary isomorphism from 
the  (${\mathbb R}^{2k+2}$-valued) Cameron-Martin space to 
$L^2 ([0,1], {\mathbb R}^{2k+2})$.
It is easy to see that $(d/d\tau) \vert_{\tau =0} F(c(\tau))$ exists.
Moreover, it
depends only on $\xi =c(0)$ and $l:= (d/d\tau) \vert_{\tau =0} c(\tau)$.
(The latter can be regarded as a Cameron-Martin path  
in ${\mathbb R}^{2k+2}$.)
It is obvious that $(d/d\tau) \vert_{\tau =0} F(c(\tau)) =0$ 
if $l \in T_{\xi} \tilde{\cal K}_{ {\bf e}, -{\bf e}}
({\mathbb S}^{2k+1} )$.

On the other hand, we can show that 
$(d/d\tau) \vert_{\tau =0} F(c(\tau)) \neq 0$
if 
$l \in T_{\xi}  {\cal K}_{ {\bf e}, -{\bf e}}^{min} ({\mathbb S}^{2k+1} )
\setminus 
T_{\xi}  \tilde{\cal K}_{ {\bf e}, -{\bf e}}^{min} ({\mathbb S}^{2k+1} )$ as follows.
Let ${\bf v} = (0, 0, p_1,q_1, \ldots, p_k, q_k) \in {\mathbb S}^{2k-1}$
be such that 
 $\xi = \xi_{{\bf v}}$ and let
$\tau \mapsto {\bf v} (\tau) = (0, q_0 (\tau),  \ldots, p_k  (\tau), q_k  (\tau))$ 
be a  smooth curve in ${\mathbb S}^{2k}$ such that 
${\bf v} (0) = {\bf v}$.
It is easy to see that $(d/d\tau) \vert_{\tau =0} {\bf v}(\tau)$
is tangent to $ {\mathbb S}^{2k-1}$ at ${\bf v}$ 
if and only if $(d/d\tau) \vert_{\tau =0} q_0 (\tau)=0$.
As we have seen, 
$(d/d\tau) \vert_{\tau =0} F(\xi_{ {\bf v} (\tau)} ) =\pi (d/d\tau) \vert_{\tau =0} 
q_0 (\tau)$.
This proves our claim in this paragraph and hence (\ref{eq.kyex}), too.
%

%
Finally, we show that the sub-Laplacian $\triangle_b$ on 
${\mathbb S}^{2k+1}$
cannot be written in a ``sum of squares" form 
of exactly $2k$ vector fields if $2k+1 \neq 3$ or $7$,
that is, 
we cannot find vector fields $A_i ~(0 \le i \le 2k)$ 
on ${\mathbb S}^{2k+1}$ such that $\triangle_b = \sum_{i=1}^{2k} A_i^2 +A_0$.
(We should recall that $2k$ is the real dimension of $T_{1,0}$.)
Consequently, we need the stochastic parallel transport 
as in \cite{kt}
to construct 
the diffusion process generated by $-\triangle_b/2$
(as long as we use $2k$-dimensional Brownian motion).
The reason why $\triangle_b$ cannot be written in such a form is as follows.
If it could,
by a formula $\triangle_{LB} = \triangle_b + T^2$ in p.137  \cite{bw},
the Laplace-Beltrami operator
$\triangle_{LB}$ of  ${\mathbb S}^{2k+1}$ as a Riemannian manifold
would have a  ``sum of squares" form
$\triangle_{LB} = \sum_{i=1}^{2k} A_i^2 + T^2+A_0$.
However,  this is impossible if $2k+1 \neq 3, 7$.
(${\mathbb S}^{j}$ is not parallelizable for $j \neq1, 3, 7$
and therefore 
$\{A_1, \ldots, A_{2k},T \}$ cannot be a global frame of the tangent bundle of
 ${\mathbb S}^{2k+1}$.)
\end{exm}


\vspace{10mm}

In Proposition 3.7 \cite{bw}, the leading term of  
short time asymptotic of $\tilde{p}_t ({\bf e}, {\bf f})$
is calculated with an explicit value of the leading constant $c_0$,
where ${\bf f} =(e^{i \theta}, 0, \ldots , 0) 
\in {\mathbb C}^{n+1} \cong {\mathbb R}^{2n+2}$ ($0 < \theta <\pi$).
However, the case $ \theta =\pi$ (i.e. ${\bf f} =-{\bf e}$) seems to be excluded.
Since all the energy-minimizers on the standard CR sphere 
is obtained in \cite{cmv, mm}, 
it may be an interesting future task to strengthen Proposition 3.7 \cite{bw}
by combining our result (\ref{tm.CRmain}) with \cite{cmv, mm}.


Grong and Thalmaier \cite{gt} recently showed that 
diffusion processes on sub-Riemannian manifolds
associated with
sub-Laplacians also admit a similar construction.
For basic information on sub-Riemannian geometry, 
we refer to \cite{mo, ri, cc} among others.

\begin{exm}\normalfont\label{exm.subR}
(Diffusions associated with sub-Laplacians on sub-Riemannian manifolds.
See Section 2, \cite{gt} for details.)
Let $({\cal N}, {\cal D}, {\bf g})$ be a compact sub-Riemannian manifold.
Here,
${\cal D}$ is a sub-bundle of the tangent bundle $T{\cal N}$
of a  smooth compact manifold 
${\cal N}$ with $\dim {\cal N} =n$
and ${\bf g}$ is a metric tensor on ${\cal D}$.
By definition, ${\cal D}$ is bracket generating.

Let $\bar{\bf g}$ be a Riemannian metric tensor of ${\cal N}$
that tames ${\bf g}$ (i.e. $\bar{\bf g} |_{{\cal D}} = {\bf g}$).
Let $\Gamma$
 be the orthogonal complement of ${\cal D}$ in $T{\cal N}$.
We write the orthogonal projection by 
${\rm pr}_{{\cal D}}$ and ${\rm pr}_{\Gamma}$, respectively.
Denote by $\bar\nabla$  the Levi-Civita connection on ${\cal N}$
with respect to $\bar{\bf g}$.
Define a sub-Laplacian by using a local orthonormal frame 
$\{ A_1, \ldots, A_r\}$ of ${\cal D}$ by
\[
\triangle_{sub} = \sum_{i=1}^r A_i^2 + \sum_{i,j =1}^r  
{\bf g} \la  {\rm pr}_{{\cal D}} \bar\nabla_{A_i} A_j   , A_i  \ra   A_i
\]
where $r~(1<r<n)$ is the rank of ${\cal D}$. 
This is a globally well-defined differential operator 
and known to be hypoelliptic.

The diffusion process on ${\cal N}$ associated with 
$\triangle_{sub} /2 + A_0$
admits Eell-Elworthy's construction
for any vector field $A_0$ on ${\cal N}$
($A_0$ need not be a section of ${\cal D}$.)
In this case, the principle bundle is 
\begin{eqnarray*}
O({\cal D})\odot O(\Gamma)
&=&
\bigcup_{y \in {\cal N}}
\bigl\{ (u, v) \mid \mbox{ $u \colon {\mathbb R}^r \to {\cal D}_y$ and 
 $v\colon  {\mathbb R}^{n-r} \to \Gamma_y$ are isometries}
\bigr\}
\end{eqnarray*}
with its structure group $O(r) \times O(n-r)$.
The projection is denoted by $\Pi$.
Since $\bar\nabla$  is a metric connection, 
it defines a Cartan-Ehresmann connection on $O({\cal D})\odot O(\Gamma)$.
Define canonical horizontal vector fields
$\bar{A_i} ~(1\le i\le r)$ on 
$O({\cal D})\odot O(\Gamma)$
so that
$(\bar{A_i})_{(u,v)}$ is the horizontal lift of $u e_i \in {\cal D}_{ \Pi (u,v)}$,
where $\{ e_i\}_{i=1}^r$ is the canonical basis of ${\mathbb R}^r$.
Also define 
$\bar{A_0}$ to be the horizontal lift of $A_0$.
Then,  the solution to 
SDE on $O({\cal D})\odot O(\Gamma)$
with the coefficient vector fields $\bar{A}_i ~(0 \le i \le r)$
is the diffusion process associated with $\triangle_{sub} /2 + A_0$.
(In \cite{gt}, only the case $A_0 \equiv 0$ is treated.
However, a modification of this first order term is easy.)
\end{exm}

In Example \ref{exm.subR} above, 
an example of a submersion
$\Pi \colon {\cal M} \to {\cal N}$ and vector fields 
which satisfy {\bf (A1)} and {\bf (A2)} is given.
However, as readers may have noticed, there is no concrete exmaple 
of a sub-Riemannian manifold  and two points on it.
(In this sense this should be called a potential example.)
It is an interesting and important future task to find such examples 
which satisfies our assumptions {\bf (B1)}-- {\bf (B3)}.

\begin{re}\normalfont
Our present paper is based on Takanobu-Watanabe \cite{tw}
and the assumptions are quite similar.
In particular, the assumptions on the set of energy-minimizers
are the exactly same.
Our main result is stronger than the one in \cite{tw} because 
we work under the partial H\"ormander condition and 
work also on manifold. 

%
On  the other hand, Barilari-Boscain-Neel \cite{bbn} 
imposes assumptions on the set of energy-minimizers 
which look quite different from our {\bf (B1)}--{\bf (B3)}
and obtained the leading term of the asymptotics.

The relation between the two types of assumptions are unclear, yet.
Neither is it clear at the moment whether
a full asymptotic expansion can be proved under the assumptions in  \cite{bbn}.
\end{re}


\section{Preliminaries }

In this section we summarize  results we will use in the proof
 of our main theorems.  
All the results in this section are either known or easily 
derived from known results.

\subsection{Preliminaries from Malliavin calculus}
\label{subsec.31}

We first recall Watanabe's theory of 
generalized Wiener functionals (i.e. Watanabe distributions) in Malliavin calculus.
Most of the contents and the notations
in this subsection are contained in Sections V.8--V.10, Ikeda and Watanabe \cite{iwbk}
with trivial modifications.
We also refer to
Shigekawa \cite{sh}, Nualart \cite{nu},  Hu \cite{hu}
and Matsumoto and Taniguchi \cite{mt}. 
For basic results of quasi-sure analysis, see Chapter II, Malliavin \cite{ma}.

Let $({\cal W}, {\cal H}, \mu)$ be the classical Wiener space as before.  
(The results in this subsection also hold on any abstract Wiener space, however.)
%
We denote by $D$  the gradient operator (${\cal H}$-derivative) 
and by $L = -D^* D$ the Ornstein-Uhlenbeck operator.
%
The following are of particular importance in this paper:
\\
\\
{\bf (a)}~ Basics of Sobolev spaces:
We denote by ${\mathbb D}_{p,r} ({\cal X})$ 
 the Sobolev space of ${\cal X}$-valued 
(generalized) Wiener functionals, 
where $p \in (1, \infty)$, $r \in {\mathbb R}$, and ${\cal X}$ is a real separable Hilbert space.
As usual, we will use the spaces 
${\mathbb D}_{\infty} ({\cal X})= \cap_{k=1 }^{\infty} \cap_{1<p<\infty} {\mathbb D}_{p,k} ({\cal X})$, 
$\tilde{{\mathbb D}}_{\infty} ({\cal X}) 
= \cap_{k=1 }^{\infty} \cup_{1<p<\infty}  {\mathbb D}_{p,k} ({\cal X})$ of test functionals 
and  the spaces ${\mathbb D}_{-\infty} ({\cal X}) = \cup_{k=1 }^{\infty} \cup_{1<p<\infty} {\mathbb D}_{p,-k} ({\cal X})$, 
$\tilde{{\mathbb D}}_{-\infty} ({\cal X}) = \cup_{k=1 }^{\infty} \cap_{1<p<\infty} {\mathbb D}_{p,-k} ({\cal X})$ of 
 Watanabe distributions as in \cite{iwbk}.
When ${\cal X} ={\mathbb R}$, we simply write ${\mathbb D}_{p, r}$, etc.
The ${\mathbb D}_{p, r} ({\cal X})$-norm is denoted by $\| \,\cdot\, \|_{p,r}$.
The precise definition of an asymptotic expansion 
up to any order can be found in Section V-9, \cite{iwbk}.
\\

{\bf (b)}~ Meyer's equivalence of Sobolev norms:
See Theorem 8.4, \cite{iwbk}. 
A stronger version can be found in Theorem 4.6 in \cite{sh},  
Theorem 1.5.1 in \cite{nu} or Theorem 5.7.1 in Bogachev \cite{bog}.
\\

{\bf (c)}~Watanabe's pullback: 
Pullback 
$T \circ F =T(F)\in \tilde{\mathbb D}_{-\infty}$ 
of a tempered Schwartz distribution $T \in {\cal S}^{\prime}({\mathbb R}^n)$
on ${\mathbb R}^n$
by a non-degenerate Wiener functional $F \in {\mathbb D}_{\infty} ({\mathbb R}^n)$. 
(See Sections 5.9 \cite{iwbk}.)
The key to prove this pullback is an integration by parts formula 
in the sense of Mallavin calculus.
(Its generalization is given in Item {\bf (d)} below.)
\\

{\bf (d)}~A generalized version of the integration by parts formula in the sense 
of Malliavin calculus
 for Watanabe distribution,
which is given as follows (see p. 377,  \cite{iwbk}):

For a non-degenerate Wiener functional
$F =(F^1, \ldots, F^n) \in {\mathbb D}_{\infty} ({\mathbb R}^n)$, 
we denote by $\sigma [F ](w) = \sigma_F (w) $
the Malliavin covariance matrix of $F$ 
whose $(i,j)$-component is given by 
$\sigma_F^{ij} (w) =  \la DF^i (w),DF^j (w)\ra_{{\cal H}}$.
We denote by $\gamma^{ij}_F (w)$ the $(i,j)$-component of the inverse matrix $\sigma^{-1}_F$.
Note that $\sigma^{ij}_F \in {\mathbb D}_{\infty} $ and
$D \gamma^{ij}_F =- \sum_{k,l} \gamma^{ik}_F ( D\sigma^{kl}_F ) \gamma^{lj}_F $.
Hence, derivatives of $\gamma^{ij}_F$ can be written in terms of
$\gamma^{ij}_F$'s and the derivatives of $\sigma^{ij}_F$'s.
Suppose 
$G \in {\mathbb D}_{\infty}$ and $T \in {\cal S}^{\prime} ({\mathbb R}^n)$.
Then, the following integration by parts holds;
\begin{align}
{\mathbb E} \bigl[
(\partial_i T \circ F )  G 
\bigr]
=
{\mathbb E} \bigl[
(T \circ F )  \Phi_i (\, \cdot\, ;G)
\bigr],
\label{ipb1.eq}
\end{align}
where $\Phi_i (w ;G) \in  {\mathbb D}_{\infty}$ is given by 
\begin{align}
\Phi_i (w ;G) &=
\sum_{j=1}^d  D^* \Bigl(   \gamma^{ij }_F(w)  G (w) DF^j(w)
\Bigr)
\nn\\
&=
-
\sum_{j=1}^d  
\Bigl\{
-\sum_{k,l =1}^d G(w) \gamma^{ik }_F (w)\gamma^{jl }_F (w)    \la D\sigma^{kl}_F (w),DF^j (w)\ra_{{\cal H}}
\nn\\
&
\qquad\qquad
+
\gamma^{ij }_F (w) \la DG (w),DF^j (w)\ra_{{\cal H}} + \gamma^{ij }_F (w) G (w) LF^j (w)
\Bigr\}.
\label{ipb2.eq}
\end{align}
Note that the expectations in (\ref{ipb1.eq}) are in fact  the generalized ones,
i.e.
the pairing of $\tilde{{\mathbb D}}_{- \infty}$ and $\tilde{{\mathbb D}}_{\infty}$.
\\

%
Watanabe's asymptotic expansion theorem is a key theorem in his distributional Malliavin calculus.
Its standard version 
can be found  in Theorem 9.4, pp. 387-388, \cite{iwbk} or Watanabe \cite{wa}.
In the present paper, however, 
we need a modified version in pp. 216--217, \cite{tw}.
Though it plays a key role, no proof is given in \cite{tw} unfortunately.
Therefore, we will prove it below.

Let $\rho >0$, $\xi \in {\mathbb D}_{\infty}$ 
and $F \in {\mathbb D}_{\infty} ({\mathbb R}^n)$
and suppose that
\begin{equation}\label{wat1.eq}
\inf_{v \in {\mathbb R}^{n} \mid  |v|=1}    v^* \sigma_F (w) v  \ge \rho
\qquad
\mbox{on  \quad$\{w \in {\cal W} \mid |\xi (w)| \le 2 \}$.}
\end{equation}
Let $\chi \colon {\mathbb R} \to {\mathbb R}$ be a smooth function 
whose support is contained in $[-1, 1]$.
Then, the following proposition holds (Proposition 6.1, \cite{tw}).

\begin{pr}\label{pr.comp1}
Assume (\ref{wat1.eq}).
For every $T \in {\cal S}^{\prime}({\mathbb R}^n)$, 
$\chi(\xi) (T \circ F) =\chi(\xi) T (F) \in  \tilde{\mathbb D}_{-\infty}$
can be defined in a unique way so that the following properties hold:
\\
\noindent
{\rm (i)}~If $T_k \to T \in {\cal S}^{\prime}({\mathbb R}^n)$ as $k \to \infty$, 
then $\chi(\xi) T_k (F) \to \chi(\xi)  T (F) \in \tilde{\mathbb D}_{-\infty}$.
\\
\noindent
{\rm (ii)}~If $T$ is given by $g \in {\cal S}({\mathbb R}^n)$,
then $\chi(\xi) T (F) = \chi(\xi)  g (F) \in {\mathbb D}_{\infty}$.
\end{pr}

\Proof
Let $\eta \colon {\mathbb R}\to{\mathbb R}$ be a smooth function
whose support is contained in $(-2,2)$.
By the assumption (\ref{wat1.eq}), 
$\sigma_F$ is invertible on $\{\eta(\xi)\ne0\}$
and $\eta(\xi)\sigma_F^{-1}$ is of class ${\mathbb D}_{\infty}$.
In fact, it is the limit of $\eta(\xi)(\sigma_F+\frac1{m} {\rm Id}_n)^{-1}$ in
${\mathbb D}_{\infty}$ as $m\to\infty$, where ${\rm Id}_n$ stands for the 
$n\times n$-identity matrix.
Hence $\Phi_i(\,\cdot\,  ;\eta(\xi)G)$ in (\ref{ipb2.eq}) is well defined
and of class ${\mathbb D}_{\infty}$.
Moreover, observe that
\begin{equation}\label{tn.eq0}
    \Phi_i(\,\cdot\,  ;\eta(\xi)G)
    =D^*\biggl(\eta(\xi) G\sum_{k=1}^n \gamma_F^{ik} 
     DF^k\biggr)
\end{equation}
for any 
$G\in{\mathbb D}_{\infty}$.
Let $f\in{\mathcal S}({\mathbb R}^n)$. Since  
\[
    \eta(\xi)(\partial_i f\circ F)
    =\eta(\xi)\sum_{k=1}^n 
    \langle D(f\circ F),DF^k\rangle_{{\mathcal H}} \gamma_F^{ik},
\]
we have the same integration by parts formula as (\ref{ipb1.eq});
\begin{equation}\label{tn.eq1}
    {\mathbb E}[\eta(\xi)(\partial_i f\circ F)G]
    ={\mathbb E}[(f\circ F)\Phi_i(\,\cdot\,  ;\eta(\xi)G)].
\end{equation}
Furthermore, for $p>1$ and $r>0$, let $C_{p,r}$ be a constant 
such that
\[
    \|D^*(GK)\|_{p,r}
    \le C_{p,r}\|G\|_{2p,r+1}
       \|K\|_{2p,r+1}
\]
for any $G\in{\mathbb D}_{2p,r+1}$ and 
$K\in{\mathbb D}_{2p,r+1}({\mathcal H})$.
Then 
\begin{equation}\label{tn.eq2}
    \|\Phi_i(\cdot;\eta(\xi)G)\|_{p,r}
    \le C_{p,r} \biggl\|\eta(\xi) \sum_{k=1}^n \gamma_F^{ik} 
         DF^k\biggr\|_{2p,r+1}
         \|G\|_{2p,r+1}.
\end{equation}

Take a sequence $\{\chi_k\}_{k=1}^\infty$ of smooth functions on 
${\mathbb R}$ such that
$\chi_1=\chi$, $\text{supp}\chi_k\subset(-2,2)$, and
$\chi_{k+1}=1$ on $\text{supp}\chi_k$, $k=1,2,\dots$
For a multi-index $\alpha=(\alpha_1,\dots,\alpha_n)$, where
$\alpha_i$'s are non-negative integers, let 
$i_\alpha=\max\{i;\alpha_i\ne0\}$ and
$\alpha^\prime=(\alpha_1-\delta_{1i_\alpha},\dots,
 \alpha_n-\delta_{ni_\alpha})$, $\delta_{ij}$ being Kronecker's
 delta.
Define $\Phi_{(\alpha)}$ by
\[
    \Phi_{(\alpha)}(\,\cdot\,  ;G)
    =\Phi_{i_\alpha}(\,\cdot\,  ;\chi_1(\xi)G) 
\]
if $|\alpha|=\sum_{k=1}^n \alpha_k=1$, and
\[
    \Phi_{(\alpha)}(\,\cdot\,  ;G)
    =\Phi_{(\alpha^\prime)  }  
    \bigl( \,\cdot\,;
       \chi_{|\alpha|}(\xi)\Phi_{ i_\alpha}(\cdot;G) \bigr)
\]
if $|\alpha|\ge2$.

It then holds that
\begin{equation}\label{tn.eq3}
    {\mathbb E}[\chi(\xi)(\partial^\alpha f\circ F)G]
    ={\mathbb E}[(f\circ F)\Phi_{(\alpha)}(\,\cdot\,  ;G)]
\end{equation}
for any $f\in{\mathcal S}({\mathbb R}^n)$ and
$G\in{\mathbb D}_{\infty}$, 
where $\partial^\alpha=(\partial_1)^{\alpha_1}\cdots
                (\partial_n)^{\alpha_n}$.
In fact, for $|\alpha|=1$, (\ref{tn.eq3}) is nothing but (\ref{tn.eq1}).
Suppose that (\ref{tn.eq3}) holds for $\alpha$ with $|\alpha|=k$.
Then for $\alpha$ with $|\alpha|=k+1$, 
since 
\[
   \chi_{k+1}(\xi)\Phi_{i_\alpha}(\,\cdot\, ;  \chi_{1}(\xi) G)
   =\Phi_{i_\alpha}(\,\cdot\, ;  \chi_{1}(\xi)G),
\]
by the assumption of induction and (\ref{tn.eq1}), we obtain
\begin{align*}
    {\mathbb E}[\chi(\xi)(\partial^\alpha f \circ F)G]
    & ={\mathbb E}[\chi_1(\xi)
    (\partial_{i_\alpha}  (\partial^{\alpha^\prime} f)\circ F)G]
   \\
    &={\mathbb E}[(\partial^{\alpha^\prime}   f\circ F)
       \Phi_{i_\alpha}(\,\cdot\,   ; \chi_1(\xi) G)]
    \\
    & ={\mathbb E}[\chi_{|\alpha|}(\xi)(\partial^{\alpha^\prime}     f\circ F)
       \Phi_{i_\alpha}(\,\cdot\,  ;  \chi_1(\xi) G)]
    ={\mathbb E}[(f\circ F)\Phi_{(\alpha)}(\,\cdot\,  ;G)].
\end{align*}

On account of (\ref{tn.eq2})
 and (\ref{tn.eq3}), repeating the standard argument 
to construct Watanabe's pullback in Section V-9, \cite{iwbk}, 
we arrive at the unique existence of the continuous mapping 
$u \colon {\mathcal S}^\prime({\mathbb R}^n)\to
     \tilde{\mathbb D}_{-\infty}$ such that 
$u(f)=\chi(\xi)f(F)$ for $f\in{\mathcal S}({\mathbb R}^n)$.
Rewriting $u(T)$ as $\chi(\xi) T(F)$, we obtain the desired
continuous linear mapping.
\QED

Next, we state the asymptotic expansion theorem, which is Proposition 6.2, \cite{tw}.
Let $\{ F_{\ve}\}_{0<\ve \le 1} \subset {\mathbb D}_{\infty} ({\mathbb R}^n)$
and $\{ \xi_{\ve}\}_{0<\ve \le 1} \subset {\mathbb D}_{\infty}$
 be families
of Wiener functionals such that the following asymptotics hold:
%
\begin{align}
 F_{\ve}  
 &\sim
 f_0 +\ve f_1 + \ve^2 f_2 +\cdots 
& \quad
 & \mbox{in ${\mathbb D}_{\infty} ({\mathbb R}^n)$ as $\ve \searrow 0$,}
  \label{wat2.eq}
  \\
 \xi_{\ve}  
 &\sim 
 a_0 +\ve a_1 + \ve^2 a_2 +\cdots 
& \quad
 &
 \mbox{in ${\mathbb D}_{\infty} $ as $\ve \searrow 0$.}
  \label{wat3.eq}
 \end{align}

\begin{pr}\label{pr.comp2}
Assume (\ref{wat2.eq}), (\ref{wat3.eq}) and $|a_0| \le 1/8$.
Moreover, assume that there exists $\rho >0$ independent of $\ve$ such that 
(\ref{wat1.eq}) with $F = F_{\ve}$ and $\xi = \xi_{\ve}$ holds for any $\ve \in (0,1]$.
Let $\chi \colon {\mathbb R} \to {\mathbb R}$ be a smooth function 
whose support is contained in $[-1, 1]$ such that $\chi (x) =1$ if $|x| \le 1/2$.
Then, we have the following asymptotic expansion:
\begin{eqnarray}
 \chi(\xi_{\ve}) T (F_{\ve}) \sim 
  \Phi_0 +\ve \Phi_1 + \ve^2 \Phi_2 +\cdots 
   \qquad
 \qquad \mbox{in $\tilde{\mathbb D}_{-\infty}$ as $\ve \searrow 0$.}
 \nn
 \end{eqnarray}
\end{pr}
In the above proposition, $\Phi_k \in \tilde{\mathbb D}_{-\infty}$ can be written as 
the $k$th coefficient of the formal Taylor expansion of $T(f_0 + [ \ve f_1 + \ve^2 f_2 +\cdots ])$.
In particular, $\Phi_0 = T(f_0)$.

\Proof
Let $T\in{\mathcal S}^\prime({\mathbb R}^n)$.
Take an $m\in{\mathbb N}$ so that 
$\phi=(1+|x|^2-\frac12\Delta)^{-m}T$ is a bounded function on
${\mathbb R}^n$ which is $k$-times continuous differentiable with
bounded derivatives up to order $k$.
By virtue of (\ref{tn.eq3}), there exists a continuous linear mapping
$\ell_\varepsilon\colon {\mathbb D}_{\infty}\to{\mathbb D}_{\infty}$ such
that
\[
    {\mathbb E} [\chi(\xi_{\ve}) T(F_{\ve}) G]
    ={\mathbb E}[\phi(F_{\ve})\ell_{\ve}(G)]
    \quad\text{for every }G\in{\mathbb D}_{\infty}.
\]
By (\ref{tn.eq0}), 
\[
     \ell_{\ve}(G)=\sum_{i=0}^{2m} 
     \langle P_i(\ve),D^iG\rangle_{{\mathcal H}^{\otimes i}},
\]
where $P_i(\ve)\in{\mathbb D}_{\infty}$, $i=0,\dots,2m$, are
polynomials in $F_{\ve}$,
$\chi_k(\xi_{\ve})\sigma_{F_{\ve}}^{-1}$, $k=1,\dots,2m$,
and their derivatives.
Since
\[
        \lim_{\ve\to0}
        \frac{1}{\ve^k} \|\chi(\xi_{\ve})-1\|_{p,r}
        =0
\]
for any $p>1,r>0$ and $k\in{\mathbb N}$, 
applying the argument used in the proof of 
Theorem~V.9.4, \cite{iwbk} to
$\chi(\xi_{\ve})\phi(F_{\ve})$ instead of $\phi(F_{\ve})$, 
we obtain the desired asymptotic expansion.
\QED


At the end of this subsection, we gather well-known facts about 
SDE (\ref{sdeX.def})
and ODE (\ref{ode.def}) for later use.
The Jacobian $J_t$
of $\phi(t, x, h)$ with respect to $x$ and its inverse $K_t$ satisfy the following ODEs:
\begin{align}
dJ_t &= \sum_{i=1}^r  \nabla V_i ( \phi_t) J_t dh_t^i 
&
\quad
&
\mbox{with \quad $J_0 = {\rm Id}_d$,}
\label{ode_J.def}
\\
dK_t &= - \sum_{i=1}^r  K_t \nabla V_i ( \phi_t)  dh_t^i 
&
\quad
&
\mbox{with \quad $K_0 = {\rm Id}_d$.}
\label{ode_K.def}
\end{align}
Here, $J, K, \nabla V_i $ are all $d \times d$ matrices.
Note that $K_t =J_t^{-1}$.
When the  dependence on $h$ and $x$ needs to be specified, 
we write $J_t (h)$ or $J(t,x,h)$, etc.
The deterministic Malliavin covariance matrix is given by
\begin{equation}\label{MCode.eq}
\sigma [\phi_1 ](h) = J_1(h)  
\Bigl\{
\int_0^1  J_{t}(h)^{-1}  {\bf V} ( \phi_t (h) )  {\bf V} ( \phi_t (h) )^*
  (J_{t}(h)^{-1})^*  dt
\Bigr\}   J_1(h)^*
\end{equation}
and $\sigma [\psi_1 ](h) = \Pi_{{\cal V}} \sigma [\phi_1 ](h)  \Pi_{{\cal V}}^{*}$,
where we set ${\bf V} :=[V_1, \ldots, V_r] \in  {\rm Mat}(d,r)$ for simplicity.

Similarly, 
the Jacobian process $J^{\ve}_t$
for SDE (\ref{sdeX.def}) and its inverse $K^{\ve}_t$ satisfy the following SDEs:
\begin{align}
dJ^{\ve}_t &= \ve \sum_{i=1}^r  \nabla V_i ( X^{\ve}_t) J^{\ve}_t  \circ dw_t^i 
+\ve^2 \nabla V_0 ( X^{\ve}_t) J^{\ve}_t dt
&
\quad
&\mbox{with \quad $J^{\ve}_0 = {\rm Id}_d$,}
\label{sde_J.def}
\\
dK^{\ve}_t &= - \ve \sum_{i=1}^r  K^{\ve}_t \nabla V_i ( X^{\ve}_t) \circ dw_t^i  
 - \ve^2 K^{\ve}_t \nabla V_0 ( X^{\ve}_t)dt
&
\quad
&\mbox{with \quad $K^{\ve}_0 = {\rm Id}_d$.}
\label{sde_K.def}
\end{align}
As before, $K^{\ve}_t = (J^{\ve}_t)^{-1}$.
When necessary 
we will write $J^{\ve}_t (w)$ or $J^{\ve} (t,x,w)$, etc.
 The  Malliavin covariance matrix is given by
\begin{equation}\label{MCsde.eq}
\ve^{-2} \sigma [X^{\ve}_1 ] = J^{\ve}_1  
\Bigl\{
\int_0^1  (J^{\ve}_{t})^{-1}  {\bf V} ( X^{\ve}_t  )  {\bf V} ( X^{\ve}_t )^*
  (J^{\ve}_{t})^{-1, *}  dt
\Bigr\}   (J^{\ve}_1)^*
\end{equation}
and $\sigma [Y^{\ve}_t  ](w) = \Pi_{{\cal V}} \sigma [X^{\ve}_t  ](w)  \Pi_{{\cal V}}^{*}$
a.s.
(The dependence on $w$ is suppressed above.)


\subsection{Preliminaries from rough path theory}

In this subsection we recall the geometric rough path space 
with H\"older or Besov norm 
and quasi-sure property of rough path lift.
For basic properties of geometric rough path space, 
we refer to Lyons, Caruana, and L\'evy \cite{lcl},
and Friz and Victoir \cite{fvbk}.
For the geometric rough path space with Besov norm, 
we refer to Appendix A.2, \cite{fvbk}.
Quasi-sure property of rough path lift is summarized in Inahama \cite{in2}.

We assume that 
the Besov parameters $(\alpha, 4m)$ satisfy the following condition:
\begin{equation}\label{eq.amam}
\frac13 <\al < \frac12, \quad m =1,2,3,\ldots, \quad
\al - \frac{1}{4m} > \frac13, \quad
\mbox{ and } \quad 4m (\frac12 -\al)   >1.
\end{equation}
We choose such a pair $(\alpha, 4m)$ and fix it throughout this paper.


%
%
We denote by  $G\Omega^H_{\alpha} ( {\mathbb R}^r) $, $1/3 <\alpha < 1/2$,  the geometric rough path space 
over ${\mathbb R}^d$ with $\alpha$-H\"older norm.
Let $C_0^{\beta-H}([0,1], {\mathbb R})$, $0<\beta \le 1$,
be the Banach space of all the ${\mathbb R}$-valued, $\beta$-H\"older continuous paths 
that start at $0$.
If $\alpha + \beta >1$, then the Young pairing 
\[
 G\Omega^H_{\alpha} ( {\mathbb R}^r) \times C_0^{\beta-H}([0,1], {\mathbb R})
  \ni ({\bf w}, \lambda) \mapsto
  ({\bf w}, \bm{\lambda}) 
     \in  G\Omega^H_{\alpha} ( {\mathbb R}^{r+1})
     \]
is a well-defined, locally Lipschitz continuous map.
(See Section 9.4, \cite{fvbk} for instance.)



Now we consider a system of RDEs driven by the Young pairing
 $({\bf w}, \bm{\lambda}) \in  G\Omega^H_{\alpha} ( {\mathbb R}^{d+1})$
of ${\bf w} \in G\Omega^H_{\alpha} ( {\mathbb R}^{d})$
and 
$\lambda \in C_0^{1-H}([0,1], {\mathbb R}^1)$.
(In most cases, we will assume $\lambda_t = \mbox{const} \times t$.)
For vector fields $V_{i}\colon {\mathbb R}^d \to {\mathbb R}^d$ ($0 \le i \le r$), consider
\begin{equation}\label{rde_x.def}
dx_t = \sum_{i=1}^r  V_i ( x_t) dw_t^i + V_0 ( x_t) d\lambda_t
\qquad
\qquad
\mbox{with  \quad $x_0 =x \in {\mathbb R}^d$.}
\end{equation}
The RDEs for the Jacobian process and its inverse are given as follows;
\begin{align}
dJ_t &= \sum_{i=1}^r  \nabla V_i ( x_t) J_t dw_t^i + \nabla V_0 ( x_t) J_t d\lambda_t
&
\quad
&
\mbox{with $J_0 ={\rm Id}_d \in {\rm Mat}(d,d)$,}
\label{rde_J.def}
\\
dK_t &= - \sum_{i=1}^r K_t  \nabla V_i ( x_t)  dw_t^i  - K_t \nabla V_0 ( x_t) d\lambda_t
&
\quad
&
\mbox{with  $K_0 ={\rm Id}_d \in {\rm Mat}(d,d)$.}
\label{rde_K.def}
\end{align}
Here, 
$J, K,$ and  $\nabla V_i $ are all ${\rm Mat}(d,d)$-valued.

Assume that $V_i$'s are of $C_b^{4}$ for a while.
Then, a global solution of (\ref{rde_x.def})--(\ref{rde_K.def}) 
exists for any ${\bf w}$ and $\lambda$ . 
Moreover, Lyons' continuity theorem holds.
(The linear growth case is complicated and will be discussed later.)
In that case, the following maps are continuous:
\begin{eqnarray}
 G\Omega^H_{\alpha} ( {\mathbb R}^r) \times C_0^{1-H}([0,1], {\mathbb R}^1)
  \ni ({\bf w}, \lambda) &\mapsto&
  ({\bf w}, {\bm \lambda}) \in  G\Omega^H_{\alpha} ( {\mathbb R}^{r+1})
      \nn\\
      &\mapsto& 
  ({\bf w},{\bm \lambda};{\bf x}, {\bf J}, {\bf K}) \in
   G\Omega^H_{\alpha} ( {\mathbb R}^{r+1} \oplus
   {\mathbb R}^d  \oplus {\rm Mat}(d,d)^{\oplus 2})
   \nn\\
   &\mapsto&
  ({\bf x}, {\bf J}, {\bf K}) \in
   G\Omega^H_{\alpha}  ( {\mathbb R}^d  \oplus {\rm Mat}(d,d)^{\oplus 2}).
       \label{3arrows.eq}
       \end{eqnarray}
Here, the first map is the Young pairing,
the second is  the Lyons-It\^o map, 
and the third is the canonical projection.
(The map $({\bf w}, \lambda) \mapsto  {\bf x}$ will be denoted by 
$\Phi \colon G\Omega^H_{\alpha} ( {\mathbb R}^r ) \times C_0^{1-H}([0,1], {\mathbb R}) 
\to G\Omega^H_{\alpha} ( {\mathbb R}^d$).)
Recall that in Lyons' formulation of rough path theory, the initial values of the first level paths 
must be adjusted.
Note that
$({\rm Id}+{\bf J}^1_{0,t})^{-1}= {\rm Id}+{\bf K}^1_{0,t}$ always holds.

When ${\bf w}$ is the natural lift $h \in {\cal H}$ and $\lambda \equiv 0$,
the first level path
\begin{equation}\label{rp.yaji.eq}
 t \mapsto ( x + {\bf x}^1_{0,t}, {\rm Id}+{\bf J}^1_{0,t}, {\rm Id}+{\bf K}^1_{0,t} )
\end{equation}
is identical to the solution of a system (\ref{ode.def}), (\ref{ode_J.def}), (\ref{ode_K.def}).
Here, ${\bf x}^1_{0,t}$ is the first level path of ${\bf x}$ evaluated at $(0,t)$, etc.
Similarly, if $({\bf w}, \lambda) = (\ve {\bf W}, \lambda^{\ve})$,
where ${\bf W} ={\cal L} (w)$
is the Brownian rough path  under $\mu$ and $ \lambda^{\ve}_t =\ve^2 t$,
then (\ref{rp.yaji.eq}) coincides
with the solution of (\ref{sdeX.def}), (\ref{sde_J.def}), (\ref{sde_K.def}) a.s.


We define a continuous function
$\Gamma\colon  G\Omega^H_{\alpha} ( {\mathbb R}^r) \times C_0^{1-H}([0,1], {\mathbb R})
\to {\rm Mat}(d,d)$ 
as follows:
Set 
\begin{equation}\label{MCrde.eq}
\Gamma ({\bf w}, \lambda)  
= ({\rm Id}+{\bf J}^1_{0,t}) \hat\Gamma ({\bf w}, \lambda) ({\rm Id}+{\bf J}^1_{0,t})^*,
\end{equation}
where
\[
\hat\Gamma ({\bf w}, \lambda)
:=
\int_0^1  ({\rm Id}+{\bf K}^1_{0,t}) {\bf V} ( x + {\bf x}^1_{0,t})  
{\bf V}( x + {\bf x}^1_{0,t})^*  ({\rm Id}+{\bf K}^1_{0,t})^* dt
\]
with ${\bf V} :=[V_1, \ldots, V_r] \in  {\rm Mat}(d,r)$.

From (\ref{MCode.eq}) and (\ref{MCsde.eq}) 
we can easily see the following:
If $\lambda_t^{\ve} = \ve^2 t$, 
then $\Gamma (\ve {\bf W},  \lambda^{\ve}) 
= \ve^{-2}\sigma [X_1^{\ve} ] (w)$ 
for $\mu$-almost all $w$, 
where  $X_1^{\ve}$ denotes the solution of SDE (\ref{sdeX.def}) at $t=1$.
If $\lambda_t \equiv 0$ and ${\bf h}={\cal L}(h)$ is the natural lift of $h \in {\cal H}$, 
then
$\Gamma ({\bf h}, 0) = \sigma [\phi_1 ](h)$, the deterministic Malliavin covariance matrix 
given in (\ref{MCode.eq}).
From these we can easily see that 
$\Pi_{{\cal V}}\Gamma (\ve {\bf w},  \lambda^{\ve}) \Pi_{{\cal V}}^*= \ve^{-2}\sigma [Y_1^{\ve} ] (w)$ a.s.
and
$\Pi_{{\cal V}}  \Gamma ({\bf h}, 0) \Pi_{{\cal V}}^*= \sigma [\psi_1 ](h)$.

\begin{re}\normalfont
In this paper we will use Lyons' continuity theorem 
only with respect to $\alpha$-H\"older topology $(1/3 <\alpha <1/2)$ and 
for $C^4_b$-coefficient vector fields. 
We do not try to extend it to the case of 
unbounded coefficient vector fields or Besov topology.
\end{re}

Now we introduce the Besov topology on the rough path space.
For $(\al, 4m)$ which satisfies (\ref{eq.amam}), 
$G\Omega^B_{\alpha, 4m} ( {\mathbb R}^r) $ denotes the geometric rough path space 
over ${\mathbb R}^r$ with $(\al, 4m)$-Besov norm.
Recall that the distance on this space is given by 
\begin{align}
d({\bf w}, \hat{\bf w}) 
&= \| {\bf w}^1- \hat{\bf w}^1 \|_{\al, 4m-B} 
+\| {\bf w}^2- \hat{\bf w}^2 \|_{2\al, 2m-B}
\nn\\
&
:=
\Bigl(
\iint_{0 \le s <t \le 1}  \frac{ | {\bf w}^1_{s,t}- \hat{\bf w}^1_{s,t}|^{4m}}
{|t-s|^{1 +4m\al }} 
dsdt
\Bigr)^{1/4m}
+
\Bigl(
\iint_{0 \le s <t \le 1}  \frac{ | {\bf w}^2_{s,t}- \hat{\bf w}^2_{s,t}|^{2m}}
{|t-s|^{1 +4m\al }} 
dsdt
\Bigr)^{1/2m}.
\nn
\end{align}
By the Besov-H\"older embedding theorem for rough path spaces,
there is a continuous embedding $G\Omega^B_{\alpha, 4m} ( {\mathbb R}^r)
 \hookrightarrow G\Omega^H_{\alpha -(1/4m)} ( {\mathbb R}^r)$.
If $\al < \al' <1/2$, there is a continuous embedding
$G\Omega^H_{\alpha'} ( {\mathbb R}^r)
\hookrightarrow G\Omega^B_{\alpha, 4m} ( {\mathbb R}^r)$.
We remark that 
we will not write these embeddings explicitly.
(For example, if we write $\Phi({\bf w}, \lambda)$ for  
$({\bf w}, \lambda) \in G\Omega^B_{\alpha, 4m} ( {\mathbb R}^r) \times C_0^{1-H}([0,1], {\mathbb R})$,
then it is actually the composition of the first embedding map above and $\Phi$ with respect to 
$\{\al -1/(4m)\}$-H\"older topology.)

Note also that the Young translation by $h \in {\cal H}$ works 
perfectly on $G\Omega^B_{\alpha, 4m} ( {\mathbb R}^r)$ under (\ref{eq.amam}).
The map $({\bf w}, h) \mapsto \tau_h ({\bf w})$
is continuous from $G\Omega^B_{\alpha, 4m} ( {\mathbb R}^r) \times {\cal H}$
to $G\Omega^B_{\alpha, 4m} ( {\mathbb R}^r)$, 
where 
$ \tau_h ({\bf w})$ is called the Young translation of ${\bf w}$ by $h$
and is defined by
\[
\tau_{h}( {\bf w})^1_{s,t} := {\bf w}^1_{s,t} + {\bf h}^1_{s,t},
\qquad
\tau_{h}( {\bf w})^2_{s,t} := {\bf w}^2_{s,t} + {\bf h}^2_{s,t} 
+ 
\int_s^t  {\bf w}^1_{s,u} dh_u + \int_s^t  {\bf h}^1_{s,u} dw_u.
\]
Here, the integrals are in the Young (or Riemann-Stieltjes) sense and we set
$w_t :={\bf w}^1_{0,t}$.
Moreover, there exists a positive constant $C=C_{\alpha, 4m}$ such that 
\begin{eqnarray}
\| \tau_{h}( {\bf w})^1 \|_{\alpha, 4m -B}
&\le& 
\|  {\bf w}^1 \|_{\alpha, 4m -B}
+C\| h \|_{{\cal H}},
\label{shift1.ineq}
\\
\| \tau_{h}( {\bf w})^2 \|_{2\alpha, 2m -B} 
&\le& 
\|  {\bf w}^2 \|_{2\alpha, 2m -B} +2 C \|  {\bf w}^1 \|_{\alpha, 4m -B}\|h\|_{{\cal H}}
+ C^2  \| h \|_{{\cal H}}^2
\label{shift2.ineq}
\end{eqnarray}
hold for all $h \in {\cal H}$ and ${\bf w} \in G\Omega^B_{\alpha, 4m}({\mathbb R}^r)$.


For $\gamma >0$ and $h \in {\cal H}$, we set 
\begin{eqnarray}
U_{h, \gamma} &=&
\{ {\bf w} \in G\Omega^B_{\alpha, 4m}({\mathbb R}^r) 
\mid  \| \tau_{-h}( {\bf w})^1 \|_{\alpha, 4m -B}^{4m}+ \|  \tau_{-h}( {\bf w})^2 \|_{2\alpha, 2m -B}^{2m}
< \gamma^{4m}\},
\label{Uhg.def}
\\
U_{h, \gamma}^{\prime} 
&=&
U_{h, 2^{-1/4m}\gamma}
=
\{ {\bf w} \in G\Omega^B_{\alpha, 4m}({\mathbb R}^r) 
\mid  \| \tau_{-h}( {\bf w})^1 \|_{\alpha, 4m-B}^{4m}+ \|  \tau_{-h}( {\bf w})^2 \|_{2\alpha, 2m -B}^{2m}
< \gamma^{4m} /2
\}.
\nn
\end{eqnarray}
When $h=0$ we simply write $U_{\gamma}$ and $U_{\gamma}^{\prime}$.
By the continuity of $\tau_h$ with respect to the $(\alpha, 4m)$-Besov rough path topology, 
$\{ U_{h, \gamma} \}_{\gamma >0}$ forms a system of open neighborhoods
around ${\bf h} ={\cal L}(h)$
and so does $\{ U_{h, \gamma}^{\prime} \}_{\gamma >0}$.
Clearly, 
$U_{h, \gamma}^{\prime} = U_{h, \gamma 2^{-1/4m}}$,
$U_{h, \gamma} = \tau_h (U_{\gamma})$ and
$U_{h, \gamma}^{\prime} = \tau_h (U_{\gamma}^{\prime})$.
%


Now we discuss quasi-sure properties of rough path lift map 
${\cal L}$ from ${\cal W}$ to $G\Omega^B_{\alpha, 4m} ( {\mathbb R}^r)$.
For $k=1,2,\ldots$ and $w \in {\cal W}$,
we denote by $w(k)$ the $k$th dyadic piecewise linear approximation of $w$ 
associated with the partition $\{ l2^{-k} \mid 0 \le l \le 2^k\}$ of $[0,1]$.
We set 
\[
{\cal Z}_{\al, 4m} := \bigl\{ w \in {\cal W}\mid
\mbox{ $\{ {\cal L} (w(k)) \}_{k=1}^{\infty}$ is Cauchy in $G\Omega^B_{\alpha, 4m} ( {\mathbb R}^r )$} 
\bigr\}.   
\] 
We define ${\cal L} \colon {\cal W} \to G\Omega^B_{\alpha, 4m} ( {\mathbb R}^r )$
by ${\cal L} (w) = \lim_{m\to \infty} {\cal L} (w(k))$ if $w \in {\cal Z}_{\al, 4m}$
and 
we do not define ${\cal L} (w)$ if $w \notin {\cal Z}_{\al, 4m}$.
(We will always use this version of ${\cal L}$.)
Note that 
${\cal H}$ and $C_0^{\beta-H}([0,1], {\mathbb R}^r)$ with $\beta \in (1/2, 1]$
are subsets of ${\cal Z}_{\al, 4m}$
and this lift for elements of such subsets
coincides 
with the direct lift by means of the Riemann-Stieltjes integral.
Under scalar multiplication and  Cameron-Martin translation,  
 ${\cal Z}_{\al, 4m}$ is invariant.
Moreover, $c {\cal L} (w) = {\cal L} (cw)$ and $\tau_h({\cal L} (w))= {\cal L} (w+h)$ 
for any $w \in {\cal Z}_{\al, 4m}$, $c \in {\mathbb R}$, and $h \in {\cal H}$.


It is known that ${\cal Z}_{\al, 4m}^c$ is slim, that is the $(p,r)$-capacity of 
this set is zero for any $p \in (1,\infty)$ and $r \in {\mathbb N}$.
(See Aida \cite{ai}, Inahama \cite{in1, in2}.)
Therefore, from a viewpoint of quasi-sure analysis, 
the lift map ${\cal L}$ is well-defined.
Moreover, 
the map ${\cal W} \ni w \mapsto {\cal L} (w) \in G\Omega^B_{\alpha, 4m} ( {\mathbb R}^r) $
is $\infty$-quasi-continuous
(Aida \cite{ai}).
We will often write ${\bf W} := {\cal L} (w)$ when it is regarded 
as a rough path space-valued random variable defined on ${\cal W}$.
Due to Lyons' continuity theorem and uniqueness of quasi-continuous modification, 
$\tilde{X}^{\ve} (\,\cdot\, , x,w) = x + \Phi(\ve {\cal L}(w), \ve^2\lambda)^1$ holds quasi-surely
if $V_i~(0 \le i \le d)$ is of $C^3_b$.
(Here, $\lambda_t = t$ and $\tilde{X}^{\ve} (\,\cdot\, , x,w)$ denotes the 
$\infty$-quasi continuous modification of the path space-valued random variable
$w \mapsto X^{\ve} (\,\cdot\, , x,w)$.)
A similar remark holds for $J^{\ve}$ and $K^{\ve}$.


Before closing this subsection, we give a remark for the coefficient vector fields 
with linear growth.
\begin{re}\normalfont
If $V_i~(0 \le i \le d)$ satisfies {\bf (A1)} and has linear growth,
it is not easy to prove the existence of a global solution of RDE (\ref{rde_x.def})
for a given $({\bf w}, \lambda)\in G\Omega^H_{\alpha} ( {\mathbb R}^r) \times C_0^{1-H}([0,1], {\mathbb R})$.
However, Bailleul \cite{be} recently proved it.
Hence, the Lyons-It\^o map $\Phi$ can be defined  
on the whole space $G\Omega^H_{\alpha} ( {\mathbb R}^r) \times C_0^{1-H}([0,1], {\mathbb R})$
and Lyons' continuity theorem holds under {\bf (A1)}, too.
However, almost no other properties of  the solution are known under {\bf (A1)}
at this moment.
Therefore, we basically assume that $V_i~(0 \le i \le r)$ is of  $C_b^{\infty}$
when we use RDEs.
The reason why we may do so is as follows.
In the proof of our main theorem
we will discard contributions from rough paths distant from ${\cal L} ({\cal K}_a^{\min})$
by using large deviation theory.
Therefore, by simple cut-off argument, 
the problem reduces to the $C_b^{\infty}$-case anyway.
\end{re}

\subsection{Some useful results in Malliavin calculus on rough path space}

First,
we  recall Taylor-like  expansion of (Lyons-)It\^o map.
In the proof of our main theorem, 
we will compute the Cameron-Martin translation of $X^{\ve}$.
For a Cameron-Martin path $h \in {\cal H}$, 
we set $X^{\ve, h}_t :=X^{\ve} (t,x,w + (h/\ve))$, which satisfies 
\begin{equation}\label{sdeX.shift}
dX^{\ve,h}_t = \sum_{i=1}^r  V_i ( X^{\ve, h}_t) \circ  (\ve dw_t^i +dh^i_t ) + \ve^2  V_0 (X^{\ve, h}_t)   dt
\qquad
\qquad
\mbox{with \quad $X^{\ve,h}_0 =x \in {\mathbb R}^d$.}
\end{equation}
Similarly,  the translations of the Jacobian process and its inverse satisfy:
\begin{equation}\label{sdeJ.shift}
dJ^{\ve, h}_t = \sum_{i=1}^r  \nabla V_i ( X^{\ve, h}_t) J^{\ve, h}_t\circ  (\ve dw_t^i  +dh^i_t)
+ \ve^2  \nabla V_0 (X^{\ve, h}_t)  J^{\ve, h}_t dt
\qquad
\mbox{with \quad $J^{\ve, h}_0 ={\rm Id}_d$}
\end{equation}
and
\begin{eqnarray}
dK^{\ve, h}_t = -\sum_{i=1}^r  K^{\ve, h}_t \nabla V_i ( X^{\ve, h}_t) \circ  (\ve dw_t^i  +dh^i_t)
- \ve^2 K^{\ve, h}_t \nabla V_0 (X^{\ve, h}_t)  dt
\quad
\mbox{with  $K^{\ve, h}_0 ={\rm Id}_d$.}
\label{sdeK.shift}
\end{eqnarray}
It is easy to see that  $(X^{\ve,h},  J^{\ve, h}, K^{\ve, h})$ coincides a.s. with  (\ref{rp.yaji.eq}) 
with the driving rough path being 
$( \tau_h (\ve {\bf W}), \lambda^{\ve}) = ( \ve {\cal L} ( w + h/\ve) , \lambda^{\ve})$.


It is known that under {\bf (A1)} the following 
asymptotic expansion holds in ${\mathbb D}_{\infty} ({\mathbb R}^d)$ as $\ve \searrow 0$:
\begin{equation}\label{asympX.1}
X^{\ve} (1,x,w + \frac{h}{\ve}) = f_0(h) + \ve f_1 (w;h) + \cdots +\ve^k f_k (w;h) + Q^{\ve}_{k+1} (w;h)
\end{equation}
where $Q^{\ve}_{ k+1} (w;h) =O(\ve^{k+1})$.
Note that $f_0(h) = \phi_1 (h)$.
It immediately  follows that
\begin{equation}\label{asympX.2}
Y^{\ve} (1,x,w + \frac{h}{\ve}) = g_0(h) + \ve g_1 (w;h) + \cdots +\ve^k g_k (w;h) + R^{\ve}_{k+1} (w;h)
\end{equation}
in ${\mathbb D}_{\infty} ({\mathbb R}^d)$ as $\ve \searrow 0$,
where we set 
$R^{\ve}_{ k+1} (w;h) =\Pi_{ {\cal V}}  (Q^{\ve}_{ k+1} (w;h))$ and $g_j (w;h) = \Pi_{ {\cal V}}  (f_j (w;h) )$
for all $j \ge 0$.
Note that $g_0(h) = \psi_1 (h)$ and $ g_1 (w;h)=D\psi_1 (h) \la w\ra$.
Moreover, 
these asymptotic expansions are uniform in $h$ as $h$
varies in an arbitrary bounded set in ${\cal H}$.

The expansion (\ref{asympX.1}) has a counterpart in rough path theory, 
which we call Taylor-like  expansion of the Lyons-It\^o map $\Phi$ 
(See Inahama and Kawabi \cite{ik} and Inahama \cite{ina10}.)
Assume $1/3 <\alpha <1/2$
and  that $V_i$ is of $C_b^{\infty}$ for a while ($0 \le i \le r$).
Then, there exist  continuous maps 
$\hat{f}_k \colon G\Omega^H_{\alpha} ( {\mathbb R}^r ) \times {\cal H} \mapsto {\mathbb R}^d$
and 
$\hat{Q}^{\ve}_{k+1} \colon G\Omega^H_{\alpha} ( {\mathbb R}^r ) \times {\cal H} \mapsto {\mathbb R}^d$
for $k =0,1,\ldots$  which satisfy the following {\rm (i)} and {\rm (ii)}:
\\
\noindent
{\rm (i)}~For all $k  =0,1,\ldots$
and $\ve \in [0,1]$,
\begin{equation}\label{tylr1.eq}
f_k (w;h) = \hat{f}_k ({\bf W};h)
\quad
\mbox{and}
\quad
Q^{\ve}_{k+1} (w;h) = \hat{Q}^{\ve}_{k+1} ({\bf W};h),
\quad
\mbox{$\mu$-a.s.,}
\end{equation}
where we set 
\begin{equation}\label{tylr_rp.eq}
\hat{Q}^{\ve}_{k+1} ({\bf w};h)
:=
x + \Phi ( \tau_h (\ve {\bf w}), \lambda^{\ve} )^1_{0,1}
-\bigl\{
\hat{f}_0(h) + \ve \hat{f}_1 ({\bf w};h) + \cdots +\ve^k \hat{f}_k ({\bf w};h)
\bigr\}.
\end{equation}
\\
\noindent
{\rm (ii)}~For any $k \in {\mathbb N}$, $\rho>0$,  $h \in {\cal H}$, 
there exist positive constants $C_k$ and $C'_{k, \rho}$ such that
\begin{eqnarray}
|  \hat{f}_k ({\bf w};h)| &\le&
 C_k (1 + \| {\bf w}^1 \|_{\alpha -H}  +  \| {\bf w}^2 \|_{2\alpha -H}^{1/2}  )^k
\label{tylr2.eq}
\end{eqnarray}
for any ${\bf w} \in G\Omega^H_{\alpha} ( {\mathbb R}^r )$
and
\begin{eqnarray}
| \hat{Q}^{\ve}_{k+1} ({\bf w};h) |
&\le&
C'_{k,\rho}  (\ve + \| (\ve {\bf w} )^1 \|_{\alpha -H}  +  \| (\ve {\bf w} )^2 \|_{2\alpha -H}^{1/2} )^{k+1}
\label{tylr3.eq}
\end{eqnarray}
if $\| (\ve {\bf w} )^1 \|_{\alpha -H}  +  \| (\ve {\bf w} )^2\|^{1/2} \le \rho$.
Of course,  $C_k$ and $C'_{k, \rho}$ also depend on $h$, 
but in fact they  depend only on $\|h\|_{{\cal H}}$.

Note that  the inequalities  (\ref{tylr2.eq}) and (\ref{tylr3.eq}) are deterministic.
It is obvious that $\hat{g}_k (w;h) := \Pi_{ {\cal V}}  (\hat{f}_k (w;h) )$
and 
$\hat{R}^{\ve}_{ k+1} ({\bf w};h) :=\Pi_{ {\cal V}}  (\hat{Q}^{\ve}_{ k+1} ({\bf w};h))$ 
have similar properties to (\ref{tylr1.eq})--(\ref{tylr3.eq}).

\begin{re}\normalfont\label{re.tylr}
{\rm (i)}~Heuristically, (\ref{tylr1.eq}) means that,
in a sense, $\hat{f}_k$ and $\hat{Q}^{\ve}_{ k+1}$
are  ``lifts" of $f_k$ and $Q^{\ve}_{ k+1}$,
respectively.
\\
{\rm (ii)}~By abusing notations, 
we will simply write $f_k, g_k, Q^{\ve}_{ k+1}, R^{\ve}_{ k+1}$ 
for $\hat{f}_k, \hat{g}_k, \hat{Q}^{\ve}_{ k+1}, \hat{R}^{\ve}_{ k+1}$ in the sequel.
\\
{\rm (iii)}~If we assume that $V_i$ is of $C_b^{\infty}$  ($0 \le i \le r$),
then we can actually obtain (\ref{asympX.1}) via
the deterministic expansion (\ref{tylr1.eq})--(\ref{tylr3.eq}) on the geometric rough path space
(See Inahama \cite{in2}.)
Under {\bf (A1)},
however, it is not known yet whether this is possible or not.
\\
{\rm (iv)}~There is a simple expression of $\hat{f}_k ({\bf w}, h)$
when the drift vector field $V_0 \equiv 0$. 
Suppose that 
 $\xi_l \in {\cal H} ~(l=1,2,\ldots)$ converges to ${\bf w}$ in $G\Omega^H_{\alpha} ( {\mathbb R}^r )$.
 Then, 
 \[
 \hat{f}_k ({\bf w}, h) = \lim_{l \to \infty} \frac{1}{k !} D^k \phi_1 (h) \la \xi_l, \ldots, \xi_l \ra,
 \qquad
 \qquad
 \mbox{($k$-times)}.
  \]
Here, $D$ stands for the Fr\'echet derivative on ${\cal H}$.
(Even when $V_0$ does not vanish, this kind of expression exists.
Since it looks quite complicated due the factor $\ve^2$ in front of $V_0$,
we omit it, however.)
 \end{re}

In order for $F (X^{\ve} (1, x, w +h/\ve))$ 
to admit expansions as above,
$F$ need not be a projection or a linear map.
In the next lemma, we show that for quite general $F$, 
it admits expansions in both senses.
The proof is straightforward. So we omit it.
We also remark that the expansions in Lemma \ref{lm.tnki_F} are
uniform in $h$ as $h$ varies in any bounded set in ${\cal H}$.
(Lemma \ref{lm.tnki_F} will be used in the proof of the manifold case, not in 
the Euclidean case.)
\begin{lm}\label{lm.tnki_F}
Let $F \colon {\mathbb R}^d \to {\mathbb R}^e$ be a smooth map 
such that $\nabla^j F$ is bounded for all $j \ge 1$.
Then, we have the following {\rm (i)}--{\rm (ii)}:
\\
{\rm (i)}~ As $\ve \searrow 0$,
\begin{equation}\nn
F( X^{\ve} (1,x,w + \frac{h}{\ve}) )
= g^F_0(h) + \ve g^F_1 (w;h) + \cdots +\ve^k g^F_k (w;h) + R^{\ve, F}_{k+1} (w;h)
\end{equation}
in ${\mathbb D}_{\infty} ({\mathbb R}^d)$.
Here, $R^{\ve, F}_{ k+1} (w;h) =O(\ve^{k+1})$ 
in ${\mathbb D}_{\infty} ({\mathbb R}^d)$
and
$g^F_j (w;h)$ are determined by the formal composition 
of (\ref{asympX.1}) and the Taylor expansion of $F$.
For example,
\begin{eqnarray*}
g^F_0(h) &=& F (f_0(h)), 
\\
g^F_1 (w;h) &=& \nabla F (f_0(h)) \la f_1 (w;h) \ra,
\\
g^F_2 (w;h) &=& \nabla F (f_0(h)) \la f_2 (w;h) \ra 
+ \frac12 \nabla^2 F (f_0(h)) \la f_1 (w;h),  f_1 (w;h)\ra
\end{eqnarray*}
and so on.
\\
{\rm (ii)}~Similarly, if  $\hat{g}^F_j ({\bf w};h)$ is determined
by  the formal composition 
of (\ref{tylr_rp.eq}) and the Taylor expansion of $F$, 
then we have
\[
F (x + \Phi ( \tau_h (\ve {\bf w}), \lambda^{\ve} )^1_{0,1})
=
\hat{g}^F_0(h) + \ve \hat{g}^F_1 ({\bf w};h) + \cdots +\ve^k \hat{g}^F_k ({\bf w};h)
+
\hat{R}^{\ve, F}_{k+1} ({\bf w};h), 
\]
where $\hat{g}^F_k$ and $\hat{R}^{\ve, F}_{k+1}$
satisfy essentially the same estimates as in (\ref{tylr2.eq}) and (\ref{tylr3.eq})
(for different positive constants $C_k$ and $C'_{k, \rho}$).
Moreover, we have 
$\hat{g}^F_j ({\bf W};h) = g^F_j (w;h)$ and 
$
R^{\ve, F}_{k+1} (w;h)=\hat{R}^{\ve, F}_{k+1} ({\bf W};h)$, $\mu$-a.s.
\end{lm}


Next we will give a simple lemma for Malliavin covariance 
matrices of $X^{\ve} (1, x, w +h/\ve)$ and 
$Y^{\ve} (1, x, w +h/\ve)$.
For any $h \in {\cal H}$,
the Malliavin covariance matrix of 
$Q^{\ve}_{ 1} (w;h)/\ve = \{ X^{\ve} (1,x,w + h/\ve) - f_0(h) \}/\ve$
is given by
\[
\sigma [Q^{\ve}_{ 1} (\,\cdot\, ;h)  /\ve](w) 
= \int_0^1 J^{\ve, h}_1 K^{\ve, h}_t  {\bf V} (X^{\ve, h}_t) 
{\bf V} (X^{\ve, h}_t)^{*} (K^{\ve, h}_t)^{*} (J^{\ve, h}_1)^{*} dt
=
\Gamma ( \tau_h (\ve {\bf W}), \lambda^{\ve}), 
\]
where $\lambda^{\ve}_t =\ve^2 t$, ${\bf W} ={\cal L} (w)$ is Brownian rough path
and $\Gamma$ is defined in (\ref{MCrde.eq}).
Note that  
$
\sigma  [R^{\ve}_{ 1} (\,\cdot\, ;h)/\ve](w) 
= \Pi_{{\cal V}}\sigma [ Q^{\ve}_{ 1}(\,\cdot\, ;h)  /\ve ](w)  \Pi_{{\cal V}}^*
$.


Let $h \in {\cal K}^{min}_a$ and let
$\{e_i\}_{i=1}^{n'}$ be an orthonormal basis of $T_h {\cal K}^{min}_a \subset {\cal H}$.
Here, $n' = \dim {\cal K}^{min}_a$.
Set 
\begin{equation}
{\bf i}_h \la  w\ra = \sum_{i=1}^{n'} \la  e_i, w\ra_{{\cal H}}  e_i.
\label{ih.def}
\end{equation}
In fact, $\la  e_i, \,\cdot\,  \ra_{{\cal H}} \in {\cal W}^*$, 
as we will see later.
Then, 
the Malliavin covariance matrix of  $(R^{\ve}_{ 1}(w;h)/\ve, {\bf i}_h ) 
\in {\mathbb D}_{\infty} ({\cal V} \times T_h {\cal K}^{min}_a )$
is given by 
\[
\sigma [(R^{\ve}_{ 1} (\,\cdot\, ;h) /\ve, {\bf i}_h) ]
=
\begin{pmatrix}
\sigma [ R^{\ve}_{ 1} (\,\cdot\, ;h) /\ve] & \Xi\\
\Xi^* & {\rm Id}_{n'}
\end{pmatrix},
\]
where $\Xi \in {\cal V} \otimes (T_h {\cal K}^{min}_a )^*$ which is defined by 
$e \mapsto \Pi_{{\cal V}}\int_0^1 J^{\ve, h}_1 K^{\ve, h}_t  {\bf V} (X^{\ve, h}_t) e^{\prime} (t) dt$.


The following lemma is essentially Lemma 6.3, \cite{tw}, rewritten in a rough path way.
Let 
\begin{equation}\label{lamb.def}
\lambda :=\inf_{h \in {\cal K}^{min}_a} \inf_{z \in {\cal V}:  |z|=1} z^* \sigma [\psi_1 ](h) z >0
\end{equation}
be the infimum over $h$ of the smallest eigenvalue of $\sigma [\psi_1] (h)$.

\begin{lm}\label{lm.tw63}
There exist $\gamma_0, \ve_0 \in (0,1]$ such that
the smallest eigenvalue of $\sigma [(R^{\ve}_{ 1} (\,\cdot\, ;h)/\ve, {\bf i}_h ) ] (w)$ 
is greater than $(\lambda \wedge 1)/2$
if $\| ( \ve {\bf W})^1 \|_{\alpha, 4m}^{4m}+ \|  (\ve  {\bf W})^2 \|_{2\alpha, 2m}^{2m}
< 
\gamma_0^{4m}$ (i.e. $ \ve {\bf W} \in U_{\gamma_0}$), 
$h \in {\cal K}^{min}_a$ and $0 \le \ve \le \ve_0$.
Moreover, both $\gamma_0$ and $\ve_0$ can be chosen independent of $h \in {\cal K}^{min}_a$.
\end{lm}

\Proof
The covariance matrix $\sigma [(R^{\ve}_{ 1} (\,\cdot\, ;h)/\ve, {\bf i}_h ) ] (w)$ 
is of the form 
$A (\tau_{h}(\ve {\bf W}), \bm\lambda^{\ve} )$, 
where 
\\
\noindent
{\rm (i)}~$A\colon  G\Omega^B_{\alpha, 4m}({\mathbb R}^{r+1}) \to 
L ( {\cal V} \times T_h {\cal K}^{min}_a, {\cal V} \times T_h {\cal K}^{min}_a) \cong
 {\rm Mat}(n+n', n+n')$
is a certain continuous map,
\\
\noindent
  {\rm (ii)}~$\tau_{h}$ stands for the Young translation by $h$, 
\\
\noindent
  {\rm (iii)}~ $( \tau_{h}(\ve {\bf W}), \bm\lambda^{\ve} )$
   stands for the Young pairing of $\tau_{h}(\ve {\bf W})$ and $\lambda^{\ve}$.
       At $({\bf h}, 0)$ (or equivalently when $\ve =0$) the covariance matrix becomes
       \[
       A ({\bf h}, 0) = \begin{pmatrix}
\sigma [\psi_1] (h) & O\\
O & {\rm Id}_{n'}
\end{pmatrix},       \]
which is clearly greater than or equal to $\lambda \wedge 1$ as a quadratic form.
(Since $T_h {\cal K}^{min}_a \subset \ker D\psi_1 (h)$,
 the  ``off-diagonal" components are zero.)
      Hence we can find $\gamma_0$ and $\ve_0$ as in the statement for each $h$. 
       
       Since $A$ and the Young translation/pairing are Lipschitz continuous on any bounded set, 
        we may take $\gamma_0$ and $\ve_0$ independently from $h \in {\cal K}^{min}_a$.
              \QED


At the end of this subsection we discuss a large deviation principle of Freidlin-Wentzell type
for conditional measures. 
Using the upper estimate of the large deviation principle, 
we will prove that  contributions from a subset 
away from ${\cal L}( {\cal K}_a^{min} )$ is negligibly small
in the proof of the asymptotic expansions.
Note that we do not need the lower estimate, 
which probably does not hold true under our assumptions.

Let $\theta_{a}^{\ve}$ be a finite Borel measure on ${\cal W}$ 
which corresponds to the positive Watanabe distribution
$\delta_a (Y^{\ve}_1)$ via Sugita's theorem \cite{su}.
Since 
${\cal L}$ is quasi-surely defined, we can lift 
this measure  to a
measure  $ \mu^{\ve}_{a} := (\ve {\cal L})_* [\theta^{\ve}_{x,a}]$
on $G\Omega^B_{\alpha, 4m} ( {\mathbb R}^d)$.
In other words, 
$\mu^{\ve}_{a}$ is the law of the random variable $\ve {\bf W} ={\cal L}(\ve w)$
under 
$\theta_{a}^{\ve}$.

Set a rate function $I \colon G\Omega^B_{\alpha, 4m} ({\mathbb R}^n) \to [0, \infty]$ as follows;
\begin{align}
I ({\bf w}) 
= 
\begin{cases}
    \|h\|^2_{{\cal H}}/2 & (\mbox{if ${\bf w}= {\cal L}(h)$ for some $h \in {\cal K}_{a}$}), \\
    \infty &  (\mbox{otherwise}).
  \end{cases}
\nn
\end{align}
This rate function $I$ is actually good under {\bf (A1)}.
(We can prove this by using Lyons' continuity theorem and 
the goodness of the rate function 
for the usual Schilder-type large deviation on the rough path space.)

\begin{tm}\label{tm.ub_ldp}
Assume {\bf (A1)} and {\bf (A2)}.
Then, 
the family $\{ \mu^{\ve}_{a}\}_{\ve >0}$ 
of finite measures
satisfies the following upper bound of 
a large deviation principle 
on $G\Omega^B_{\alpha, 4m} ({\mathbb R}^d)$ as $\ve \searrow 0$ with a good rate function $I$, that is, 
for any closed set $A \subset G\Omega^B_{\alpha, 4m} ({\mathbb R}^d)$,
\begin{align}
 \limsup_{\ve \searrow 0 } \ve^2 \log \mu^{\ve}_{a} (A)
\le 
- \inf_{{\bf w} \in A} I  ({\bf w}).
\nn
\end{align}
\end{tm}

\Proof
This was  proved in Inahama \cite{in3} (under slightly stronger assumptions).
The keys of the proof 
are Kusuoka-Stroock's estimate (\ref{ks_bnd.ineq}) and 
the integration by parts formula (\ref{ipb1.eq}) for Watanabe distributions.
Hence, the same proof works under  {\bf (A1)} and {\bf (A2)}.
(The assumptions in
\cite{in3} are stronger for the lower bound of the LDP.)
\QED


\subsection{Basic differential geometry on ${\cal K}^{min}_a$}

In this subsection 
we recall some basic results on ${\cal K}^{min}_a$
and  fix notations.
All the ingredients in this subsection  can be found in Section 6.2, \cite{tw}.
We assume {\bf (A1)},  {\bf (B1)}, {\bf (B2)} in this subsection.

Let $\iota \colon {\cal K}^{min}_a \hookrightarrow {\cal H}$
denote the inclusion map.
The induced Riemannian metric $g$ on ${\cal K}^{min}_a$
is given by $g_h \la u,v \ra = \la   (\iota_*)_h u, (\iota_*)_h v\ra_{{\cal H}}$
for $u, v \in T_h {\cal K}^{min}_a$
and 
$h \in {\cal K}^{min}_a$.
The Riemannian volume measure is denoted by $\omega$.

Let ${\cal  A}_h \colon T_h  {\cal K}^{min}_a \times T_h  {\cal K}^{min}_a
\to  (T_h  {\cal K}^{min}_a)^{\bot}$
be the second fundamental form
defined by
\[
{\cal  A}_h \la u,v\ra = \pi^{\bot} [ \nabla_{(\iota_*)_h u} (\iota_*)_h v ]
\qquad
\mbox{ for $u, v \in T_h  {\cal K}^{min}_a$
and 
$h \in {\cal K}^{min}_a$.}
\]
Here, $\nabla$ denotes the flat connection on ${\cal H}$ 
and 
$\pi^{\bot}$ is a short hand for the orthogonal projection 
from ${\cal H}$ onto the orthogonal complement of 
$(\iota_*)_h [ T_h  {\cal K}^{min}_a]$.
(Precisely, $v$ should be extended to a vector field near $h$.)

Set ${\bf i}_h \in {\mathbb D}_{\infty} ( T_h  {\cal K}^{min}_a)$ by
\[
\la  u, {\bf i}_h \la w \ra \ra_{T_h  {\cal K}^{min}_a}
 = \la (\iota_*)_h u, w\ra_{{\cal H}}
\qquad \quad
 \mbox{ for $u \in T_h  {\cal K}^{min}_a$,
$h \in {\cal K}^{min}_a$.}
\]
Precisely,
the right hand side above is the stochastic extension of 
the bounded linear map
$\la (\iota_*)_h u,  \,\cdot\, \ra_{{\cal H}} \colon {\cal H} \to {\mathbb R}$.
However, 
we will see in the next section (Proposition \ref{pr.smth}) that
it in fact 
extends to a bounded linear map from ${\cal W}$ to ${\mathbb R}$.

Define ${\bf a}_h  
\in {\mathbb D}_{\infty} ( T^*_h ( {\cal K}^{min}_a) \otimes T^*_h ( {\cal K}^{min}_a))$
by
\[
{\bf a}_h (w)\la u,v\ra
= 
\la  {\cal  A}_h \la u,v\ra , w\ra_{{\cal H}}
\qquad 
 \mbox{ for $u, v \in T_h {\cal K}^{min}_a$,
$h \in {\cal K}^{min}_a$.}
\]
For the same reason as above, 
this also extends to a bounded linear map from ${\cal W}$
(see a local expression in (\ref{loc.a}) below).

An important remark is that
\begin{equation}\label{akk.eq}
{\bf a}_h (h)\la u,v\ra
= 
- \la  u, v \ra_{T_h  {\cal K}^{min}_a}
\qquad 
 \mbox{ for $u, v \in T_h {\cal K}^{min}_a$,
$h \in {\cal K}^{min}_a$.}
\end{equation}
Note that (\ref{akk.eq}) is non-trivial and 
is a consequence of the
 fact that ${\cal K}^{min}_a$ is a subset of a sphere in ${\cal H}$
centered at $0$.
It can be checked as follows.
First, the problem reduces to the case of two-dimensional sphere 
embedded in three (or higher) dimensional vector space
in a standard way.
Next, the case of the two-dimensional sphere
can be shown by straightforward computation.


We will write these quantities in a local coordinate chart
$(U; \theta^1, \ldots, \theta^{n'})$
of ${\cal K}^{min}_a$.
Set 
\begin{eqnarray*}
(\iota_*)_h  \Bigl(  \frac{\partial}{ \partial \theta^{i}}\Bigr)
=
 \Bigl(  \frac{\partial \iota}{ \partial \theta^{i}}\Bigr)_h,
 \quad
 g_{ij} (h) = 
 \Bigl\la
 \Bigl(  \frac{\partial \iota}{ \partial \theta^{i}}\Bigr)_h,
 \Bigl(  \frac{\partial \iota}{ \partial \theta^{j}}\Bigr)_h 
  \Bigr\ra_{{\cal H}}
  \\
 G(h) = ( g_{ij} (h))_{1\le i,j \le n'}
 \quad
 \mbox{and}
 \quad
 (e_i)_h :=  \{ G(h)^{-1/2} \}^{ij}  \Bigl(  \frac{\partial}{ \partial \theta^{j}}\Bigr)_h. \end{eqnarray*}
(Summation over repeated indices are omitted.)
By way of construction, 
$\{ (e_i)_h\}_{i=1}^{n'}$ is an orthonormal basis of $T_h  {\cal K}^{min}_a$.
Using this we can write down ${\bf i}_h$ and ${\bf a}_h$
as follows:
\begin{eqnarray}
{\bf i}_h (w)
&=&
 \la (\iota_*)_h   (e_i)_h , w \ra_{{\cal H}} (e_i)_h,
 \nn\\
{\bf a}_h (w)\la  (e_i)_h, (e_j)_h\ra
&=&
 G(h)^{-1/2}_{ik} 
{\bf a}_h (w)
\Bigl\la \Bigl(  \frac{\partial \iota}{ \partial \theta^{k}}\Bigr)_h,
 \Bigl(  \frac{\partial \iota}{ \partial \theta^{l}}\Bigr)_h    \Bigr\ra
  G(h)^{-1/2}_{lj} 
\nn \\
 &=&
  G(h)^{-1/2}_{ik}
 \Bigl[
 \Bigl\la
 \Bigl(  \frac{\partial^2  \iota}{ \partial \theta^{k} \partial \theta^{l}}\Bigr)_h,  w 
  \Bigr\ra_{{\cal H}}
\nn  \\
  &&\qquad
  -
    \Bigl\la
 \Bigl(  \frac{\partial^2  \iota}{ \partial \theta^{k} \partial \theta^{l}}\Bigr)_h,  
  (\iota_*)_h (e_m)_h
      \Bigr\ra_{{\cal H}}
        \la  (\iota_*)_h (e_m)_h ,w \ra_{{\cal H}}
             \Bigr]
   G(h)^{-1/2}_{lj}. 
   \label{loc.a}
      \end{eqnarray}


Define 
\begin{eqnarray*}
\Delta^{(h)} (w)  =  \delta_0 ({\bf i}_h (w))
\in 
\tilde{\mathbb D}_{- \infty}
\quad \mbox{and} \quad
D^{(h)} (w)  = \det[ - {\bf a}_h (w) ]   \in {\mathbb D}_{ \infty},
\end{eqnarray*}
where $\delta_0$ stands for the delta function at $0$ on $T_h  {\cal K}^{min}_a$.
Note that 
Malliavin covariance matrix of ${\bf i}_h (w)$ is the identity matrix.
In a local coordinate chart, $\Delta^{(h)}$ and $D^{(h)}$ can be written as follows:
\begin{equation*}
\Delta^{(h)} (w)  =
\sqrt{\det G(h)} \cdot 
\delta_0 \Bigl( \Bigl\{  \bigl\la
 \Bigl(  \frac{\partial \iota}{ \partial \theta^{i}}\Bigr)_h, w
  \bigr\ra_{{\cal H}}  \Bigr\}_{i=1}^{n'} \Bigr),
\end{equation*}
where $\delta_0$ stands for the delta function at $0$ on ${\mathbb R}^{n'}$,
and
\begin{eqnarray*}
D^{(h)} (w)  
&=&
\det G(h)^{-1} 
\\
&\times&
 \det \Bigl[
\Bigl\{
- \Bigl\la
 \Bigl(  \frac{\partial^2  \iota}{ \partial \theta^{i} \partial \theta^{j}}\Bigr)_h,  w 
  \Bigr\ra_{{\cal H}}
+
\Bigl\la
 \Bigl(  \frac{\partial^2  \iota}{ \partial \theta^{i} \partial \theta^{j}}\Bigr)_h,  
  (\iota_*)_h (e_m)_h
      \Bigr\ra_{{\cal H}}
        \la  (\iota_*)_h (e_m)_h ,w \ra_{{\cal H}}\Bigr\}_{i, j=1}^{n'}
\Bigr].
\end{eqnarray*}
Hence, 
\begin{eqnarray}
D^{(h)} (w)  \Delta^{(h)} (w) 
&=&
\frac{1}{\sqrt{\det G(h)} }
\nn\\
&&
\det \Bigl[
\Bigl\{
- \Bigl\la
 \Bigl(  \frac{\partial^2  \iota}{ \partial \theta^{i} \partial \theta^{j}}\Bigr)_h,  w 
  \Bigr\ra_{{\cal H}}
\Bigr\}_{i, j=1}^{n'}
\Bigr]
\cdot
\delta_0 \Bigl( \Bigl\{  \bigl\la
 \Bigl(  \frac{\partial \iota}{ \partial \theta^{i}}\Bigr)_h, w
  \bigr\ra_{{\cal H}}  \Bigr\}_{i=1}^{n'} \Bigr)
   \in \tilde{\mathbb D}_{- \infty}.
   \nn
     \end{eqnarray}
Since ${\cal K}^{min}_a \ni h 
\mapsto D^{(h)}  \Delta^{(h)}  \in \tilde{\mathbb D}_{- \infty}$
is continuous,
the $\tilde{\mathbb D}_{- \infty}$-valued  integration 
\[
\int_U D^{(h)}  \Delta^{(h)}  \omega (dh)
\]
is well-defined, where $\omega (dh) = \sqrt{\det G(h)} d\theta^1 \cdots d\theta^{n'}$
is the Riemannian volume measure.


\section{Some computations of skeleton ODE}
\label{sec.skltn}

In this section we compute the skeleton ODE under {\bf (A1)} and {\bf (B1)}.
We do not assume {\bf (A2)}, {\bf (B2)} or {\bf (B3)} for a while.
The main point is that if $h \in {\cal K}_a^{min}$ then
$\phi (h)$ satisfies (the configuration component of) a Hamiltonian ODE
for a naturally defined Hamiltonian.
The key in the proof of this fact is the Lagrange multiplier method.
This kind of argument has been known for a long time 
(See probabilistic literatures such as
Bismut \cite{bi}, Ben Arous \cite{bena2} 
or sub-Riemannian geometric literatures such as \cite{cc,mo,ri}.)
Our exposition basically follows Section 2.2, Rifford \cite{ri}.
However, we take a new look at 
this well-known argument  from a viewpoint of rough path theory
and slightly modify it. 
For example,  the Lagrange multiplier $q(h)$,
the key quantity in this argument,
can be understood in a rough path way.

First we give some formulae for the skeleton ODE (\ref{ode.def}).
In the sequel  $D$ also stands for the Fr\'echet derivative on ${\cal H}$.
Let $\phi_t (h)$ be the solution of the skeleton ODE driven by $h \in {\cal H}$
and set $\psi_t (h) =\Pi_{{\cal V}} \phi_t (h)$.
For simplicity, we set ${\bf V} = [ V_1, \ldots, V_r]$ which is a $(d \times r)$-matrix.

By straightforward computation, 
\begin{equation}
 D \phi_t (h)\la k  \ra=  J_t (h) \int_0^t  J_s (h)^{-1} {\bf V} (\phi_s (h)) dk_s
\qquad
\qquad
(k \in {\cal H}).
\label{D1phi.eq}
 \end{equation}
From this expression, $D \phi_t (h)$ and $D \psi_t (h)= \Pi_{{\cal V}} D \phi_t (h)$ 
extend to continuous linear maps 
from ${\cal W}$ to ${\mathbb R}^d$ and to ${\cal V}$, respectively. 
In a similar way, we have
\begin{eqnarray}
 D^2 \phi_t (h)\la k, \hat{k}  \ra
 &=& 
  J_t (h) \int_0^t  J_s (h)^{-1}
 \Bigl\{
 \nabla^2 {\bf V} (\phi_s (h)) \la  D \phi_s (h)\la k  \ra, D \phi_s (h)\la \hat{k}  \ra ,dh_s\ra
   \nn\\
   &&
   + \nabla {\bf V} (\phi_s (h)) \la  D \phi_s (h)\la k\ra  , d\hat{k}_s\ra 
    \nn\\
   &&   +\nabla {\bf V} (\phi_s (h)) \la  D \phi_s (h)\la \hat{k} \ra  , dk_s\ra    \Bigr\}
   \qquad\qquad
(k,k' \in {\cal H})
\label{D2phi.eq}
 \end{eqnarray}
and
$D ^2\psi_t (h)= \Pi_{{\cal V}} D^2 \phi_t (h)$.
Note that 
$| D^2 \phi_t (h)\la k, \hat{k}  \ra| \le C \|k\|_{{\cal H}}  \|\hat{k}\|_{{\cal W}}$
for some positive constant $C=C(h)$ 
which actually depends only on $\|h\|_{{\cal H}}$.


We apply the Lagrange multiplier method to $\| \,\cdot\, \|^2_{{\cal H}} /2$
under the constraint $\psi_1 =a$.
Due to non-degeneracy of $D \psi_1 (h)$, 
we have the following;
for any $h \in {\cal K}^{min}_a$, there exists a unique 
$q=q(h) \in {\cal V}^* \cong {\cal V}$ such that
the function 
\[
{\cal H}
\ni \hat{h}   \mapsto \frac12 \|\hat{h}\|^2_{{\cal H}}- \la q(h), \psi_1 (\hat{h}) -a\ra_{{\cal V} } 
\]
is stationary at $h$,
that is, 
\begin{equation}\label{Lag.def}
\la h, k \ra_{{\cal H} }   = \la q(h), D \psi_1 (h) \la k\ra \ra_{{\cal V} } 
\qquad \qquad
(k\in {\cal H},~h \in {\cal K}^{min}_a).
\end{equation}
Hence, 
$\la h, \,\cdot\, \ra_{{\cal H} }  \in {\cal W}^*$
for any $h \in {\cal K}^{min}_a$.
%
%

%
If we choose an orthonormal basis of ${\cal V}$ and set $n =\dim{\cal V} $,
then we can easily see
from the above relation and {\bf (A1)} that
 $q(h)$ has the following explicit expression in this coordinate system:
\begin{equation}
q(h)^i   = \sum_{j=1}^n  \{\sigma [\psi_1](h)^{-1}\}_{ij} D\psi_1^j (h) \la h\ra
\qquad
\qquad
(1 \le i \le n).
\label{qofh.eq}
\end{equation}
Therefore, if we assume {\bf (B2)} in addition,  then $h \to q(h)$ is smooth from ${\cal K}^{min}_a$ 
to ${\cal V}$.
(We will see later that $h \to q(h)$ extends to a continuous map from 
the geometric rough path space under {\bf (A1)} and {\bf (B1)}.)


We now introduce a Hamiltonian
and a Hamiltonian ODE which are naturally associated with the skeleton ODE.
Define 
$$
H (x,p) = \frac12  \sum_{i=1}^r \la p, V_i (x) \ra^2 
\qquad
\qquad
((x,p) \in {\mathbb R}^d \times {\mathbb R}^d\cong T^* {\mathbb R}^d),
$$
where the bracket denotes the pairing of a covector
and a tangent vector.
This is called a Hamiltonian function and
clearly smooth on ${\mathbb R}^d \times {\mathbb R}^d \cong T^* {\mathbb R}^d$.

The Hamiltonian ODE  is given by
\begin{eqnarray}
\frac{dx_t}{dt} &=& \frac{\partial H}{\partial p} (x_t, p_t) 
=
\sum_{i=1}^r  V_i (x_t )  \la  p_t, V_i (x_t )\ra,
\label{Ham1.eq}
\\
\frac{dp_t}{dt} &=& - \frac{\partial H}{\partial x} (x_t, p_t) 
=
- \sum_{i=1}^r  \la p_t, \nabla V_i (x_t ) \ra   \la  p_t, V_i (x_t )\ra.
\label{Ham2.eq}
\end{eqnarray}
Since we always assume $x_0 =x$, dependence on $x$ will be suppressed.
The solution with the initial condition $(x, p_0)$ is denoted by $(x_t (p_0),  p_t (p_0))$.

In the natural coordinate $x=(x^1, \ldots, x^d)$ of ${\mathbb R}^d$, we can write
$p = \sum_{j=1}^d p_j dx^j$
and 
$V_i (x) = \sum_{j=1}^d V_i^j (x)  (\partial / \partial x^j)$.
Moreover, 
$$\la  p, V_i (x )\ra = \sum_{j=1}^d p_j V_i^j (x) \in {\mathbb R}
\quad \mbox{and} \quad 
\la  p, \nabla V_i (x )\ra = \sum_{j,k=1}^d   p_j   
\frac{\partial V_i^j }{  \partial x^k}(x)  dx^k. $$
Therefore, the Hamiltonian ODE 
(\ref{Ham1.eq})--(\ref{Ham2.eq}) above can also be expressed as follows:
\begin{eqnarray*}
\frac{dx_t^k }{dt}   &=& \frac{\partial H}{\partial p_k} (x_t, p_t) 
=
\sum_{i=1}^r  V_i^k (x_t ) \sum_{j=1}^d p_{j,t} V_i^j (x_t), 
\\
\frac{dp_{k,t}}{dt} &=& - \frac{\partial H}{\partial x^k} (x_t, p_t) 
=
- \sum_{i=1}^r  
\sum_{j=1}^d   p_{j,t}   \frac{\partial V_i^j }{  \partial x^k} (x_t)
 \sum_{l=1}^d p_{l,t} V_i^l (x_t).
\end{eqnarray*}
We do not prove that a unique global solution exists for a given initial condition
because it is unnecessary for our purpose.


A remarkable fact is that the solution $\phi_t (h)$ of the skeleton ODE 
for $h \in {\cal K}^{min}_a$
can be expressed as the  ``configuration component" of the Hamiltonian ODE.

\begin{pr}\label{pr.hml}
Assume {\bf (A1)} and  {\bf (B1)}. Let $h \in {\cal K}_a^{min}$ 
and denote by $q(h)$  the Lagrange multiplier given in (\ref{Lag.def}). 
Then, 
a unique global solution of the Hamiltonian ODE with 
the initial condition $p_0 = \la q(h),  \Pi_{{\cal V}} J_1(h) \bullet \ra_{{\cal V}} 
\in T^*_x {\mathbb R}^d (\cong  {\mathbb R}^d)$ 
exists and
satisfies that $x_t (p_0) \equiv \phi_t (h)$.
Moreover,  $h_t^{i} \equiv \int_0^t   \la  p_s (p_0), V_i (x_s (p_0) )\ra ds ~(1 \le i \le r)$.
\end{pr}

\Proof
For simplicity of notations
we assume that ${\cal V} = {\mathbb R}^n \times \{ {\bf 0}_{d-n} \}$, where $n =\dim {\cal V}$.
(Without loss of generality we may do so.)
Then, $\Pi_{{\cal V}} =[ {\rm Id}_n | {\bf 0}_{n \times (d-n)} ]$.
$J$ and $J^{-1}$ are $d \times d$ matrices
and ${\bf V}$ is a $d \times r$ matrix.
$q(h)$ and $p$ are regarded as an $n$-dimensional 
and 
a $d$-dimensional row vector, respectively.
%
Since $H$ is smooth, the uniqueness is obvious 
if a solution exists.
%

From (\ref{D1phi.eq}) and (\ref{Lag.def}) we have 
\[
\int_0^1 \la h_s^{\prime}, k_s^{\prime} \ra ds
=
\int_0^1
q(h) \Pi_{{\cal V}}  J_1(h)  J_s^{-1} (h) {\bf V} ( \phi_s (h) ) k_s^{\prime }ds. 
\]
Since $k \in {\cal H}$ is arbitrary,  we can easily see that 
\[
(h^1_t, \ldots,  h^r_t) \equiv 
\int_0^t q(h) \Pi_{{\cal V}}  J_1(h)  J_s^{-1} (h) {\bf V} ( \phi_s (h) ) ds.
\]
Define $x_t = \phi_t (h)$
and 
$p_t = q(h) \Pi_{{\cal V}}  J_1(h)  J_t^{-1} (h)$.
Clearly, 
$p_0 =q(h) \Pi_{{\cal V}}  J_1(h)$ and $(h_t^1, \ldots,  h^r_t) \equiv \int_0^t p_s {\bf V} ( \phi_s (h) ) ds$.
The skeleton ODE can be rewritten as
\[
d\phi_t (h) = \sum_{i=1}^r V_i ( \phi_t (h) )  dh^i_t 
= 
 \sum_{i=1}^r V_i ( \phi_t (h) )  [p_t  V_i ( \phi_t (h) )] dt,
\]
which means that (\ref{Ham1.eq}) is satisfied.
From the ODE for $J^{-1}_t =K_t$, we can see that
\begin{eqnarray*}
p_t^{\prime} 
&=&
q(h) \Pi_{{\cal V}}  J_1(h)  
\Bigl[  
 - J_t^{-1} (h) \sum_{i=1}^r  \nabla V_i (\phi_t (h)  )  ( h^i_t)^{\prime}
\Bigr]
\\
&=&
- p_t
\sum_{i=1}^r  \nabla V_i (\phi_t (h)  ) [p_t  V_i ( \phi_t (h) )].
\end{eqnarray*}
Thus, (\ref{Ham2.eq}) is also satisfied.
\QED

\begin{re}\normalfont\label{re.rglr}
In Proposition \ref{pr.hml} above, 
the Hamiltonian $H$ is of course constant 
along the trajectory of $t \mapsto (x_t (p_0), p_t (p_0))$.
From the explicit form of $h$ we can easily see that
the constant is $d_a^2 /2 = \|h\|_{{\cal H}}^2 /2 >0$.
Furthermore, this implies that 
the smooth path
$t \mapsto x_t (p_0) = \phi_t (h)$ is regular, that is,
its velocity vector never vanishes.
Indeed, 
\[
\la  p_t (p_0),  \phi^{\prime}_t (h)\ra
=
\Bigl\la  p_t (p_0),   \sum_{i=1}^r V_i ( \phi_t (h) )  \la p_t (p_0),  V_i ( \phi_t (h) )\ra 
\Bigr\ra
=
2H ( x_t (p_0), p_t (p_0)) 
\equiv
d^2_a>0
 \]
 for any $t$.
%
\end{re}

Thanks to the above proposition, 
we know that if $h \in {\cal K}^{min}_a$, then both $h$ and $\phi(h)$
are smooth in $t$.
This fact is highly non-trivial 
since a generic element of ${\cal H}$ is not even of $C^1$.
Moreover, we have the following proposition:
\begin{pr}\label{pr.smth}
Assume {\bf (A1)},  {\bf (B1)} and  {\bf (B2)}.
Then, the mappings
$(t, h) \mapsto \phi_t (h)$ and $(t, h) \mapsto h_t$
are smooth from $[0,1] \times {\cal K}^{min}_a$ to ${\mathbb R}^d$ and to ${\mathbb R}^r$,
respectively.
\end{pr}

\Proof
Assume that 
 Hamiltonian ODE (\ref{Ham1.eq})--(\ref{Ham2.eq}) has a unique global solution
for an initial condition $(x, \hat{p}_0)$. 
Then, a standard cut-off technique  and smoothness of $H$
shows that 
there is an open neighborhood $U$ of $\hat{p}_0$ such that
{\rm (i)} a unique global solution exists for any $(x,p_0)$ with $p_0 \in U$
and 
{\rm (ii)} the mapping $[0,1] \times U \ni (t,p_0) \mapsto (x_t (p_0), p_t (p_0)) 
\in {\mathbb R}^{d} \times {\mathbb R}^{d}$
is smooth.

Since ${\cal K}^{min}_a \ni h \mapsto p_0=\la q(h),  
\Pi_{{\cal V}} J_1(h) \bullet \ra_{{\cal V}} $ is smooth 
under {\bf (B2)},
we see that $(t, h) \mapsto x_t(p_0) =\phi_t (h)$  is also smooth.
Smoothness of $(t, h) \mapsto h_t$ is immediate if we recall
\[
\frac{d }{dt} 
(h^1_t, \ldots,  h^r_t) \equiv 
q(h) \Pi_{{\cal V}}  J_1(h)  J_t^{-1} (h) {\bf V} ( \phi_t (h) ).
\]
This completes the proof.
\QED


The continuity of 
the natural lift map ${\cal L}|_{{\cal H}} \colon {\cal H} \hookrightarrow G\Omega_{\alpha, 4m}^{B} ({\mathbb R}^r)$
means that the rough path topology, restricted on ${\cal L}({\cal H})$, 
is weaker than ${\cal H}$-topology.
The next proposition claims that
on ${\cal L} ({\cal K}^{min}_a )$ the two topologies are in fact the same.
\begin{pr}
\label{pr.top=}
Assume {\bf (A1)} and  {\bf (B1)}.
Let $h_1, h_2, \ldots, h_{\infty} \in {\cal K}^{min}_a$.
Then, $\lim_{j \to \infty} h_j = h_{\infty}$ in ${\cal H}$
if and only if 
$\lim_{j \to \infty}  {\cal L} (h_j) = {\cal L}(h_{\infty})$ in $G\Omega_{\alpha, 4m}^{B} ({\mathbb R}^r)$.
\end{pr}

\Proof
By setting $({\bf w}, \lambda) =( {\cal L} (h), 0)$ in 
(\ref{3arrows.eq}),
we see that $h \mapsto (h, \phi(h), J(h), J^{-1}(h))$ extends 
to a continuous map from $G\Omega_{\alpha, 4m}^{B} ({\mathbb R}^r)$
to 
$G\Omega_{\alpha'}^{H} ({\mathbb R}^r \oplus {\mathbb R}^d \oplus  {\rm Mat}(d,d) \oplus  {\rm Mat}(d,d))$
for 
any $\alpha' \in (1/3, \alpha - 1/(4m))$.
From this and (\ref{D1phi.eq}) we see that 
$$h \mapsto D \psi_1 (h)\la h \ra = \Pi_{{\cal V}}J_1 (h) \int_0^1  J_s (h)^{-1} {\bf V} (\phi_s (h)) dh_s$$
extends to a continuous map from $G\Omega_{\alpha, 4m}^{B} ({\mathbb R}^r)$
since the right hand side 
can be regarded as a rough path integral along  the natural lift of 
$(h, \phi(h), J(h), J^{-1}(h) )$.
Also, the deterministic Malliavin covariance matrix
$h \mapsto \sigma [\psi_1] (h)$
extends to a continuous map from $G\Omega_{\alpha, 4m}^{B} ({\mathbb R}^r)$.
Hence, from the explicit expression (\ref{qofh.eq}),
 we see that 
${\cal K}^{min}_a \ni h \mapsto q (h)$
is  continuous with respect to the rough path topology.
Therefore, 
as $h \in  {\cal K}^{min}_a$ varies continuously with respect to the rough path topology,
\[
\frac{d}{dt} (h^1_t, \ldots,  h^r_t) \equiv 
 q(h) \Pi_{{\cal V}}  J_1(h)  J_t^{-1} (h) {\bf V} ( \phi_t (h) )
\]
varies continuously 
with respect to $L^2 ([0,1], {\mathbb R}^r)$-topology.
Thus, 
$\lim_{j \to \infty}  {\cal L} (h_j) = {\cal L}(h_{\infty})$ in $G\Omega_{\alpha, 4m}^{B} ({\mathbb R}^r)$
implies 
$\lim_{j \to \infty} \| h_j - h_{\infty} \|_{{\cal H}} =0$.
The converse implication is obvious. 
\QED

\begin{re}\normalfont\label{re.top=}
Since the rate function in a Schilder-type large 
deviation principle on the rough path space is good, 
${\cal K}^{min}_a$ is compact with respect to $G\Omega_{\alpha, 4m}^{B} ({\mathbb R}^r)$-topology
under {\bf (A1)}.
If we assume  {\bf (B1)} in addition, the above proposition implies that
it is also compact with respect to ${\cal H}$-topology.
From this viewpoint, the compactness assumption in {\bf (B2)} is not strong. 
\end{re}


The following lemma is essentially a rough path version of  Lemma 6.2, \cite{tw}.
We use $(\alpha, 4m)$-Besov rough path topology 
 instead of the norm $\{ |w_1|^2 + \int_0^1 |w_s|^2 ds \}^{1/2}$.
In the statement of this lemma, a  ``neighborhood" means a neighborhood
with respect to ${\cal H}$-topology.
Due to Proposition \ref{pr.top=},
the assertion {\rm (i)} clearly holds if $\gamma$ is small enough
since $\{ U_{h, \gamma} \mid \gamma >0\}$ forms a fundamental system of neighborhood
of ${\cal L}(h)$
with respect to $(\alpha, 4m)$-Besov rough path topology.
\begin{lm}\label{lm.tw_62}
Assume {\bf (A1)},   {\bf (B1)} and  {\bf (B2)}.
For any $h \in {\cal K}^{min}_a$, there exists a coordinate neighborhood 
$O \subset {\cal K}^{min}_a$
of 
 $h$ which satisfies the following property:
 For any subneighborhood $O^{\prime}$ of $O$, there exists $\gamma >0$ such that
\\
{\rm (i)}~$\{ k \in  {\cal K}^{min}_a \mid{\cal L}(k)\in U_{h, \gamma}\} \subset O^{\prime}$,
\\
{\rm (ii)}~$\int_{O^{\prime}} D^{(k)} \Delta^{(k)} \omega(dk) =1$ on 
$\{ w \in  {\cal W} \mid   {\cal L}(w)\in U_{h, \gamma} \}$,
that is, 
for any $F \in \tilde{\mathbb D}_{\infty}$ such that $F \cdot {\bf 1}_{ \{  {\cal L}(w)\in U_{h, \gamma} \}  } =F$ a.s.,
it holds that
${\mathbb E} [F \int_{O^{\prime}} D^{(k)} \Delta^{(k)} \omega(dk) ] = {\mathbb E} [F]$.
\end{lm}

\Proof 
The proof here is essentially the same as the proof of Lemma 6.2, \cite{tw}.
We use Proposition \ref{pr.smth} and the inverse function theorem.
We denote by ${\cal W}^B_{\alpha,4m}$
the Banach space of all the elements in ${\cal W}$  with finite $(\alpha,4m)$-Besov norm.
The injection 
${\cal K}_a^{min} \hookrightarrow {\cal H}$ is
denoted by $\iota$, which will be sometimes omitted 
when the notation gets too heavy.

Take $h \in {\cal K}_a^{min}$ and a coordinate neighborhood $O$
of $h$ arbitrarily.
Let $Q$ be an open subset of ${\mathbb R}^{n'}$ such that
$Q \ni (\theta^1, \ldots, \theta^{n'})
  \mapsto k(\theta) =k (\theta^1, \ldots, \theta^{n'})\in O$ is diffeomorphic
(namely, $k\colon Q \to O$ is a chart map).

We consider 
\[
\frac{\partial \iota}{\partial \theta^{i} } \Big|_{k(\theta)}
=
\Bigl\{ 
\frac{\partial }{\partial \theta^{i} }  k(\theta)_t
\Bigr\}_{0 \le t \le 1}
\qquad
\mbox{and}
\qquad
\frac{\partial^2 \iota}{\partial \theta^{i} \partial \theta^{j}} \Big|_{k(\theta)}
=
\Bigl\{ 
\frac{\partial^2 }{\partial \theta^{i} \partial \theta^{j}}  k(\theta)_t
\Bigr\}_{0 \le t \le 1}.
\]
An integration by parts on $[0,1]$ immediately yields 
\begin{eqnarray*}
\Bigl\langle 
\frac{\partial \iota}{\partial \theta^{i} } \Big|_{k(\theta)}, w
\Bigr\rangle_{{\cal H}}
&=&
\Bigl\la
w_1,  \frac{\partial }{\partial \theta^{i} }  k(\theta)_1^{\prime} 
\Bigr\ra_{{\mathbb R}^r }
- 
\int_0^1  
\Bigl\la w_t,  \frac{\partial }{\partial \theta^{i} }  k(\theta)_t^{\prime\prime}
 \Bigr\ra_{{\mathbb R}^r } dt,
 \\
 \Bigl\langle 
\frac{\partial^2 \iota}{\partial \theta^{i} \partial \theta^{j}} \Big|_{k(\theta)}, w
\Bigr\rangle_{{\cal H}}
&=&
 \Bigl\langle
 w_1, \frac{\partial^2 }{\partial \theta^{i} \partial \theta^{j}}  k(\theta)_1^{\prime}
\Bigr\ra_{{\mathbb R}^r }
- 
 \int_0^1 \Bigl\langle  w_t,  \frac{\partial^2 }{\partial \theta^{i} \partial \theta^{j} }  k(\theta)_t^{\prime\prime}
\Bigr\ra_{{\mathbb R}^r }
 dt.
 \end{eqnarray*}
From Proposition  \ref{pr.smth},
these two functions on the left hand sides above
are defined for all $(\theta, w) \in Q \times {\cal W}^B_{\alpha,4m}$
and continuous with respect to the product topology.
(In particular, 
for each fixed $\theta$
these two linear functionals on ${\cal H}$ continuously extends to ones on ${\cal W}^B_{\alpha,4m}$ or ${\cal W}$.)
Moreover, 
they are smooth in $\theta$ for any fixed $w$ and
their derivatives are all continuous in $(\theta, w)$, too.


Set 
\begin{eqnarray*}
f( \theta, w) 
&=&
- \Bigl[
\Bigl\langle 
\frac{\partial \iota}{\partial \theta^{i} } \Big|_{k(\theta)}, w
\Bigr\rangle_{{\cal H}}
\Bigr]_{ 1 \le i \le n'},
\\
J_f( \theta, w) 
&=&
- \Bigl[
\Bigl\langle 
\frac{\partial^2 \iota}{\partial \theta^{i}  \partial \theta^{j}} \Big|_{k(\theta)}, w
\Bigr\rangle_{{\cal H}}
\Bigr]_{ 1 \le i, j \le n'},
\end{eqnarray*}
where the latter is the Jacobian matrix of the former with respect to the $\theta$-variable.
Then, we can easily check that
$
f(\theta, k(\theta)) =0
$
and 
$
J_f (\theta, k(\theta))_{ij} = g_{ij} ( k(\theta))
$,
where $(g_{ij} )$ is the Riemannian metric tensor on ${\cal K}_a^{min}$.
(We omit the proof because it is a routine. See pp. 215--216, \cite{tw}.
However, we note that the fact that ${\cal K}_a^{min}$ is a subset of a sphere is crucial here.)
For each fixed $\theta$, 
$f$ is an ${\mathbb R}^{n'}$-valued linear functional 
which is non-degenerate in the sense of Malliavin
since its Malliavin covariance matrix is given by $(g_{ij} )$. 
Also, $J_f$ is linear in $w$  and $\det J_f$ is just a polynomial in $w$.


For each fixed $w$,
we use the inverse function theorem for $f(\,\cdot\,, w)$ 
and carefully  keep track of the dependence on $w$.
We will write $h = k (\theta_0)$,
$Q_{\rho} (\theta_0)= \{ \theta \mid  |\theta -\theta_0 | <\rho  \} \subset Q$
and 
$B_{\gamma} (h)
 = \{ w \mid \| w-h\|_{\alpha,4m -B}  < \gamma \}$ for $\rho, \gamma >0$.

Then, we can show the existence of $\rho' >0$ and $\gamma'  >0$ such that
\\
\\
{\rm (a)}~for each $w \in B_{\gamma' } (h)$, $Q_{\rho'} (\theta_0) \ni \theta \mapsto
f(\theta,w)\in
 f (Q_{\rho'} (\theta_0),w)$ is diffeomorphic and
\\
{\rm (b)}~for each $\rho \in (0,\rho')$, there exists $\gamma \in (0,\gamma')$ such that
 $\cap \{  f (Q_{\rho} (\theta_0),w) \mid w \in  B_{\gamma} (h)\}$
 contains $0 \in {\mathbb R}^{n'}$ as an interior point
 and the assertion {\rm (i)} holds. 
\\
\\
(Though it is not difficult,  it is not so obvious, either. 
Find a nice textbook on calculus and modify a proof in it.
Note also that  the implicit function theorem  is not used here.)
We write $O_{\rho} = \{ k ( \theta) \mid \theta \in Q_{\rho} (\theta_0)\}$,
which is the subneighborhood $O'$ in the statement of the lemma.

Let $\chi \colon {\mathbb R}^{n'} \to [0,\infty)$ be a smooth, radial function 
with compact support
such that $\int_{ {\mathbb R}^{n'}} \chi =1$.
It is well-known that $\chi_{\kappa} := \kappa^{-n'} \chi ( \,\cdot\, /\kappa ) \to \delta_0$ as $\kappa \searrow 0$
in ${\cal S}^{\prime} ({\mathbb R}^{n'})$.
If $\kappa >0$ is sufficiently small, 
\[
\int_{Q_{\rho} (\theta_0) }   \det  J_f(\theta, w)  \chi_{\kappa} (f(\theta, w))  d\theta =1
\qquad
\mbox{for all $w \in B_{\gamma} (h)$.}
\]
If $F$ is supported in $\{ w \mid   {\cal L}(w)\in U_{h, \gamma} \} \subset B_{\gamma} (h)$,
\[
{\mathbb E} \Bigl[
F
\int_{Q_{\rho} (\theta_0) }   \det  J_f(\theta, \,\cdot\,)  \chi_{\kappa} (f(\theta, \,\cdot\,))  d\theta
\Bigr]
=
{\mathbb E} [F].
\]
There exists constant $l>0$ independent of $\theta$ such that
$ \chi_{\kappa} (f(\theta, \,\cdot\,)) \to \delta_0 (f(\theta, \,\cdot\,) ) = \delta_0 (-f(\theta, \,\cdot\,) ) $
in ${\mathbb D}_{p', -l}$-norm
for any $p' \in (1, \infty)$ uniformly in $\theta$.
Since $F \in {\mathbb D}_{p, l}$ for some $p =p(l) \in (1, \infty)$,
we have  
\[
{\mathbb E} \Bigl[
F
\int_{Q_{\rho} (\theta_0) }   \det  J_f(\theta, \,\cdot\,)  \delta_0 (- f(\theta, \,\cdot\,))  d\theta
\Bigr]
=
{\mathbb E} [F].
\]
Since $\omega(dk) = \sqrt{ \det G(k(\theta)) } d\theta$,
where $G = (g_{ij})$ is the Riemannian metric tensor,
the left hand side is equal to
\[
{\mathbb E} \Bigl[
F \int_{O_{\rho}  }  D^{(k)} \Delta^{(k)}\omega(dk) \Bigr].
\]
Thus, we have shown {\rm (ii)}.
\QED

\begin{re}\normalfont\label{re.tw_62}
(1)~
By carefully examining the proof,
we can slightly strengthen the assertion {\rm (ii)}, Lemma \ref{lm.tw_62}
as follows:
For any $\ve \in (0,1]$,
\[
\mbox{
$\int_{O^{\prime}} D^{(k)} (\ve w) \Delta^{(k)} (\ve w)  \omega(dk) =1$ on 
$\{ w  \mid   {\cal L}(\ve w)\in U_{h, \gamma} \}$
}
\]
for the same $\gamma >0$ as in Lemma \ref{lm.tw_62}.
The proof is essentially the same, but we should note that 
$f(\theta, \ve w)$, $J_f(\theta, \ve w)$, etc. are well-defined 
since $ f(\theta, \,\cdot\,), J_f(\theta, \,\cdot\,)$, etc. are continuous maps from ${\cal W}$.
\\
(2)~It is almost obvious that if {\rm (i)}, {\rm (ii)} in Lemma \ref{lm.tw_62}
hold for some $\gamma >0$,
then they still hold for any $\hat\gamma \in (0, \gamma)$. 
In other words, we may replace $\gamma$ in Lemma \ref{lm.tw_62} by any smaller positive constant. 
\end{re}



\section{Exponential integrability lemmas}

In this section
 we will see that the positivity of Hessian  assumed in {\bf (B3)} implies 
exponential integrability of a corresponding quadratic Wiener functional.
In this section we assume {\bf (A1)} and {\bf (B1)}--{\bf (B3)}.


Let $(- \tau_0, \tau_0) \ni \tau \mapsto c(\tau) \in {\cal K}_a$
 be a smooth curve in ${\cal K}_a$
such that $c(0) = h \in {\cal K}_a^{min}$ and $
0\neq k: =c^{\prime} (0)  
\in {\cal H}_0 (h)$,
where we set 
${\cal H}_0 (h)=
T_h {\cal K}_a  \cap (T_h {\cal K}^{min}_a )^{\bot}$ 
as in {\bf (B3)}.
Then, a straightforward calculation shows that
\begin{align}
\lefteqn{
\frac{d^2}{d \tau^2} \Big|_{\tau=0} \frac{ \|   c(\tau) \|^2_{{\cal H}} }{2}
=
\frac{d^2}{d \tau^2} \Big|_{\tau=0}
\Bigl(  \frac{ \|   c(\tau) \|^2_{{\cal H}} }{2}  - \la q(h), \psi_1(c(\tau)) -a \ra  \Bigr)
}
\nn\\
&
=
\|   c^{\prime}(0) \|^2_{{\cal H}}+ \la  c^{\prime\prime}(0),  c(0)  \ra_{{\cal H}}
 - \bigl\la q(h) , D\psi_1(c(0)) \la c^{\prime\prime}(0) \ra \bigr\ra 
 - \bigl\la q(h) , D^2\psi_1(c(0) ) \la c^{\prime}(0), c^{\prime}(0) \ra \bigr\ra 
    \nn\\
&
=
\|   k \|^2_{{\cal H}}
-  \bigl\la q(h), D^2\psi_1(h) \la k, k \ra \bigr\ra.
   \label{hess.cmp1}
    \end{align}
The cancellation above was due to the Lagrange multiplier method (\ref{Lag.def}).

It is known that the symmetric bounded bilinear form 
$\bigl\la q(h), D^2\psi_1(h) \la \bullet, \star \ra \bigr\ra$
is Hilbert-Schmidt
and so is $\bigl\la q(h), D^2\psi_1(h) \la \pi^h \bullet, \pi^h\star \ra \bigr\ra$,
where $\pi^h \colon {\cal H} \to {\cal H}_0 (h)$ is the orthogonal projection.
(An explicit form of $\pi^h$ will be given in (\ref{pih.def}) below.)
As a result, the spectra of their corresponding  symmetric operators are discrete except at $0$.
Consequently, the condition {\bf (B3)} implies that
\begin{eqnarray}
C(h)
:=
\sup
\Bigl\{
\frac12 \bigl\la q(h), D^2\psi_1(h) \la \pi^h v, \pi^h v' \ra  \bigr\ra_{{\cal V}}  ~\Big|~
 v,  v' \in {\cal H}, \| v \|_{{\cal H}}  =\| v' \|_{{\cal H}}  =1
\Bigr\}
<\frac 12.
\label{hess.cmp2}
 \end{eqnarray}
From the explicit expressions of $D^2\psi_1(h)$, $q(h)$ and $\pi^h$
in (\ref{D2phi.eq}), (\ref{qofh.eq}) and (\ref{pih.def}),   
we can easily see that 
\begin{equation}\label{up12.ineq}
{\cal K}_a^{min} \ni h
\mapsto 
\frac12 \bigl\la q(h), D^2\psi_1(h) \la \pi^h \bullet, \pi^h\star \ra \bigr\ra
\end{equation}
is continuous with respect to 
the topology of bounded bilinear forms.
This implies the continuity of $h \mapsto C(h)$ and 
\begin{equation}\label{cmax1/2.ineq}
C_{max} :=\max\{C(h) \mid h \in {\cal K}_a^{min} \}< \frac12.
\end{equation}

Therefore, if $\Xi_h \in {\cal C}_2$ corresponds to the symmetric Hilbert-Schmidt 
bilinear form in (\ref{up12.ineq}), namely the bilinear form is equal to $(1/2) D^2\Xi_h$, 
then 
${\mathbb E} [ e^{c_1 \Xi_h  } ] \le M <\infty$
for certain constants $c_1 >1$ and $M >0$  independent of $h \in {\cal K}_a^{min} $.
(We denote by ${\cal C}_i$  the $i$th order homogeneous Wiener chaos.)


The following lemma is Lemma 6.1, \cite{tw},
which 
states that the quadratic Wiener functionals that appear as the second term of 
the Taylor-like  expansion of It\^o map is exponentially integrable with respect to the 
conditional Gaussian measure.  
Though $g_2$ can be written as a iterated stochastic  integral of second order,
we do not use it.
Our proof uses rough path theory and  
can be found in \cite{in2} in  a more detailed way.
For the definition of ${\bf i}_h$ and $g_i~(i=0,1,2,\ldots)$
see (\ref{ih.def}) and  (\ref{asympX.2}), respectively.
\begin{lm}\label{lm.tw61}
Assume {\bf (A1)} and {\bf (B1)}--{\bf (B3)}.
Then, there exists a constant $c_1 >1$ which is 
 independent of $h \in {\cal K}^{min}_a$
 and satisfies that
\begin{equation}\label{tw_eq_6.4}
\sup_{h \in {\cal K}^{min}_a}
{\mathbb E} \Bigl[
\exp \bigl( c_1 \langle q(h), g_2 (w;h)\rangle_{{\cal V}} \bigr)
\delta_0 ( g_1 (w;h),  {\bf i}_h \la w\ra)
\Bigr]
<\infty.
\end{equation}
Here, $\delta_0$ is the Dirac delta function on ${\cal V} \times T_h {\cal K}^{min}_a\cong {\mathbb R}^{n+ n'}$.
\end{lm}

\Proof
Let $h \in {\cal K}^{min}_a$ and let
$\{e_i (h)\}_{i=1}^{n'}$ be an orthonormal basis of $T_h {\cal K}^{min}_a \subset {\cal H}$.
Due to Proposition \ref{pr.smth}, $\la e_i (h), \,\cdot\, \ra_{{\cal H}} \in {\cal H}^*$ naturally extends to 
an element in ${\cal W}^*$.
Recall that $g_1 (\,\cdot\, ; h ) \in {\cal W}^*$ is the continuous 
extension of $D\psi_1 (h) \in {\cal H}^*$. 
Then,
$$
{\cal H}_0 (h)
=
\ker g_1 (\,\cdot\, ; h )  \cap \bigl( \cap_{i=1}^{n'}  \ker  \la e_i (h), \,\cdot\,\ra_{{\cal H}} \bigr)
=
\ker g_1 (\,\cdot\, ; h )  \cap \ker {\bf i}_h,
$$
which may depend on $h$, but not on the choice of the orthonormal basis.
It is easy to see that 
${\cal H}_0 (h)^{\bot}$ is a finite dimensional subspace of dimension $n+ n'$.
Denote by ${\cal W}_0 (h)$  the closure of 
${\cal H}_0 (h)$ with respect to the topology of ${\cal W}$.
The orthogonal projection $\pi^h \colon {\cal H} \to {\cal H}_0 (h)$
naturally extends to 
a continuous 
projection $\pi^h \colon {\cal W} \to {\cal W}_0 (h)$, which we denote by the same symbol.

The finite Borel measure corresponding to 
the positive Watanabe distribution $\delta_0 ( g_1 (\,\cdot\, ; h ), {\bf i}_h )$
via Sugita's theorem \cite{su}
is 
$(2\pi)^{-(n+n')/2} \{ \det \sigma [\psi_1] (h) \}^{-1/2} \cdot \pi_{*}^h \mu$, 
which is a constant multiple of a Gaussian probability measure $\pi_{*}^h \mu$ 
supported on ${\cal W}_0 (h)$.
It is clear that $( {\cal W}_0 (h), {\cal H}_0 (h), \pi_{*}^h \mu)$ is an abstract Wiener space.


Let us now identify ${\cal V} \cong {\mathbb R}^n$ by choosing an orthonormal basis of ${\cal V}$.
(The choice is independent of $h$.)
Then we have explicitly that
\begin{eqnarray}
\pi^h w &=& w - \sum_{l, l'=1}^n   \{\sigma [\psi_1] (h)^{-1}\}_{ll'}   g_1 (w ; h )^l 
\cdot {}^{\sharp} g_1 (\,\cdot\, ; h )^{l'}
-\sum_{i=1}^{n'} \la e_i (h),w \ra e_i (h)
\nn\\
&=&
 w- \sum_{l, l'=1}^n   \{\sigma [\psi_1] (h)^{-1}\}_{ll'}   g_1 (w ; h )^l 
 \cdot {}^{\sharp} g_1 (\,\cdot\, ; h )^{l'}
- {\bf i}_h \la w\ra
\label{pih.def}
\end{eqnarray}
and we set $(\pi^h)^{\bot} w = w -\pi^h w \in {\cal H}$.
Here, 
$g_1 (w ; h )^l$ is the $l$th component of $g_1 (w ; h )= D\psi_1 (h) \la w\ra$,
$\{ \sigma [\psi_1] (h)^{-1} \}_{ll'}$ is the $(l,l')$th component of the inverse of 
 the deterministic covariance matrix 
$\sigma [\psi_1] (h)$ at $h$,
${}^{\sharp} g_1 (\,\cdot\, ; h )^{l}$ is the unique element in ${\cal H}$
that corresponds to $g_1 (\,\cdot\, ; h )^{l} \in {\cal W}^* \subset {\cal H}^*$ via the Riesz isometry.
Note that 
${}^{\sharp} g_1 (\,\cdot\, ; h )^{l}$ and $e_i (h)$ are orthogonal for any $l$ and $i$.

By abusing notations we write 
$\pi^h {\bf W} = {\cal L} (\pi^h w)=
\lim_{l\to \infty} {\cal L} ((\pi^h w)(l))$.
Here, 
$w(l)$ denotes the $l$th dyadic 
piecewise linear approximation of $w$ and
${\cal L}$ denotes the rough path lift map defined by  those dyadic approximations.
By the explicit formula (\ref{pih.def}) for $\pi^h w$,
we can easily see that $\pi^h {\bf W}$ is well-defined a.s. with respect to $\mu$
and in fact equal to the Young translation of ${\bf W}$ 
by 
$- \sum_{l, l'=1}^n  \{ \sigma [\psi_1] (h)^{-1} \}_{ll'}   \cdot g_1 (w ; h )^l 
\cdot {}^{\sharp} g_1 (\,\cdot\, ; h )^{l'}
-{\bf i}_h \la w\ra$.

It is shown in \cite{in2} that $g_2$ actually has the following form:
For some continuous map $\hat{g}_2 \colon G\Omega^B_{\alpha, 4m} ({\mathbb R}^d) \times {\cal H} \to {\cal V}$
of  ``quadratic order,"
\begin{eqnarray*}
g_2(w; h) &=& \hat{g}_2 ( {\bf W}, h)
\nn\\
&=&
\lim_{l \to \infty}
\Bigl\{
\frac12 D^2 \psi_1 (h)\la  w(l),w(l)   \ra 
+ 
\Pi_{{\cal V}}  J_1 (h) \int_0^1  J_s (h)^{-1} V_0 (\phi_s (h)) ds
\Bigr\}.
\end{eqnarray*}

Then, 
$g_2(\pi^h w; h) = \hat{g}_2 ( \pi^h{\bf W}, h)$ is defined 
with respect to $\mu$.
This in turn implies that 
$g_2(w; h) = \hat{g}_2 ( {\bf W}, h)$ is well-defined with respect to $\pi_{*}^h \mu$.
Hence, if 
\begin{equation}\label{expint1.ineq}
\int_{{\cal W}_0 (h)} 
\exp ( c_1 \langle q(h), g_2 (w;h)\rangle_{{\cal V}})  \pi_{*}^h \mu (dw)
=
\int_{{\cal W}}
 \exp ( c_1 \langle q(h), g_2 (\pi^h w;h)\rangle_{{\cal V}})  \mu (dw)
\le M
\end{equation}
holds for some $M >0$,
then the proof of the lemma is done.
(Since $g_2 (w;h) = \hat{g}_2 ( {\bf W}, h)$ is $\infty$-quasi continuous in $w$, 
the left hand side above is a constant multiple of the integral in (\ref{tw_eq_6.4}).)

Now we show (\ref{expint1.ineq}).
It is straightforward to check that 
${\mathbb E}[D_k \langle q(h), g_2 (\pi^h w;h)\rangle_{{\cal V}} ]= 0$ for all $k \in {\cal H}$
and
$D^3 \langle q(h), g_2 (\pi^h w;h)\rangle_{{\cal V}} =0$.
This shows that 
$$ \langle q(h), g_2 (\pi^h w;h)\rangle_{{\cal V}}  
- {\mathbb E} [\langle q(h), g_2 (\pi^h w;h)\rangle_{{\cal V}} ]$$
belongs to ${\cal C}_2$.
The corresponding symmetric Hilbert-Schmidt bilinear form is easily calculated as
\[
(k, k')
\quad\mapsto \quad
\frac12 D^2_{k, k'}  \langle q(h), g_2 (\pi^h w;h)\rangle_{{\cal V}} 
=
\frac12 \bigl\la q(h), D^2\psi_1(h) \la \pi^h k, \pi^h k' \ra \bigr\ra_{{\cal V}}.
\]
From (\ref{cmax1/2.ineq})
and boundedness of ${\mathbb E} [\langle q(h), g_2 (\pi^h w;h)\rangle_{{\cal V}} ]$ in $h$,
we obtain (\ref{expint1.ineq}).
\QED


We give two technical lemmas (Lemmas \ref{lm.tw64} and  \ref{lm.tw65}) on integrability of certain Wiener functionals 
which will appear in the asymptotic expansion.
Since the solution of RDE is involved in these lemmas, 
we assume for safety that the coefficients are of $C_b^{\infty}$.
In the proofs of these lemmas 
we will write $g_j =g_j(w;h)$, $R^{\ve}_{ j}= R^{\ve}_{ j} (w;h)$ 
for simplicity of notations.

The first one below corresponds to Lemma 6.4, \cite{tw}. 
Note that the constant $\gamma$ depends on $r$, but not on $h$ with $\|h\|_{{\cal H}} \le r$.

\begin{lm}\label{lm.tw64}
Assume {\bf (A1)}, {\bf (B1)}--{\bf (B3)} and boundedness of $V_i~(0 \le i \le r)$.
Then, the following assertions {\rm (i)} and {\rm (ii)} hold:
\\
{\rm (i)}~For any $r>0$ and $c>0$, there exists $\gamma >0$ such that
\[
\sup_{0 <\ve \le 1}\sup_{ \|h\|_{{\cal H}} \le r}
{\mathbb E}
\Bigl[
\exp \Bigl( c  \Bigl|  \frac{R^{\ve}_{ 2} (\,\cdot\, , h)  }{\ve^2} \Bigr| \Bigr) 
 {\bf 1}_{U_{\gamma}} (\ve {\bf W})
\Bigr] <\infty.
\]
{\rm (ii)}~For any $r>0$ and $c>0$, there exists $\gamma >0$ such that
\[
\sup_{0 <\ve \le 1}\sup_{ \|h\|_{{\cal H}} \le r}
{\mathbb E}
\Bigl[
\exp \Bigl( c   \Bigl| \frac{R^{\ve}_{ 3}(\,\cdot\, , h) }{\ve^2}  \Bigr| \Bigr)  
{\bf 1}_{U_{\gamma} }(\ve {\bf W})
\Bigr] <\infty.
\]
\end{lm}

\Proof
We use 
the deterministic estimates (\ref{tylr1.eq})--(\ref{tylr3.eq}) 
for the Taylor-like expansion of the Lyons-It\^o map
and Besov-H\"older embedding for rough path spaces with $\alpha' = \alpha -1/(4m)$.
If $\ve {\bf W} \in U_{\gamma'}$, then 
\[
|R^{\ve}_{ j}| \le C (\ve + \| (\ve {\bf W})^1 \|_{\alpha, 4m}+ \| (\ve {\bf W})^2 \|_{2\alpha, 2m}^{1/2})^j 
\]
for some positive constant $C=C(j, \gamma', r)$.
We will take $\gamma' =1$ below.
For $0<\gamma <1$,
\[
c |R^{\ve}_{ 2}/\ve|^2 
\le 
c C (\ve +\gamma)^2 (1 + \| {\bf W}^1 \|_{\alpha, 4m}+ \| {\bf W}^2 \|_{2\alpha, 2m}^{1/2})^2.
\]
Hence, we have 
\[
\sup_{0 <\ve \le \gamma}\sup_{ \|h\|_{{\cal H}} \le r}
{\mathbb E}
\Bigl[
\exp ( c   |R^{\ve}_{ 2}/\ve|^2 )  {\bf 1}_{U_{\gamma} (\ve {\bf W})}
\Bigr] 
\le 
{\mathbb E}
\Bigl[
\exp \bigl( 
4c C\gamma^2 (1 + \| {\bf W}^1 \|_{\alpha, 4m}+ \| {\bf W}^2 \|_{2\alpha, 2m}^{1/2})^2
\bigr) 
\Bigr].
\]
By a Fernique-type theorem for Brownian rough path,
the right hand side  is integrable
if $\gamma$ is chosen sufficiently small. 
Once $\gamma$ is  fixed,
we can easily estimate the integral for $\ve \in [\gamma, 1]$.
Thus, we have shown assertion {\rm (i)}. 
The proof of assertion {\rm (ii)} is essentially the same.
\QED


The other one is a key technical lemma.
It corresponds to Lemma 6.5, \cite{tw}.
\begin{lm}\label{lm.tw65}
Assume {\bf (A1)}, {\bf (B1)}--{\bf (B3)} and boundedness of $V_i~(0 \le i \le r)$.
Let $c_1 >1$ be as in (\ref{tw_eq_6.4}) in Lemma \ref{lm.tw61}. 
Then, for any $c_2 \in (1,c_1)$, there exists a constant $\gamma_1 >0$ 
which is independent of $h \in  {\cal K}^{min}_a$ and 
satisfies that
\[
\sup_{0 <\ve \le 1}\sup_{ h \in  {\cal K}^{min}_a}
{\mathbb E}
\Bigl[
\exp ( c_2  \la  q (h), R^{\ve}_{ 2} (\,\cdot\, , h) \ra_{{\cal V}} /\ve^2)
  {\bf 1}_{U_{\gamma_1} }(\ve {\bf W})
~;~
| {\bf i}_h |^2 + |R^{\ve}_{1} (\,\cdot\, , h)/\ve|^2 \le \kappa^2
\Bigr] 
<\infty
\]
for any $\kappa >0$.
\end{lm}

\Proof
Thanks to Lemma \ref{lm.tw64}, {\rm (ii)} and H\"older's  inequality,
it is sufficient to show the following integrability.
For any $c_2 \in (1,c_1)$, there exists $\gamma >0$ such that
\begin{equation}\label{lm.tw64.suff}
\sup_{0 <\ve \le 1}\sup_{ h \in  {\cal K}^{min}_a}
{\mathbb E}
\Bigl[
\exp ( c_2  \la  q (h), g_{ 2}\ra_{{\cal V}} )  {\bf 1}_{U_{\gamma} }(\ve {\bf W})
~;~
| {\bf i}_h |^2 + |R^{\ve}_{1}/\ve|^2 \le \kappa^2
\Bigr] 
<\infty
\end{equation}
for any $\kappa >0$.

By straightforward computation, we see that
\begin{eqnarray*}
g_2 (w;h) 
&=&
\lim_{l \to \infty} g_2 (w(l);h) 
\nn\\
&=&
\lim_{l \to \infty}
\Bigl\{
\frac12 D^2 \psi_1 (h)\la  w(l),w(l)   \ra 
+ 
\Pi_{{\cal V}}  J_1 (h) \int_0^1  J_s (h)^{-1} V_0 (\phi_s (h)) ds
\Bigr\}
\\
&=&
\lim_{l \to \infty}
\Bigl\{
\frac12 D^2 \psi_1 (h)\la  (\pi^h w)(l), (\pi^h w)(l)   \ra 
+ 
\Pi_{{\cal V}}  J_1 (h) \int_0^1  J_s (h)^{-1} V_0 (\phi_s (h)) ds
\Bigr\}
\\
&&
+
\lim_{l \to \infty}
\Bigl\{
 D^2 \psi_1 (h)\la  (\pi^h w)(l), ((\pi^h)^{\bot} w)(l)   \ra 
+ 
\frac12 D^2 \psi_1 (h)\la ((\pi^h)^{\bot} w)(l)  , ((\pi^h)^{\bot} w)(l)   \ra 
\Bigr\}
\\
&=&
g_2 (\pi^h w;h) 
+
D^2 \psi_1 (h)\la  \pi^h w, (\pi^h)^{\bot} w  \ra 
+\frac12
D^2 \psi_1 (h)\la   (\pi^h)^{\bot} w  , (\pi^h)^{\bot} w  \ra. 
\end{eqnarray*}
Here, we used basic properties of Young translation and
the fact that
$g_2 (\,\cdot\,; h)$ is actually a continuous function in ${\bf W}={\cal L}(w)$.

There exists a positive constant $C$ (which depends only on $\|h\|_{{\cal H}}$
and may vary from line to line)
such that 
\begin{eqnarray}
|g_2 (w;h) - g_2 (\pi^h w;h) | 
&\le&
C (  \| \pi^h w \|_{{\cal W}}  \| (\pi^h)^{\bot} w\|_{{\cal H}} +  \| (\pi^h)^{\bot} w\|_{{\cal H}}^2 )
\nn\\
&\le&
C (  \|  w \|_{{\cal W}}  \| (\pi^h)^{\bot} w\|_{{\cal H}} 
+  \| (\pi^h)^{\bot} w\|_{{\cal H}}^2 )
\nn\\
&\le&
C \bigl\{ \frac{\rho^2}{2} \|  w \|_{{\cal W}}^2
+ 
( \frac{1}{2 \rho^2} +2 ) \| (\pi^h)^{\bot} w\|_{{\cal H}}^2 \bigr\},
\label{lm64.ineq2}
\end{eqnarray}
where $\rho >0$ is a small constant which will be determined later
and we used the estimate for $D^2 \phi_1$ 
given just below (\ref{D2phi.eq}).

Under the condition that $| {\bf i}_h |^2 + |R^{\ve}_{1}/\ve|^2 \le \kappa^2$,
we can easily see that 
$|g_1 (w;h)| \le | R^{\ve}_{ 1}/\ve - R^{\ve}_{ 2}/\ve| \le \kappa + |R^{\ve}_{ 2}/\ve|$.
From these, we have
$\| (\pi^h)^{\bot} w\|_{{\cal H}} \le C (\kappa + |R^{\ve}_{ 2}/\ve|)$
and 
\[
|g_2 (w;h) - g_2 (\pi^h w;h) | 
\le 
\frac{C\rho^2}{2} \|  w \|_{{\cal W}}^2
+ 
C ( \frac{1}{2 \rho^2} +2 ) (\kappa^2 + \Bigl|\frac{R^{\ve}_{ 2} }{\ve}\Bigr|^2 ).
\]
Note that $C >0$ is independent of $\kappa, \ve, \rho, w, h$.


Using H\"older's inequality with $p=c_1/c_2 \in (1,\infty)$ 
and $1/p +1/p' =1$, we have
\begin{eqnarray}
\lefteqn{
{\mathbb E}
\Bigl[
\exp ( c_2  \la  q (h), g_{ 2}\ra_{{\cal V}} )  {\bf 1}_{U_{\gamma} }(\ve {\bf W})
~;~
| {\bf i}_h |^2 + |R^{\ve}_{1}/\ve|^2 \le \kappa^2
\Bigr] 
}
\nn\\
&\le&
{\mathbb E}
\Bigl[
\exp ( c_1  \la  q (h), g_{ 2}\circ \pi^h \ra_{{\cal V}} )  {\bf 1}_{U_{\gamma}} (\ve {\bf W})
~;~
| {\bf i}_h |^2 + |R^{\ve}_{1}/\ve|^2 \le \kappa^2
\Bigr]^{\frac{1}{p}} 
\nn\\
&&
\times 
{\mathbb E}
\Bigl[
\exp ( p' c_2  \max_{h \in {\cal K}^{min}_a}|q (h)|  |g_{ 2} - g_{ 2}\circ \pi^h| )  
{\bf 1}_{U_{\gamma} }(\ve {\bf W})
~;~
| {\bf i}_h |^2 + |R^{\ve}_{1}/\ve|^2 \le \kappa^2
\Bigr]^{\frac{1}{p'}} 
\nn\\
&\le&
(2\pi)^{(n+n' )/2p} \{ \det \sigma [\psi_1] (h) \}^{1/2p}
{\mathbb E} \Bigl[
\exp (  c_1 \langle q(h), g_2\rangle_{{\cal V}})
\delta_0 ( g_1,  {\bf i}_h )
\Bigr]^{\frac{1}{p}} 
\nn\\
&& 
\times 
{\mathbb E}
\Bigl[
\exp \Bigl(  
\frac{C^{\prime} \rho^2}{2} \|  w \|_{{\cal W}}^2
+ 
C^{\prime} ( \frac{1}{2 \rho^2} +2 ) (\kappa^2 + \Bigl|\frac{R^{\ve}_{ 2} }{\ve}\Bigr|^2 )
\Bigr)  
{\bf 1}_{U_{\gamma}} (\ve {\bf W})
\Bigr]^{\frac{1}{p'}},
\label{tw64.ineq2}
\end{eqnarray}
where we set $C^{\prime} = Cp' c_2  \max_{h \in {\cal K}^{min}_a}|q (h)|$.
The first factor on the right hand side of (\ref{tw64.ineq2}) above is dominated by
a positive constant independent of $h \in {\cal K}^{min}_a$,
due to (\ref{tw_eq_6.4}) or (\ref{expint1.ineq}).
Now, we choose $\rho$ so small that $\exp ( C^{\prime} \rho^2 \|  w \|_{{\cal W}}^2)$
is integrable, which is possible by Fernique's theorem, 
and use Schwarz' inequality for the second factor.
Then, by Lemma \ref{lm.tw64}, {\rm (i)},
we can choose $\gamma >0$ so that $\sup_{\ve} \sup_{h}$ of the second factor is finite.
Thus, we have shown (\ref{lm.tw64.suff}), which completes the proof.
\QED


\begin{co}\label{co.newadd}
Keep  the same notations and assumptions as in Lemma \ref{lm.tw65}.
Then, the following {\rm (i)}--{\rm (iii)} hold:
\\
\noindent
{\rm (i)}~
For any $c_2 \in (1,c_1)$ and any $\kappa >0$,  
$$
\sup_{ h \in  {\cal K}^{min}_a}
{\mathbb E}
\Bigl[
\exp ( c_2  \la  q (h), g_{ 2} (\,\cdot\, , h) \ra_{{\cal V}} )
~;~
| {\bf i}_h |^2 + |g_{1} (\,\cdot\, , h)|^2 \le \kappa^2
\Bigr] 
<\infty.
$$
\noindent
{\rm (ii)}~For any smooth function $f$ on ${\cal V} \times T_h {\cal K}^{min}_a$
with compact support,
$$
\exp ( c_2  \la  q (h), g_{ 2} (\,\cdot\, , h) \ra_{{\cal V}} )
f( g_{1} (\,\cdot\, , h), {\bf i}_h ) \in \tilde{\mathbb D}_{\infty}.
$$
Moreover, 
as $h$ varies in ${\cal K}^{min}_a$, 
these Wiener functionals  form  a bounded set in $\tilde{\mathbb D}_{\infty}$.
\\
\noindent
{\rm (iii)}~Let $T \in {\cal S}^{\prime} ( {\cal V} \times T_h {\cal K}^{min}_a)$
with compact support.
Take any two smooth function $f$ and $\hat{f}$ 
on ${\cal V} \times T_h {\cal K}^{min}_a$ with compact support
such that $f \equiv 1 \equiv  \hat{f}$ on the support of $T$.
Then, 
\begin{eqnarray*}
\lefteqn{
\exp ( c_2  \la  q (h), g_{ 2} (\,\cdot\, , h) \ra_{{\cal V}} )
f ( g_{1} (\,\cdot\, , h), {\bf i}_h )  T ( g_{1} (\,\cdot\, , h), {\bf i}_h )
}\\
&=&
\exp ( c_2  \la  q (h), g_{ 2} (\,\cdot\, , h) \ra_{{\cal V}} )
\hat{f} ( g_{1} (\,\cdot\, , h), {\bf i}_h )  T ( g_{1} (\,\cdot\, , h), {\bf i}_h )
\in 
{\mathbb D}_{-\infty}.
\end{eqnarray*}
In particular, the generalized expectations of 
these two Watanabe distributions also coincide.
Moreover, as $h$ varies in ${\cal K}^{min}_a$,
these  Watanabe distributions form  a bounded set in ${\mathbb D}_{-\infty}$.
\end{co}

\Proof
We can easily prove  {\rm (i)}  by applying Fatou's lemma 
to the inequality in Lemma \ref{lm.tw65}. 
Noting that 
\[
{\bf 1}_{ \{ | {\bf i}_h |^2 + |g_{1} (\,\cdot\, , h)|^2 \le (\kappa/2)^2 \} }
\le 
\liminf_{\ve\searrow 0} 
{\bf 1}_{ 
\{ | {\bf i}_h |^2 + | R_1^{\ve} (\,\cdot\, , h) /\ve|^2 \le \kappa^2 \} 
}
\]
for almost all $w$,  we have the desired inequality with $\kappa$
being replaced by $\kappa/2$.
Since, $\kappa >0$ is arbitrary, we have shown {\rm (i)}.

It is straightforward to check  {\rm (ii)} from {\rm (i)}.
We now prove {\rm (iii)}.
Since $T ( g_{1} (\,\cdot\, , h), {\bf i}_h )$ is a well-defined element in 
$\tilde{\mathbb D}_{-\infty}$,
we see from {\rm (ii)} that both sides belong to ${\mathbb D}_{-\infty}$.
To check the equality, 
we just have to use 
$f \cdot T = T = \hat{f} \cdot T$ as finite-dimensional distributions.
\QED

\begin{re}\normalfont\label{re.add}
We may and will write 
$\exp ( c_2  \la  q (h), g_{ 2} (\,\cdot\, , h) \ra_{{\cal V}} )
T( g_{1} (\,\cdot\, , h), {\bf i}_h ) \in {\mathbb D}_{-\infty}$ for simplicity
in the situation of Corollary \ref{co.newadd}, {\rm (iii)}.
Later, we will usually choose 
$T$  to be a partial derivative of the delta function.
Note that since 
$\exp ( c_2  \la  q (h), g_{ 2} (\,\cdot\, , h) \ra_{{\cal V}} ) 
\notin \tilde{\mathbb D}_{\infty}$,
the product $\exp ( c_2  \la  q (h), g_{ 2} (\,\cdot\, , h) \ra_{{\cal V}} ) \Phi$
or its generalized expectation 
cannot be defined for a general element
$\Phi \in  \tilde{\mathbb D}_{-\infty}$.
\end{re}


\section{ Asymptotic partition of unity   }
\label{sec.apu}

In this section, by using rough path theory,
we modify the argument on an asymptotic partition of unity
in Section 6.6,  \cite{tw}.



Let $f \colon {\mathbb R} \to {\mathbb R}$
be an even smooth
 function such that $f(s) =1$ if $|s| \le 1/2$ and $f(s)=0$ if $|s| \ge 1$.
We also assume that $f$ is non-increasing on $[0, \infty)$.
Define, for $\gamma >0$, $\ve \in (0,1]$ and $h \in {\cal H}$,
\begin{eqnarray}
\zeta^{\ve}_{\gamma} (w) 
&=& 
f \Bigl(
\frac{ \| (\ve {\bf W})^1 \|_{\alpha, 4m -B}^{4m}+ \| (\ve {\bf W})^2 \|_{2\alpha, 2m -B}^{2m}}{\gamma^{4m}}
\Bigr),
\nn
\\
\zeta^{\ve,h}_{\gamma} (w) 
&=&
\zeta^{\ve}_{\gamma} (w -\frac{h}{\ve})
= 
f \Bigl(
\frac{ \| \tau_{-h}(\ve {\bf W})^1 \|_{\alpha, 4m -B}^{4m}+ \|  \tau_{-h}(\ve {\bf W})^2 \|_{2\alpha, 2m -B}^{2m}}{\gamma^{4m}}
\Bigr).
\nn
\end{eqnarray}
Here, $\tau_{-h}$ stands for the Young translation by $-h$ on the geometric rough path space.
Note that these are ${\mathbb D}_{\infty}$-functionals.
It is clear that {\rm (i)}
$\zeta^{\ve}_{\gamma} =0$ if  
$\ve {\bf W} \notin U_{\gamma}$
and $\zeta^{\ve,h}_{\gamma}=0$
if   
$\ve {\bf W} \notin U_{h, \gamma}$
and 
{\rm (ii)} $\zeta^{\ve}_{\gamma}=1$ 
if $\ve {\bf W} \in U^{\prime}_{\gamma}$ 
and $\zeta^{\ve,h}_{\gamma}=1$
if 
$\ve {\bf W} \in U^{\prime}_{h, \gamma}$.


%
For
$h_1, \ldots, h_N \in {\cal K}^{min}_a$
and $\gamma_1, \ldots, \gamma_N >0$ 
satisfying that 
${\cal L}({\cal K}^{min}_a) 
\subset \cup_{\nu=1}^N U_{h_{\nu}, \gamma_{\nu}}^{\prime}$,
we define
\begin{eqnarray}
\chi^{\ve} (w) 
&=& 
1 - \prod_{\nu=1}^N (1 - \zeta^{\ve,h_{\nu}}_{\gamma_{\nu}} (w) ),
\label{apu1.def}
\\
\tilde{\chi}^{\ve}_{\nu} (w) 
&=& 
 \zeta^{\ve,h_{\nu}}_{\gamma_{\nu}} (w)\prod_{\mu=1}^{\nu -1} (1 - \zeta^{\ve,h_{\mu}}_{\gamma_{\mu}} (w) )
\label{apu2.def}
\end{eqnarray}
for $1 \le \nu \le N$.
Loosely speaking, $\chi^{\ve} (w)=1$ means that $\ve{\bf w}$ is close to ${\bf h}_{\nu}$ for some $\nu$.
Similarly, $\tilde{\chi}^{\ve}_{\nu} (w)=1$ loosely means that 
$\ve{\bf w}$ is close to ${\bf h}_{\nu}$, but distant from ${\bf h}_{1}, \ldots,{\bf h}_{\nu -1}$.
It is clear that 
$0 \le \chi^{\ve}\le 1$ and easy to see that $\sum_{\nu=1}^N \tilde{\chi}^{\ve}_{\nu}=  \chi^{\ve}$.


It is easy to see that 
$$
\zeta^{\ve,h_{\nu}}_{\gamma_{\nu}} (w + \frac{k}{\ve}) 
=
1 + O(\ve^l)
\quad 
\mbox{in ${\mathbb D}_{\infty}$ as $\ve \searrow 0$ for any $l >0$}
$$
uniformly in $k \in {\cal L}^{-1} (U_{h_{\nu}, \gamma_{\nu}}^{\prime}) \cap {\cal K}^{min}_a$.
(Note that the Wiener functional on the left hand is the composition of $f$
and a polynomial in $\ve$ with the coefficients from an inhomogeneous Wiener chaos.
Therefore, it clearly has an asymptotic expansion in ${\mathbb D}_{\infty}$-topology
and we just need to check that the coefficients vanish by formal differentiation.)
In a similar way, we have
\begin{equation}\label{asy_chi_ve}
\chi^{\ve} (w + \frac{k}{\ve}) 
=
1 + O(\ve^l)
\quad 
\mbox{in ${\mathbb D}_{\infty}$ as $\ve \searrow 0$ for any $l >0$}
\end{equation}
uniformly in $k \in  {\cal K}^{min}_a$.
Note that we have repeatedly used the fact that $f(x)=1$ on $[-1/2,1/2]$.


Now we get back to our SDE/RDE and choose such 
$N$, $h_{\nu}$ and $\gamma_{\nu}$ ($1 \le \nu \le N$).
Before going into details, we will explain below why our problem 
is reduced to the case where $V_i~(0 \le i \le r)$ is of $C_b^{\infty}$.  
(This remark is unnecessary if $V_i$ is of $C_b^{\infty}$
in the first place.)
\begin{re}\normalfont\label{re.bdd_ok}
Since the deterministic It\^o map is continuous
under {\bf (A1)},
\[
R := \sup \{ | x + \Phi ({\bf h}, \lambda^{\ve})^1_{0,t} | 
\mid t \in [0,1], \, \ve \in [0,1], \, h \in {\cal K}^{min}_a \}  \in (0, \infty).
\]
From a result in \cite{be} (or a standard cut-off argument, alternatively)
we see the following:
There exists a neighborhood $\Lambda$ of ${\cal L} ({\cal K}^{min}_a )$
in $G\Omega^B_{\alpha, 4m -B} ({\mathbb R}^r)$ such that 
\[
\sup \{ | x + \Phi ({\bf w}, \lambda^{\ve})^1_{0,t} | 
\mid t \in [0,1], \, \ve \in [0,1], \, {\bf w} \in \Lambda \}  
\le 
2R.
\]
Therefore, if we work on the neighborhood,
we only need information of $V_i$ restricted to the ball of radius $2R$.

In the argument in the rest of this section, 
taking $\gamma =\gamma_h >0$ smaller if necessary,
we have $U_{h,\gamma} \subset \Lambda$.
For these reasons, we may and will assume in the sequel
that $V_i$ is of $C_b^{\infty}$.

(When we throw away contributions from the complement set of 
a neighborhood of ${\cal L} ({\cal K}^{min}_a)$ in Lemma \ref{lm.ldp_up} below,
we use the upper bound estimate of the LDP in Theorem \ref{tm.ub_ldp}.
Hence, it is important that Theorem \ref{tm.ub_ldp} holds under {\bf (A1)}.)
\end{re}


For any $h \in {\cal K}^{min}_a$, 
there exists a coordinate neighborhood $O_h$ of $h$ (with respect to ${\cal H}$-topology)
and a constant $\gamma= \gamma_h >0$ such that the following conditions holds:
\begin{eqnarray}
{\cal K}^{min}_a \cap U_{ h,\gamma}
=
\{ k \in  {\cal K}^{min}_a \mid  k \in U_{ h,\gamma} \} \subset O_h.
\label{223_1.cond}
\\
\mbox{For any $\ve \in (0,1]$,}
\quad
\int_{O_h} D^{(k)} (\ve w)\Delta^{(k)} (\ve w) \omega(dk) =1
\quad
\mbox{on $\{ w \in W \mid \ve {\bf W} \in U_{ h,\gamma}  \}$.}
\label{223_2.cond}
\\
\sup_{k \in O_h}
 \Bigl( \sum_{j=1}^{n'} \bigl| 
 \bigl( 
 \frac{\partial \iota}{\partial \theta^j} 
 \bigr)_k
 \bigr|_{{\cal H}}  
 +
 \sum_{j, j'=1}^{n'} \bigl| 
 \bigl(
  \frac{\partial^2 \iota}{\partial \theta^j \partial \theta^{j'}} 
 \bigl)_k
  \bigr|_{{\cal H}}
   \Bigr) <\infty
\qquad
\mbox{and}
\qquad
\inf_{k \in O_h} \det G(k) >0.
\label{223_5.cond}
\end{eqnarray}
Here, $(\theta^1, \ldots, \theta^{n'})$ is the local coordinate on $O_h$.
Note that (\ref{223_1.cond}) and  (\ref{223_2.cond})
are immediate from Lemma \ref{lm.tw_62} and Remark \ref{re.tw_62}.
(Note that $\gamma$ above can be taken smaller if necessary. See Remark  \ref{re.tw_62}.)


We choose any  $1< c_2 < c_1$ as in Lemmas \ref{lm.tw61} and \ref{lm.tw65} and fix them in what follows.
Let $\lambda >0$, $\gamma_0 >0$ and $\ve_0 \in (0,1]$
be the constants defined in (\ref{lamb.def})
and  Lemma \ref{lm.tw63}, respectively.
Similarly, let $\gamma_1 >0$ be the
constant that appeared in Lemma \ref{lm.tw65}.

We write $\rho= \rho_h:= \sup_{k \in O_h} \| k-h \|_{{\cal H}} >0$.
Then, we can take 
$O_h$
and $\gamma=\gamma_h >0$ so small that 
the following inequalities hold:
\begin{align}
\sigma=
\sigma_{h} := \frac{16^{1/4m} C \rho}{\gamma}  &\ge 2^{1/4m},
\label{csigma.eq}
\\
2 ( 1 + 32^{1/4m})^{4m}  C^{4m} \rho^{4m} 
=
2( 2^{1/4m} \sigma \gamma +C\rho )^{4m} 
  &\le \gamma_0^{4m},
\label{crho.eq}
\\
2 (\gamma +C\rho )^{4m} &\le \gamma_1^{4m}.
\label{csr.eq}
\end{align}
Here, 
$C=C_{\alpha, 4m} >0$ is the smallest constant that satifies the 
estimates (\ref{shift1.ineq})--(\ref{shift2.ineq}).
It is easy to check that such $\rho$ and $\gamma$ exist.
First, take $O_h$ (or $\rho$) small enough and fix it 
and then we can retake $\gamma >0$ smaller if necessary for this fixed $\rho$.


Note that the condition (\ref{crho.eq}) immediately implies the following:
If $\ve {\bf W} \in U_{ h-k, 2^{1/4m}\sigma\gamma}$,
that is,
\begin{equation}\label{k_h_gam.eq}
\| \tau_{k-h}(\ve {\bf W})^1 \|_{\alpha, 4m -B}^{4m} + \| \tau_{k-h}(\ve {\bf W})^2 \|_{2\alpha, 2m -B}^{2m}
\le 
2   \sigma^{4m} \gamma^{4m} 
\end{equation}
holds for $k \in O_h$,
then $\| (\ve {\bf W})^1 \|_{\alpha, 4m -B}^{4m} 
+ \| (\ve {\bf W})^2 \|_{2\alpha, 2m -B}^{2m} \le \gamma_0^{4m}$, 
because of 
(\ref{shift1.ineq}), (\ref{shift2.ineq}), (\ref{crho.eq}) 
and the obvious fact that $\ve {\bf W} = \tau_{-k+h}\circ \tau_{k-h}(\ve {\bf W})$.
(Thus, the assumption of Lemma \ref{lm.tw63} is satisfied.)
Therefore, under (\ref{k_h_gam.eq})  we have 
\begin{equation}
\inf_{k \in O_h}
 \inf_{z : |z|=1} 
z^*    \sigma [(R^{\ve}_{ 1}(\,\cdot\, ;k) /\ve, {\bf i}_k ) ] (w) z 
\ge \frac{\lambda \wedge 1}{2}
\qquad
\mbox{if $\ve \in (0,\ve_0]$ and $\ve {\bf W} \in U_{ h-k, 2^{1/4m}\sigma\gamma}$}.
\label{223_3.cond}
\end{equation}

Likewise, if 
$\| \tau_{k-h}(\ve {\bf W})^1 \|_{\alpha, 4m -B}^{4m} + \| \tau_{k-h}(\ve {\bf W})^2 \|_{2\alpha, 2m -B}^{2m}
\le \gamma^{4m}$
for $k \in O_h$,
then (\ref{shift1.ineq}), (\ref{shift2.ineq}) and  (\ref{csr.eq}) 
imply that  
$\| (\ve {\bf W})^1 \|_{\alpha, 4m -B}^{4m} 
+ \| (\ve {\bf W})^2 \|_{2\alpha, 2m -B}^{2m} \le \gamma_1^{4m}$. 
In other words, $U_{h-k, \gamma} \subset U_{\gamma_1}$ for $k \in O_h$.
Therefore, by Lemma \ref{lm.tw65},
\begin{equation}
\sup_{0 <\ve \le 1}\sup_{ k \in  O_h}
{\mathbb E}
\Bigl[
\exp \Bigl( \frac{c_2}{\ve^2}  \la  q (k), R^{\ve}_{ 2}(w ;k)\ra_{{\cal V}} \Bigr)  
{\bf 1}_{U_{ h-k, \gamma} } (\ve {\bf W})
~;~
| {\bf i}_k |^2 + \Bigl| \frac{ R^{\ve}_{1}(w ;k)}{\ve} \Bigr|^2 \le \kappa^2
\Bigr] 
<\infty
\label{223_4.cond}
\end{equation}
for any $\kappa >0$.

Thus, we have shown the following lemma:
\begin{lm}\label{lm.prm}
Assume {\bf (A1)} and {\bf (B1)}--{\bf (B3)}.
For any $h \in {\cal K}^{min}_a$, we can find $O_h$ and $\gamma_h$ 
such that (\ref{223_1.cond}),   (\ref{223_2.cond}), 
 (\ref{223_5.cond}), 
 (\ref{223_3.cond}), 
  (\ref{223_4.cond}) hold. 
\end{lm}

Since ${\cal K}^{min}_a$ is compact, there exist finitely many  
$h_1, \ldots, h_N \in {\cal K}^{min}_a$ and 
$\gamma_1, \ldots, \gamma_N >0$
such that 
${\cal L}({\cal K}^{min}_a) 
\subset \cup_{\nu=1}^N U_{h_{\nu}, \gamma_{\nu}}^{\prime}$,
where we wrote $\gamma_{\nu} =\gamma_{h_{\nu}}$ 
for simplicity ($1 \le \nu \le N$).
In the next section 
we denote by $\chi^{\ve}$ and $\tilde{\chi}^{\ve}_{\nu}$
the asymptotic partition of unity defined in (\ref{apu1.def}), (\ref{apu2.def})
associated with
these data $\{(h_{\nu}, \gamma_{\nu})\}_{\nu=1}^N$.

\begin{re}\normalfont
The definitions of $\chi^{\ve}$ and $\tilde{\chi}^{\ve}_{\nu}$ 
in (\ref{apu1.def}) and (\ref{apu2.def})
depend only on the data $N$
and $\{(h_{\nu}, \gamma_{\nu})\}_{\nu=1}^N$.
Information of the SDE is involved only in choosing such data.
Therefore, 
our version of the asymptotic 
partition of unity is much simpler than the one in \cite{tw}
and looks more powerful.
This is one of the advantages of our rough path approach.
(See Section 6.6 \cite{tw} and compare.)
\end{re}


\section{ Proof of the main thorem: The Euclidean case }
\label{sec_pf}

In this section we prove our main theorem in the Euclidean setting
(Theorem \ref{tm.main}).
We assume {\bf (A1)}--{\bf (A2)}, {\bf (B1)}--{\bf (B3)} and  {\bf (C1)}
throughout this section.
We set $G^{\ve} (w)=G(\ve, w)$ below.
Define
\begin{eqnarray}
{\mathbb E}[G^{\ve} \delta_a (Y^{\ve}_1) ]
=
{\mathbb E}[(1 - \chi^{\ve}) G^{\ve}  \delta_a (Y^{\ve}_1) ]
+
{\mathbb E}[\chi^{\ve} G^{\ve}  \delta_a (Y^{\ve}_1) ]
=:
I_1 (\ve)+ I_2 (\ve).
\end{eqnarray}
Due to the large deviation upper bound (Theorem \ref{tm.ub_ldp}), 
the first term $I_1 (\ve)$ does not 
contribute to the asymptotic expansion.
\begin{lm}\label{lm.ldp_up}
There exist positive constants $c, C$ such that 
\[
|I_1 (\ve)| \le  C \exp \bigl( - \frac{d_a^2 +c}{2\ve^2} \bigr) 
\qquad
\mbox{as $\ve \searrow 0$.}
\]
Here, $d_a = \inf \{ \|h\|_{{\cal H}} \mid \psi(1,x, h) =a\} >0$.
\end{lm}

\Proof
We use Theorem \ref{tm.ub_ldp}.
By the way of construction,  
$1 - \chi^{\ve}$ is dominated by ${\bf 1}_{ \{ \ve {\bf W} \in O^c\} }$
for some open set 
$O \subset G\Omega^B_{\alpha, 4m} ({\mathbb R}^d)$
containing ${\cal K}_a^{min}$.
Therefore, we have
\[
|I_1 (\ve)| =
\Bigl| \int_{{\cal W}}  (1 - \chi^{\ve}) \tilde{G}^{\ve}  \theta_{a}^{\ve}(dw) \Bigr|
\le 
 \mu_{a}^{\ve}( O^c)^{1/p}
  \Bigl\{
 \int_{{\cal W}}  | \tilde{G}^{\ve} |^{p'} \theta_{a}^{\ve}(dw) 
 \Bigr\}^{1/p'}
\]
for any $p, p' \in (1,\infty)$ with $1/p +1/p' =1$.
Here, 
$ \tilde{G}^{\ve}$ is an $\infty$-quasi continuous modification of ${G}^{\ve}$.

Since $O^c$ is closed and the rate function $I$ is good, 
it is easy to see that $\inf_{{\bf w} \in O^c} I  ({\bf w}) > d_a^2 /2$.
Set $\kappa := (\inf_{{\bf w} \in O^c} I  ({\bf w}) -  d_a^2 /2)/4 >0$.
By Theorem \ref{tm.ub_ldp} we have
\[
 \mu_{a}^{\ve}( O^c)^{1/p} 
=
O \Bigl(  \exp \bigl( \frac{ - \inf_{{\bf w} \in O^c}  I  ({\bf w})+\kappa}{ p \ve^2} \bigr) \Bigr)
\le O \Bigl( 
\exp \bigl( \frac{ - \inf_{{\bf w} \in O^c}  I  ({\bf w})+2\kappa}{  \ve^2} \bigr)\Bigr)
 \]
as $\ve \searrow 0$ if $p >1$ is sufficiently close to $1$.

On the other hand, when $p' =2l$ with arbitrarily large  $l \in {\mathbb N}$, 
\[
 \int_{{\cal W}}  | \tilde{G}^{\ve} |^{2l} \theta_{a}^{\ve}(dw) 
 =
 {\mathbb E} [ (G^{\ve} )^{2l}  \delta_a (Y^{\ve}_1) ]
 \le 
\| (G^{\ve} )^{2l}  \|_{j, 2}   \| \delta_a (Y^{\ve}_1) \|_{-j, 2}
\]
for some $j \in {\mathbb N}$ such that $ \delta_a (Y^{\ve}_1)$ belongs to the Sobolev space ${\mathbb D}_{-j,2}$.
($j$ is independent of $l$.)
Moreover, it is well-known that 
$ \| \delta_a (Y^{\ve}_1) \|_{-j, 2} = O (\ve^{ - \nu})$ as $\ve \searrow 0$ for some $\nu >0$.
(This comes from {\bf (A2)} and the integration by parts formula.) 
Since we assumed in {\bf (C1)} that
 $G^{\ve}$ is bounded in ${\mathbb D}_{\infty}$,
$(G^{\ve} )^{2l} $ is bounded in ${\mathbb D}_{j,2}$  for any $l, j$.
Combining these all, we obtain 
\[
|I_1 (\ve)| = O \Bigl( 
\exp \bigl( \frac{ - \inf_{{\bf w} \in O^c}  I  ({\bf w})+3\kappa}{  \ve^2} \bigr)\Bigr)
 \]
as $\ve \searrow 0$.
If we set $c =2\kappa$ this is the statement of the lemma.
\QED


Hence, the problem reduces to the asymptotic expansion for 
$I_2 (\ve)
=
\sum_{\nu=1}^N
{\mathbb E}[\tilde\chi^{\ve}_{\nu} G^{\ve} \delta_a (Y^{\ve}_1) ].
$
For the reason we stated in  Remark \ref {re.bdd_ok},
we may additionally assume that $V_i~(0 \le i \le r)$ is bounded.
Now we compute each summand.

Take $\eta \in C_0^{\infty} ({\cal V}, {\mathbb R})$
with $\eta \ge 0$ and $\int_{{\cal V}} \eta (y) dy =1$
and
set 
$\eta_l (y) = l^n \eta (l y)$ for $l\ge 1$.
Then, 
$\eta_l (Y^{\ve}_1 -a) \to \delta_a  (Y^{\ve}_1)$ as $l \to \infty$
in $\tilde{\mathbb D}_{-\infty}$.
Noting that
$\tilde\chi^{\ve}_{\nu} (w) = \tilde\chi^{\ve}_{\nu} (w)
 {\bf 1}_{ \{\ve {\bf W} \in U_{h_{\nu}, \gamma_{\nu}}  \} }$
by definition,
we have from (\ref{223_2.cond}) that
\begin{eqnarray}
{\mathbb E}[\tilde\chi^{\ve}_{\nu} G^{\ve} \cdot \eta_l (Y^{\ve}_1 -a) ]
&=&
{\mathbb E}[\tilde\chi^{\ve}_{\nu} G^{\ve} \cdot \eta_l (Y^{\ve}_1 -a)  
\cdot
\underbrace{ \int_{O_{\nu} } D^{(k)} (\ve w)\Delta^{(k)} (\ve w) \omega(dk) }_{=1}]
\nn\\
&=&
\int_{O_{\nu}}
{\mathbb E}[\tilde\chi^{\ve}_{\nu} G^{\ve} \cdot \eta_l (Y^{\ve}_1 -a)  
\cdot
 D^{(k)} (\ve w)\Delta^{(k)} (\ve w)  ]
\omega(dk),
\label{sonyu.eq}
\end{eqnarray}
where we set $O_{\nu} = O_{h_{\nu}}$.

Let 
$(\theta^1, \ldots, \theta^{n'})$ be the local coordinate on $O_{\nu}$
and 
let $\hat{O}_{\nu} \subset {\mathbb R}^{n'}$ be the image of $O_{\nu}$.
As before
$k \in O_{\nu}$ is denoted by $k(\theta), ~ \theta \in \hat{O}_{\nu}$.
Set 
$G (\theta) = G(k(\theta))$ and
$[e_{\beta}  (\theta)]_{\beta =1}^{n'} = G (\theta)^{-1/2}
 [ (\partial k/\partial \theta^{\beta})  (\theta)]_{\beta =1}^{n'}$.
 Then, $\{e_{\beta}  (\theta) \}_{\beta =1}^{n'}$ 
 is an orthonormal basis of  $T_{k(\theta)} {\cal K}^{min}_a\subset {\cal H}$.
We set ${\bf e}(\theta,w) = [\la e_{\beta}  (\theta), w\ra]_{\beta =1}^{n'}$,
which is a continuous linear map from ${\cal W}$ to ${\mathbb R}^{n'}$.
(Recall that 
${\bf i}_{k(\theta)} \la w\ra
=\sum_{\beta =1}^{n'}  \la e_{\beta}  (\theta), w\ra e_{\beta}(\theta)$, 
which is almost the same as ${\bf e}(\theta,w)$, but it takes values in
$T_{k(\theta)} {\cal K}^{min}_a$.)
By definition we have 
$
D^{(k(\theta))} (\ve w)\Delta^{(k(\theta))} (\ve w)
=
\det[-\ve  {\bf a}_{k(\theta)} (w)]  \delta_0 (\ve {\bf e}(\theta,w)).
$

Combining these with the Cameron-Martin translation by $k(\theta) /\ve$, we see that 
\begin{eqnarray}
\lefteqn{
{\mathbb E}[\tilde\chi^{\ve}_{\nu} G^{\ve}  \eta_l (Y^{\ve}_1 -a) ]
}\nn\\
&=&
\int_{ \hat{O}_{\nu}} 
{\mathbb E}
\bigl[
\tilde\chi^{\ve}_{\nu} G^{\ve}  \eta_l (Y^{\ve}_1 -a) 
\det[-\ve  {\bf a}_{k(\theta)} ] \delta_0 (\ve {\bf e}(\theta,\,\cdot\,))
\bigr]
\sqrt{ \det G(\theta)}
d\theta
\nn\\
&=&
\int_{ \hat{O}_{\nu}} 
{\mathbb E}
\Bigl[
\exp \Bigl( -\frac{1}{\ve} \la  k(\theta), w\ra_{{\cal H}} -  \frac{1}{2\ve^2}  |k(\theta)|^2_{{\cal H}} 
\Bigr)
\tilde\chi^{\ve}_{\nu} ( w +  \frac{k(\theta)}{\ve})
G(\ve, w +  \frac{k(\theta)}{\ve})
 \eta_l (Y^{\ve, k(\theta)}_1 -a) 
 \nn
\\
&& 
\times
\det[-\ve  {\bf a}_{k(\theta)} (w) -  {\bf a}_{k(\theta)} (k(\theta))] 
 \delta_0 \bigl(\ve {\bf e}(\theta,w) + {\bf e}(\theta, k(\theta)) \bigr)
\Bigr]
\sqrt{ \det G(\theta)}
d\theta
\nn\\
&=&
\exp \bigl( - \frac{d_a^2}{2\ve^2} \bigr)
\int_{ \hat{O}_{\nu}} 
{\mathbb E}
\Bigl[
e^{ - \la  k(\theta), w\ra_{{\cal H}} /\ve }
\tilde\chi^{\ve}_{\nu} ( w +  \frac{k(\theta)}{\ve})
G(\ve, w +  \frac{k(\theta)}{\ve})
\nn\\
&& 
\times
 \eta_l ( R^{\ve}_{1} (w; k(\theta))) 
 \det[ {\rm Id}_{n'} -\ve  {\bf a}_{k(\theta)} (w) ] 
\zeta^{\ve,h_{\nu}}_{ \sigma_{\nu} \gamma_{\nu}} (w + \frac{k(\theta) }{\ve})
 \delta_0 \bigl(\ve {\bf e}(\theta,w)  \bigr)
\Bigr]
\sqrt{ \det G(\theta)}
d\theta.
\label{226.1.eq}
\end{eqnarray}
Here, we used the following facts: (i)~ ${\bf e}(\theta, k(\theta)) =0$,
(ii)~ ${\bf a}_{k(\theta)} (k(\theta)) =-  {\rm Id}_{n'}$
by (\ref{akk.eq}),
and 
(iii)~ $\tilde\chi^{\ve}_{\nu} ( w +  k(\theta)/\ve)
=
\tilde\chi^{\ve}_{\nu} ( w +  k(\theta)/\ve) 
\zeta^{\ve,h_{\nu}}_{ \sigma_{\nu} \gamma_{\nu}} (w + k(\theta) /\ve)
$,
which is immediate from $f(s) = f(s) f(s/ \sigma_{\nu}^{4m})$ for all $s \in {\mathbb R}$,
due to (\ref{csigma.eq}).
(Recall that $f$, $\zeta^{\ve,h}_{\gamma}$  and $\tilde\chi^{\ve}_{\nu}$ were 
 introduced at the beginning of Section 6.)


Now we apply Proposition \ref{pr.comp1} with 
$F (w)= ( R^{\ve}_{1} (w; k(\theta))/\ve ,{\bf e}(\theta,w) )$, 
$\chi = f$,
\[
\xi (w) =   \xi_{\ve}^{(\nu)} (w) :=\frac{
\| \tau_{k(\theta)-h_{\nu} }(\ve {\bf W})^1 \|_{\alpha, 4m -B}^{4m} 
+ \| \tau_{k(\theta)-h_{\nu} }(\ve {\bf W})^2 \|_{2\alpha, 2m -B}^{2m}
}
{\sigma_{\nu}^{4m} \gamma_{\nu}^{4m}}.
\]
Note that $\xi (w)$ is just a polynomial of order $4m$
and $\chi (\xi) = \zeta^{\ve,h_{\nu}}_{ \sigma_{\nu} \gamma_{\nu}} (w + k(\theta)/\ve)$, 
which is on the right hand side of (\ref{226.1.eq}).
As we have seen in (\ref{k_h_gam.eq}) and (\ref{223_3.cond}),
$|\xi (w)| \le 2$ implies that 
the smallest eigenvalue of $\sigma [F]$ is greater than or equal to 
$(\lambda \wedge 1)/2$.
Hence, we can use Proposition \ref{pr.comp1}
 to obtain that, for each $0< \ve \le \ve_0$,
\begin{eqnarray*}
&&\lim_{l \to \infty} 
 \zeta^{\ve,h_{\nu}}_{\sigma_{\nu}\gamma_{\nu}} (w + \frac{ k(\theta) }{\ve})
 \eta_l ( R^{\ve}_{1} (w; k(\theta)))  \delta_0 \bigl(\ve {\bf e}(\theta,w)  \bigr)
\\
&& \quad=
 \zeta^{\ve,h_{\nu}}_{\sigma_{\nu}\gamma_{\nu}} (w + \frac{ k(\theta)}{\ve})
  \delta_0 \bigl( R^{\ve}_{1} (w; k(\theta)),\ve {\bf e}(\theta,w)  \bigr)
   \qquad \mbox{in $\tilde{\mathbb D}_{- \infty}$ }
   \end{eqnarray*}   
uniformly in $\theta \in \hat{O}_{\nu}$.
Here, the delta function on the right hand side is defined on ${\cal V} \times {\mathbb R}^{n'}$.


Therefore, we have 
\begin{eqnarray}
\lefteqn{
{\mathbb E}[\tilde\chi^{\ve}_{\nu} G^{\ve}  \delta_a (Y^{\ve}_1) ]
}\nn\\
&=&
e^{- d_a^2 /2\ve^2 }
\int_{ \hat{O}_{\nu}} 
{\mathbb E}
\Bigl[
e^{ - \la  k(\theta), w\ra_{{\cal H}} /\ve }
\tilde\chi^{\ve}_{\nu} ( w +  \frac{k(\theta)}{\ve})
G(\ve, w +  \frac{k(\theta)}{\ve})
\nn\\
&& 
\quad \times
 \det[ {\rm Id} -\ve  {\bf a}_{k(\theta)} (w) ] 
\zeta^{\ve,h_{\nu}}_{\sigma_{\nu} \gamma_{\nu}} (w + \frac{k(\theta) }{\ve})
 \delta_0 \bigl(  R^{\ve}_{1} (w; k(\theta)) ,\ve {\bf e}(\theta,w)  \bigr)
\Bigr]
\sqrt{ \det G(\theta)}
d\theta
\nn\\
&=&
\ve^{-(n+n')} e^{- d_a^2 /2\ve^2 }
\int_{ \hat{O}_{\nu}} 
{\mathbb E}
\Bigl[
e^{ - \la  k(\theta), w\ra_{{\cal H}} /\ve }
\tilde\chi^{\ve}_{\nu} ( w +  \frac{k(\theta)}{\ve})
G(\ve, w +  \frac{k(\theta)}{\ve})
\nn\\
&& 
\quad \times
 \det[ {\rm Id} -\ve  {\bf a}_{k(\theta)} (w) ] 
\zeta^{\ve,h_{\nu}}_{\sigma_{\nu} \gamma_{\nu}} (w + \frac{k(\theta) }{\ve})
 \delta_0 \Bigl( \frac{ R^{\ve}_{1} (w; k(\theta)) }{\ve}, {\bf e}(\theta,w)  \Bigr)
\Bigr]
\sqrt{ \det G(\theta)}
d\theta
\nn\\
&=&
\ve^{-(n+n' )} e^{- d_a^2 /2\ve^2 }
\int_{ O_{\nu}} 
{\mathbb E}
\Bigl[
e^{ - \la  k, w\ra_{{\cal H}} /\ve }
\tilde\chi^{\ve}_{\nu} ( w +  \frac{k}{\ve})
G(\ve, w +  \frac{k}{\ve})
\nn\\
&& 
\quad \times
 \det[ {\rm Id} -\ve  {\bf a}_{k} (w) ] 
\zeta^{\ve,h_{\nu}}_{ \sigma_{\nu}\gamma_{\nu}} (w + \frac{k }{\ve})
 \delta_0 \Bigl( \frac{ R^{\ve}_{1} (w; k) }{\ve}, {\bf i}_k \la w\ra  \Bigr)
\Bigr]
\omega(dk),
\nn
\end{eqnarray}
where in the last line $\delta_0$ stands for the delta function on ${\cal V} \times T_k {\cal K}^{min}_a$.

Under the condition that $R^{\ve}_{1} (w; k) =0$,
$\la k,w \ra_{{\cal H}} = \la q(k), D\psi_1 (k) \la w\ra \ra_{{\cal V}}
 = \la q(k), g_1 (w;k) \ra_{{\cal V}}$
is equal to 
$- \la q(k), R^{\ve}_2 (w;k) \ra_{{\cal V}} /\ve$.

Since 
$f ((|y|^2 +|z|^2) /\kappa^2) \delta_0 (y,z) = \delta_0 (y, z)$ 
as a distribution on ${\cal V} \times T_k {\cal K}^{min}_a$ for any $\kappa >0$,
where $(y, z)$ denotes a generic element in ${\cal V} \times T_k {\cal K}^{min}_a$,
we can see that 
\begin{eqnarray*}
\zeta^{\ve,h_{\nu}}_{ \sigma_{\nu}\gamma_{\nu}} (w + \frac{k}{\ve})
f \Bigl( \frac{\kappa^{\ve} (w;k) }{\kappa^2}   \Bigr)
\delta_0 \Bigl( \frac{ R^{\ve}_{1} (w; k) }{\ve}, {\bf i}_k \la w\ra  \Bigr)
=
\zeta^{\ve,h_{\nu}}_{\sigma_{\nu}\gamma_{\nu}} (w + \frac{k}{\ve})
\delta_0 \Bigl( \frac{ R^{\ve}_{1} (w; k) }{\ve}, {\bf i}_k \la w\ra  \Bigr)
\end{eqnarray*}
in $\tilde{\mathbb D}_{- \infty}$, 
where we set $\kappa^{\ve} (w;k) =   
|R^{\ve}_{1} (w; k) /\ve|^2  + | {\bf i}_k \la w\ra  |^2$ for simplicity.
%
%


Summing up these all, we have the following lemma:
\begin{lm}\label{lm.tw67}
For any $1 \le \nu \le N$, $0 \le \ve \le \ve_0$, $\kappa >0$, it holds that
\begin{eqnarray}
\lefteqn{
{\mathbb E}[\tilde\chi^{\ve}_{\nu} G^{\ve}  \delta_a (Y^{\ve}_1) ]
}\nn\\
&=&
\ve^{-(n+n')} e^{- d_a^2 /2\ve^2 }
\int_{ O_{\nu}} 
{\mathbb E}
\Bigl[
e^{  \la  q(k),R^{\ve}_2 (w;k) /\ve^2 \ra_{{\cal V}}  }
\tilde\chi^{\ve}_{\nu} ( w +  \frac{k}{\ve})
G(\ve, w +  \frac{k}{\ve})
\nn\\
&& 
\quad \times
 \det[ {\rm Id} -\ve  {\bf a}_{k} (w) ] 
f \Bigl( \frac{\kappa^{\ve} (w;k) }{\kappa^2}   \Bigr)
\zeta^{\ve,h_{\nu}}_{\sigma_{\nu} \gamma_{\nu}} (w + \frac{k }{\ve})
 \delta_0 \Bigl( \frac{ R^{\ve}_{1} (w; k) }{\ve}, {\bf i}_k \la w\ra  \Bigr)
\Bigr]
\omega(dk).
\label{tw645.eq}
\end{eqnarray}
Here, $\ve_0 \in (0,1]$ is the constant given in Lemma \ref{lm.tw63}.
\end{lm}


Now it suffices to compute a full
asymptotic expansion of the generalized expectation in Lemma \ref{lm.tw67} above.
We will write $A_{\nu} (\ve, k) := {\mathbb E}[\cdots]$ on the right hand side 
in (\ref{tw645.eq}).
Note that all the asymptotics below are {\it uniform in $k \in {\cal K}_a^{min}$.}
So, we do not explicitly write uniformity or boundedness in $k$.

Now we list the asymptotic expansion
of each factor in (\ref{tw645.eq}). 
First, by Assumption {\bf (C1)}, 
\begin{equation}\label{asympG.eq}
G(\ve, w +  \frac{k}{\ve})
\sim
\Gamma_0 (k) +\ve \Gamma_1 (w;k)  +\ve^2 \Gamma_2 (w;k) + \cdots 
\qquad
\mbox{in ${\mathbb D}_{ \infty}$ as $\ve \searrow 0$.}
\end{equation}
Next,  $\det[ {\rm Id} -\ve  {\bf a}_{k} (w) ]$ is just a polynomial:
\begin{equation}\label{asympdet.eq}
\det[ {\rm Id} -\ve  {\bf a}_{k} (w) ]
=
{\cal J}_0 (w;k)  + \ve {\cal J}_1 (w;k) + \cdots  +\ve^{n'} {\cal J}_{n'} (w;k)
\end{equation}
for some ${\cal J}_j (w;k) \in {\mathbb D}_{ \infty}$  $(1 \le j \le n')$.
Note that 
${\cal J}_0 (w;k) \equiv 1$.


From (\ref{shift1.ineq}), (\ref{shift2.ineq}) and the definition of $\sigma_{\nu}$ in (\ref{csigma.eq}),
we can easily see that 
\[
 0 \le \lim_{\ve \searrow 0}  \xi_{\ve} (w)
=
\frac{
\| {\cal L}(k -h_{\nu} )^1 \|_{\alpha, 4m -B}^{4m} 
+ \|{\cal L}(k -h_{\nu} )^2 \|_{2\alpha, 2m -B}^{2m}
}
{\sigma_{\nu}^{4m} \gamma_{\nu}^{4m}}
\le 
\frac18,
\]
which is one of the assumptions in Proposition \ref{pr.comp2}.
It is clear that 
\[
(  R^{\ve}_{1} (w; k) /\ve, {\bf i}_k \la w\ra  )
\sim 
( g_1 (w;k) , {\bf i}_k \la w\ra)
+ 
\ve 
( g_2 (w;k) , 0)
%
+\cdots
\qquad
\mbox{in ${\mathbb D}_{ \infty} ({\mathbb R}^{n+n' })$ as $\ve \searrow 0$.}
\]

Hence, we can use Proposition \ref{pr.comp2}  to obtain
\begin{eqnarray}\label{asympdel.eq}
\lefteqn{
\zeta^{\ve,h_{\nu}}_{\sigma_{\nu} \gamma_{\nu}} (w + \frac{k }{\ve})
 \delta_0 \Bigl( \frac{ R^{\ve}_{1} (w; k) }{\ve}, {\bf i}_k \la w\ra  \Bigr)
}
\nn\\
&&
\sim
\Phi_0 (w; k) +\ve \Phi_1 (w;k)  +\ve^2 \Phi_2 (w;k) + \cdots 
\qquad
\mbox{in $\tilde{\mathbb D}_{ - \infty}$ as $\ve \searrow 0$.}
\end{eqnarray}
The coefficient
$\Phi_j (w;k)$ is obtained by the formal asymptotic expansion 
for the composition of $ \delta_0$ and $(  R^{\ve}_{1} (w; k) /\ve, {\bf i}_k \la w\ra  )$
as explained in Proposition \ref{pr.comp2}.
In particular, $\Phi_0 (w;k) = \delta_0 ( g_1 (w;k) , {\bf i}_k \la w\ra)$.


Let $L \in {\mathbb N}$ be arbitrarily large.
Since $e^{  \la  q(k),R^{\ve}_2 (w;k) /\ve^2 \ra_{{\cal V}}  }
\tilde\chi^{\ve}_{\nu} ( w +  k/\ve)
f ( \kappa^{\ve} (w;k) /\kappa^2   )
$
is bounded 
in  $\tilde{\mathbb D}_{\infty}$ due to Lemma \ref{lm.tw65},
we have 
\begin{eqnarray}\label{new1.eq}
A_{\nu} (\ve, k) 
&=&
{\mathbb E}
\Bigl[
e^{  \la  q(k),R^{\ve}_2 (w;k) /\ve^2 \ra_{{\cal V}}  }
\tilde\chi^{\ve}_{\nu} ( w +  \frac{k}{\ve})
G(\ve, w +  \frac{k}{\ve})
\nn\\
&& 
\quad \times
 \det[ {\rm Id} -\ve  {\bf a}_{k} (w) ] 
f \Bigl( \frac{\kappa^{\ve} (w;k) }{\kappa^2}   \Bigr)
\Bigl(
\sum_{j=0}^{L} \ve^j \Phi_j (w; k)
\Bigr)
\Bigr]
+O(\ve^{L+1})
\end{eqnarray}
as $\ve \searrow 0$.


Then, we  expand 
$e^{  \la  q(k),R^{\ve}_2 (w;k) /\ve^2 \ra_{{\cal V}}  }
\tilde\chi^{\ve}_{\nu} ( w +  k/\ve)
f ( \kappa^{\ve} (w;k) /\kappa^2   )
$
in  $\tilde{\mathbb D}_{\infty}$
by using Lemma \ref{lm.tw65}.
It is clear that 
\[
\frac{R^{\ve}_{3} (w; k)}{\ve^2}
= 
\frac{R^{\ve}_{2} (w; k)}{\ve^2} -g_2 (w;k) 
\sim
 \ve g_3 (w;k) + \ve^2 g_4 (w;k)
 +\cdots
\qquad
\mbox{in ${\mathbb D}_{ \infty} ({\mathbb R}^{n})$ as $\ve \searrow 0$.} 
 \]
In particular, $R^{\ve}_{3} (w; k)/\ve^2 =O(\ve)$ 
in ${\mathbb D}_{ \infty} ({\mathbb R}^{n})$ as $\ve\searrow 0$.
From the Taylor expansion of $e^x$ at $x=0$, we have
\begin{eqnarray}
\lefteqn{
\tilde\chi^{\ve}_{\nu} ( w +  \frac{k}{\ve})
f ( \frac{\kappa^{\ve} (w;k) }{\kappa^2 }  )
e^{  \la  q(k),R^{\ve}_2 (w;k) /\ve^2 \ra_{{\cal V}}  }
}
\nn\\
&=&
\tilde\chi^{\ve}_{\nu} ( w +  \frac{k}{\ve})
f ( \frac{\kappa^{\ve} (w;k) }{\kappa^2 }  )
\Bigl\{
e^{  \la  q(k), g_2 (w;k)  \ra_{{\cal V}}  }
\sum_{j=0}^{L} \frac{ \la  q(k),R^{\ve}_3 (w;k) /\ve^2 \ra_{{\cal V}}^j }{j!}
\nn\\
&&
+ 
\int_0^1 
d\tau
\exp \bigl(   \la  q(k),  (1-\tau)g_2 (w;k) + \tau \frac{R^{\ve}_2 (w;k)}{\ve^2} 
 \ra_{{\cal V}}  \bigr)
\cdot
\frac{ \la  q(k),R^{\ve}_3 (w;k) /\ve^2 \ra_{{\cal V}}^{L+1} }{(L+1)!}
\Bigr\}.
\label{tks_key.eq}
\end{eqnarray}
Using Lemma  \ref{lm.tw65}, (\ref{lm.tw64.suff}) in its proof
and H\"older's inequality with 
the exponents $1/\tau$ and $1/(1 -\tau)$, 
we can show that 
the last term on the right hand side of (\ref{tks_key.eq}) 
is $O(\ve^{L+1})$ in $\tilde{\mathbb D}_{ \infty}$-topology.
Here, we have used (\ref{223_4.cond}). 
Similarly, 
$e^{  \la  q(k), g_2 (w;k)  \ra_{{\cal V}}  }
\tilde\chi^{\ve}_{\nu} ( w +  k/\ve)
f ( \kappa^{\ve} (w;k) /\kappa^2   )
\in \tilde{\mathbb D}_{\infty}$, due to (\ref{lm.tw64.suff}).

Denote by ${\cal B}_j (w;k) \in {\mathbb D}_{ \infty} ~(j \in {\mathbb N})$ 
the coefficients
that appear in the formal asymptotic expansion of  
\[
\sum_{j=0}^{\infty} \frac{ \la  q(k),R^{\ve}_3 (w;k) /\ve^2 \ra_{{\cal V}}^j }{j!}
\sim 
{\cal B}_0 (w;k)
+
\ve{\cal B}_1 (w;k)
+\ve^2
{\cal B}_2 (w;k) + \cdots
\quad
\mbox{in ${\mathbb D}_{ \infty}$ as $\ve \searrow 0$.}
\]
Note that ${\cal B}_0 (w;k) \equiv 1$.
As $\ve \searrow 0$, it holds that
\begin{eqnarray}
\lefteqn{
\tilde\chi^{\ve}_{\nu} ( w +  \frac{k}{\ve})
f ( \frac{\kappa^{\ve} (w;k) }{\kappa^2 }  )
e^{  \la  q(k),R^{\ve}_2 (w;k) /\ve^2 \ra_{{\cal V}}  }
}\nn\\
&=&
\tilde\chi^{\ve}_{\nu} ( w +  \frac{k}{\ve})
f ( \frac{\kappa^{\ve} (w;k) }{\kappa^2 }  )
e^{  \la  q(k), g_2 (w;k)  \ra_{{\cal V}}  }
\sum_{j=0}^{L} 
\ve^j
{\cal B}_j (w;k)
+ 
O(\ve^{L+1})
\quad
\mbox{in $\tilde{\mathbb D}_{ \infty}$} 
\label{2tks_key.eq}
\end{eqnarray}
%
%
and therefore 
\begin{eqnarray}\label{new2.eq}
A_{\nu} (\ve, k) 
&=&
{\mathbb E}
\Bigl[
e^{  \la  q(k), g_2 (w;k)  \ra_{{\cal V}}  }
\tilde\chi^{\ve}_{\nu} ( w +  \frac{k}{\ve})
G(\ve, w +  \frac{k}{\ve})
 \det[ {\rm Id} -\ve  {\bf a}_{k} (w) ] 
 \nn\\
&& 
\quad \times
f \Bigl( \frac{\kappa^{\ve} (w;k) }{\kappa^2}   \Bigr)
\Bigl(
\sum_{j_1=0}^{L} \ve^{j_1}  {\cal B}_{j_1} (w; k)
\Bigr)
\Bigl(
\sum_{j_2=0}^{L} \ve^{j_2} \Phi_{j_2} (w; k)
\Bigr)
\Bigr]
+O(\ve^{L+1}).
\end{eqnarray}

By Proposition \ref{pr.comp2}, 
each $\Phi_j (w;k)$ is a finite sum 
of  terms of the following form;
\[
\mbox{
(a ${\mathbb D}_{\infty}$-functional)
$\times$
$(\partial^{\beta} \delta_0) ( g_{1} (w; k), {\bf i}_k \la w\ra)$,
}
\]
where $\beta$ is a certain
$(n+ n')$-dimensional multi-index.
Hence, by Corollary \ref{co.newadd}, {\rm (iii)} 
and Remark \ref{re.add},
$\exp ( \la  q(k), g_2 (w;k)  \ra_{{\cal V}}  ) \Phi_j (w;k) \in {\mathbb D}_{- \infty}$.
Note that all
the other factors in the generalized expectation
in (\ref{new2.eq}) belong to ${\mathbb D}_{\infty}$.
(Because $\exp ( \la  q(k), g_2 (w;k)  \ra_{{\cal V}}  ) \notin 
\tilde{\mathbb D}_{\infty}$, a little  care was needed above.)
So, it suffices to expand them in ${\mathbb D}_{\infty}$-topology.


By the general theory, $f (\kappa^{\ve} (w;k)/\kappa^2   )$ 
admits asymptotic expansion in ${\mathbb D}_{ \infty}$-topology as follows:

\begin{eqnarray}\label{asympf.eq}
f ( \frac{\kappa^{\ve} (w;k) }{\kappa^2 }  )
&=&
f\Bigl(
\frac{  |R^{\ve}_{1} (w; k) /\ve|^2  + | {\bf i}_k \la w\ra  |^2 }{\kappa^2 }
\Bigr)
\nn\\
&\sim& 
{\cal C}_0 (w;k)
+
\ve{\cal C}_1 (w;k)
+\ve^2
{\cal C}_2 (w;k)+ \cdots
\quad
\mbox{in ${\mathbb D}_{ \infty}$ as $\ve \searrow 0$.} 
\end{eqnarray}
Note that 
${\cal C}_0 = f (\{|g_{1} (w; k) |^2  + | {\bf i}_k \la w\ra  |^2 \}/\kappa^2)$.

Let us see that $\{{\cal C}_i (w;k)\}$ does not contribute to 
the asymptotic expansion.
In other words, 
$f (\kappa^{\ve} (w;k) /\kappa^2)$ is merely a dummy factor which is introduced for technical reasons.
Each term in 
${\cal C}_i (w;k)$ has a factor of the form 
$f^{(l)} ( \{| g_{1} (w; k) |^2 +| {\bf i}_k \la w\ra|^2  \} /\kappa^2)$
for some $l \in {\mathbb N}$.
Similarly, by Proposition \ref{pr.comp2}, 
each term of
$\Phi_j (w;k)$ has a factor of the form 
$(\partial^{\beta} \delta_0) ( g_{1} (w; k), {\bf i}_k \la w\ra)$
for some 
$(n+ n')$-dimensional multi-index $\beta$.
Noting that
\begin{eqnarray*}
f^{(l)} \Bigl( \frac{| g_{1} (w; k) |^2 +| {\bf i}_k \la w\ra|^2  } {\kappa^2} \Bigr)
\cdot
(\partial^{\beta} \delta_0) ( g_{1} (w; k), {\bf i}_k \la w\ra)
\nn\\
=
\begin{cases}
    0 & (\mbox{if $l \ge 1$}), \\
    (\partial^{\beta} \delta_0) (g_{1} (w; k) , {\bf i}_k \la w\ra)
  & (\mbox{if $l =0$}),
  \end{cases}
\end{eqnarray*}
we can easily see that
${\cal C}_i (w;k)  \Phi_j (w;k) =0$ in $\tilde{\mathbb D}_{ - \infty}$ if $i \ge 1$
and hence
\[
\sum_{i=0}^{\infty} \ve^i {\cal C}_i (w;k) \sum_{j=0}^{\infty}  \ve^j \Phi_j (w;k)
=
\sum_{j=0}^{\infty}\ve^j  \Phi_j (w;k)
\]
holds as a formal asymptotic expansion.


Again by the general theory, $\tilde\chi^{\ve}_{\nu} ( w +  k /\ve)$ also
admits an asymptotic expansion in ${\mathbb D}_{ \infty}$:
\begin{eqnarray}\label{asympf2.eq}
\tilde\chi^{\ve}_{\nu} ( w +  \frac{k}{\ve})
&\sim& 
{\cal D}_0^{(\nu)} (w;k)
+
\ve{\cal D}_1^{(\nu)} (w;k)
+\ve^2
{\cal D}_2^{(\nu)} (w;k)+ \cdots
\quad
\mbox{in ${\mathbb D}_{ \infty}$ as $\ve \searrow 0$.} 
\end{eqnarray}
However, we should note here that the coefficient $\{ {\cal D}_i^{(\nu)} (w;k)\}$
may depend on $\nu$.
If $k \notin O_{\nu}$, then $k \notin U_{h_{\nu}, \gamma_{\nu}}$,
which implies that ${\cal D}_i^{(\nu)} (w;k) \equiv 0$ for all $i \ge 0$.
Moreover, since 
 $\sum_{\nu=1}^N \tilde{\chi}^{\ve}_{\nu}=  \chi^{\ve}$
 and 
 $\chi^{\ve}$ admits the asymptotic expansion  (\ref{asy_chi_ve}),
  we easily see that, for all $k \in  {\cal K}^{min}_a$,  
   $\sum_{\nu} {\cal D}_i^{(\nu)} (w;k) \equiv 0$ for all $i \ge 1$
  and $\sum_{\nu} {\cal D}_0^{(\nu)} (w;k) \equiv 1$.


We now set a notation.  
Define ${\cal E}_i (w;k) \in \tilde{\mathbb D}_{- \infty}$
by the following formal asymptotic series: 
\begin{eqnarray}
\lefteqn{
\sum_{j=0}^{\infty} \ve^j {\cal E}_j (w;k) 
}
\nn\\
&=&
\Bigl(
\sum_{j_1=0}^{\infty} \ve^{j_1} \Gamma_{j_1} (w;k) 
\Bigr)
\Bigl(
\sum_{j_2=0}^{m} \ve^{j_2} {\cal J}_{j_2} (w;k) 
\Bigr)
\Bigl(
\sum_{j_3=0}^{\infty} \ve^{j_3} {\cal B}_{j_3} (w;k) 
\Bigr)
\Bigl(
\sum_{j_4=0}^{\infty} \ve^{j_4} \Phi_{j_4} (w;k) 
\Bigr).
\label{calE.def}
\end{eqnarray}
Note that $\Gamma_{j} (w;k),  {\cal J}_{j} (w;k),  {\cal B}_{j} (w;k), \Phi_{j} (w;k)$ 
are all defined not just for $k \in O_{\nu}$, but for any $k \in {\cal K}^{min}_a$.
Hence, so is ${\cal E}_{j} (w;k)$, which means that
the definition of ${\cal E}_{j} (w;k)$ is independent of $\nu$.
It is easy to see that 
${\cal E}_0 (w;k) = \Gamma_0 (k) \delta_0 ( g_1 (w;k) , {\bf i}_k \la w\ra)$.


Combining these all, we have the following asymptotics: for any $L \in {\mathbb N}$,
\begin{eqnarray}\label{new3.eq}
A_{\nu} (\ve, k) 
&=&
{\mathbb E}
\Bigl[
e^{  \la  q(k), g_2 (w;k)  \ra_{{\cal V}}  }
\Bigl(
\sum_{j=0}^{L} \ve^{j} {\cal E}_{j} (w;k) 
\Bigr)
\Bigl(
\sum_{l=0}^{L} \ve^{l} {\cal D}_{l}^{(\nu)} (w;k) 
\Bigr)
\Bigr]
+O(\ve^{L+1})
\end{eqnarray}
as $\ve \searrow 0$ uniformly in $k \in {\cal K}_a^{min}$.
Therefore, 
\begin{eqnarray}
\lefteqn{
\sum_{\nu =1}^{N}
{\mathbb E}[\tilde\chi^{\ve}_{\nu} G^{\ve}  \delta_a (Y^{\ve}_1) ]
}
\nn\\
&=&
\sum_{\nu =1}^{N}
\Bigl\{
\int_{ O_{\nu}} 
{\mathbb E}
\Bigl[
e^{  \la  q(k), g_2 (w;k)  \ra_{{\cal V}}  }
\Bigl(
\sum_{j=0}^{L} \ve^{j} {\cal E}_{j} (w;k) 
\Bigr)
\Bigl(
\sum_{l=0}^{L} \ve^{l} {\cal D}_{l}^{(\nu)} (w;k) 
\Bigr)
\Bigr]
\omega(dk)
+O (\ve^{L+1})
\Bigr\}
\nn
\\
&=&
\sum_{\nu =1}^{N}
\Bigl\{
\int_{   {\cal K}^{min}_a } 
{\mathbb E}
\Bigl[
e^{  \la  q(k), g_2 (w;k)  \ra_{{\cal V}}  }
\Bigl(
\sum_{j=0}^{L} \ve^{j} {\cal E}_{j} (w;k) 
\Bigr)
\Bigl(
\sum_{l=0}^{L} \ve^{l} {\cal D}_{l}^{(\nu)} (w;k) 
\Bigr)
\Bigr]
\omega(dk)
+O (\ve^{L+1})
\Bigr\}
\nn\\
&=&
\sum_{j=0}^{L}
\ve^j
\int_{     {\cal K}^{min}_a} 
{\mathbb E}
\Bigl[
e^{  \la  q(k), g_2 (w;k)  \ra_{{\cal V}}  }
 {\cal E}_{j} (w;k) 
\Bigr]
\omega(dk)
+O (\ve^{L+1})
\label{asFin.eq}
\end{eqnarray}
as $\ve \searrow 0$, which is the desired asymptotic expansion.
Note that $\{{\cal D}_i^{(\nu)} (w;k)\}$ does not contribute to 
the asymptotic expansion, either.
The leading term is given by
\begin{equation}\label{tp_trm.eq}
c_0 =
\int_{     {\cal K}^{min}_a} 
{\mathbb E}
\Bigl[
e^{  \la  q(k), g_2 (w;k)  \ra_{{\cal V}}  }
 \delta_0 ( g_1 (w;k) , {\bf i}_k \la w\ra)
 \Bigr]
\Gamma_0 (k)\omega(dk),
\end{equation}
which is positive 
if $\Gamma_0 $ is non-negative, but not identically zero
on ${\cal K}_a^{min}$.

To complete the proof of our main theorem (Theorem \ref{tm.main}),
we show that $c_{2j +1} =0~(j=0,1,2,\ldots)$ under {\bf (C2)}.
First, note that 
$f_{2j } (w;h)$ and $f_{2j +1} (w;h)$ are even and odd in $w$, respectively.
(For simplicity, we will say that $\{f_j\}$ satisfies the even-odd property.)
Obviously, so does $\{g_j\}$.
(See  (\ref{asympX.1}) and  (\ref{asympX.2}) for the 
definitions of  $\{f_j\}$ and $\{g_j\}$.)
It immediately follows from this that $\{ {\cal B}_j\}$ satisfies this property.
Although not so obvious, 
it is still straightforward to check the even-odd property of $\{ \Phi_j\}$.
From the explicit form of ${\bf a}_k (\ve w)$,  it is easy to see that
$\{ {\cal J}_j\}$ has the same property, too.
Therefore, 
$\{ {\cal E}_j\}$ satisfies the even-odd property
and
\[
{\mathbb E}
\Bigl[
e^{  \la  q(k), g_2 (w;k)  \ra_{{\cal V}}  }
 {\cal E}_{2j +1} (w; k)
  \Bigr]
 =0
\]
by the invariance of the Wiener measure $\mu$ and $g_2 (w;k)$
under $w \mapsto -w$.
(Although it is not necessary in the proof of this Euclidean case, 
we note that $\{ {\cal D}_j^{(\nu)}\}$ also satisfies this property
due to the explicit construction of $\tilde{\chi}^{\ve}_{\nu}$.)
Thus, we have shown Theorem \ref{tm.main}.

%
\section{Proof of the manifold case}
\label{sec.pr_mfd}

In this section we prove Theorem \ref{tm.main_mfd}.
As we mentioned before, 
we choose Riemannian metrics on ${\cal M}$ and ${\cal N}$ 
arbitrarily and fix them.
The measure  ${\rm vol}$ is the Riemannian measure on ${\cal N}$.
It is often useful to embed ${\cal M}$
into ${\mathbb R}^N$ with sufficiently high dimension $N$ 
and solve differential equations in the ambient space ${\mathbb R}^N$.
For instance, 
when we prove  continuity or differentiability 
of various maps from  (rough) path spaces,  
this extrinsic view has an advantage.

Choose such an embedding ${\cal M} \hookrightarrow {\mathbb R}^N$
and 
extend $V_i ~(0 \le i \le r)$ 
so that it becomes $C_b^{\infty}$-vector field on ${\mathbb R}^N$.
(By abusing notations
we denote the extended vector fields  by the same symbols.)
Then, consider Stratonovich SDE (\ref{sdeX.def})
and skeleton ODE (\ref{ode.def})
on ${\mathbb R}^N$.
If $x \in {\cal M}$, 
then the solutions stay in ${\cal M}$
and coincide with the solutions
 (\ref{sde.mfd})
and 
 (\ref{ode.mfd}), respectively.

The case of RDE is slightly different.
We do not define RDE on ${\cal M}$ in an intrinsic way.
We just solve RDE  (\ref{rde_x.def}) on ${\mathbb R}^N$ 
and call the first level path of the solution 
the solution of RDE on ${\cal M}$.
(Note that if $x \in {\cal M}$,
then
the first level path of the solution stays in ${\cal M}$.)
Since we only need the first level path 
of the solution of the RDE, this is enough for our purpose.

Recall that $J (t, x, h)$ defined by (\ref{ode_J.def}) on the ambient space
is the Jacobian of the smooth map 
$x \mapsto \phi (t, x, w)$.
Hence,
 $J (t, x, h)$ is a linear map from $T_x ({\mathbb R}^N)$
to $T_{\phi (t,x,h)} ({\mathbb R}^N)$.
When restricted to $T_x ({\cal M})$, 
it is from $T_x ({\cal M})$ to $T_{\phi (t,x,h)} ({\cal M})$
and coincides with the Jacobian of $x \mapsto \phi (t, x, w)$ defined by (\ref{ode.mfd})
on ${\cal M}$.
A similar fact holds for $X^{\ve} (t,x,w)$ and $J^{\ve} (t,x,w)$.


We remark that all the results in Section \ref{sec.skltn}
remain valid 
in the manifold setting with trivial modifications.
To see this, let us write key quantities in an intrinsic way.
As usual 
we will identify ${\cal H}$ and ${\cal H}^*$ via the Riesz isometry.
We denote by 
${\mathbb L} ({\cal X}, {\cal Y})$ the space of bounded linear maps from 
from ${\cal X}$ to ${\cal Y}$.

First, $D\phi_t (h) \in {\mathbb L} ({\cal H}, T_{\phi_t(h)} {\cal M})$
and 
$D\psi_t (h) =  
(\Pi_*)_{\phi_t (h)}  \circ D\phi_t (h)
\in {\mathbb L} ({\cal H}, T_{\psi_t(h)} {\cal N})$,
where 
$(\Pi_*)_{x} \colon T_{x} {\cal M} \to T_{\Pi (x)} {\cal N}$ is the tangent map of $\Pi$ 
at $x \in {\cal M}$.
The deterministic Malliavin covariance matrix is defined by
\begin{equation}\label{dmc_mfd}
\sigma [\psi_1] (h) := D\psi_1 (h)
\circ D\psi_1 (h)^* 
\in
{\mathbb L}
 ( T^*_{\psi_1 (h)} {\cal N}, T_{\psi_1  (h)} {\cal N}),
\end{equation}
where the  Riesz isometry ${\cal H} ={\cal H}^*$ is implicit.

Similarly,
$DX^{\ve}_t (w)
\in {\mathbb L} ({\cal H}, T_{X^{\ve}_t(w)} {\cal M})$
and 
$DY^{\ve}_t  (w)=  
(\Pi_*)_{X^{\ve}_t (w)}  \circ DX^{\ve}_t  (w)
\in {\mathbb L} ({\cal H}, T_{Y^{\ve}_t (w)} {\cal N})$ for a.a. $w$.
The Malliavin covariance matrix of $Y^{\ve}_1$ is as follows:
\begin{equation}\label{mc_mfd}
\sigma [Y^{\ve}_1 ] (w) := DY^{\ve}_1 (w)
\circ D Y^{\ve}_1 (w)^* 
\in
{\mathbb L}
 ( T^*_{Y^{\ve}_1 (w) } {\cal N}, T_{ Y^{\ve}_1 (w)} {\cal N}).
\end{equation}
The (deterministic) Malliavin covariance matrix is 
viewed as a bilinear form on the cotangent space.
Since the cotangent space is equipped with an inner product,
 its eigenvalues and determinant make sense.


Here, we discuss the integration by parts formula 
for manifold-valued smooth Wiener functionals.
This was implicitly used in \cite{ta}
and is not very difficult. 
However, since no proof is given in \cite{ta},
we now give a sketch of proof for the reader's convenience.
This formula is used in {\rm (i)}
pullback of a Schwartz distribution on the manifold 
by a non-degenerate 
manifold-valued smooth Wiener functional, 
i.e. a manifold version of Item {\bf (c)} in Subsection \ref{subsec.31}
and 
{\rm (ii)} upper estimate of the large deviation, i.e.
a manifold version of Theorem \ref{tm.ub_ldp}.

Leaving SDE (\ref{sde.mfd}) aside for a moment,
we discuss in a general setting.
Let ${\cal N}$ be a compact Riemannian manifold 
and let $F$ be a ${\mathbb D}_{\infty}$ -Wiener functional 
that takes values in ${\cal N}$.
(For basics of manifold-valued Malliavin calculus, see \cite{ta}.)
The sets of smooth functions
and of smooth vector fields on ${\cal N}$ are denoted by 
$C^{\infty} ({\cal N})$ and $\Gamma^\infty (T{\cal N})$, respectively.
In a standard way, we
identify $T_y {\cal N}$ and ${\cal H}$
with their dual spaces $T^*_y {\cal N}$ and ${\cal H}^*$, 
respectively.

Suppose that 
$D F(w) \circ D F(w)^* \colon
   T_{F(w)}{\cal N}\to T_{F(w)}{\cal N}$ is non-singular, i.e.
\begin{equation}\label{ass.int.smfd.01}
    \det [ D F(w) \circ D F(w)^* ] \ne0
    \quad\text{for a.a. }w\in{\cal W}.
\end{equation}
For $f \in C^{\infty} ({\cal N})$  and 
$Z \in \Gamma^\infty (T{\cal N})$, substituting
\[
    h= DF(w)^* \circ
        \{  D F(w) \circ D F(w)^* \}^{-1}Z_{F(w)}
\]
into 
\[
D(f(F))(w) \la h \ra = \la  D(f(F))(w), h \ra_{{\cal H}}
    =(\nabla f) (w)
       \bigl\la DF (w) \la h\ra \bigr\ra,
       \]
we obtain that
\begin{equation}\label{eq.int.smfd.02}
    \langle D(f(F)) (w) ,
    DF(w)^* \circ
        \{  D F(w) \circ D F(w)^*     \}^{-1} Z_{F(w)}
        \rangle_{{\cal H}}
    =(Zf)(F(w)).
\end{equation}

In addition, assume that
\begin{equation}\label{ass.int.smfd.02}
     (DF)^* \circ 
        \{ DF \circ ( DF )^*\}^{-1} Z_F
     \in{\mathbb D}_{\infty}({\cal H})
     \quad\text{for any }Z\in\Gamma^{\infty} (T{\cal N}).
\end{equation}
Setting 
\begin{align*}    
    & \Phi_1(Z^1;G)
       =D^*\Bigl[  G \cdot  [ (DF)^* \circ
        \{(DF) \circ (DF)^*\}^{-1} Z_F^1 ] \Bigr],
     \\
     & \Phi_k(Z^1,\dots,Z^k;G)
         =\Phi_1(Z^k;\Phi_{k-1}(Z^1,\dots,Z^{k-1};G)),
         \quad k=2,3,\dots
\end{align*}
for $Z^1,\dots,Z^k \in\Gamma^{\infty} (T{\cal N})$ 
and $G\in{\mathbb D}_{\infty}$, and then
using \eqref{eq.int.smfd.02} successively, we arrive at 
the desired integration by parts formula
for a manifold-valued Wiener functional:
\begin{equation}\label{ibp.mfd}
    {\mathbb E}[(Z^1\dots Z^k f)(F)G]
    ={\mathbb E}[f(F)\Phi_k(Z^1,\dots,Z^k;G)]
\end{equation}
for any $Z^1,\dots,Z^k \in\Gamma^{\infty} (T{\cal N})$ 
and $G\in{\mathbb D}_{\infty}$, 
$f\in C^\infty({\cal N})$.
For the same reason as in the Euclidean-valued case, 
this formula extends to the case when 
$f$ is a Schwartz distribution on ${\cal N}$.

The assumptions \eqref{ass.int.smfd.01} and \eqref{ass.int.smfd.02}
are satisfied provided that 
\[
    \|\{DF \circ (D F)^*\}^{-1}\|_{\text{op},F}
    \in \bigcap_{1<p <\infty } L^p  (\mu)
        \quad
    \mbox{or}
    \quad
    \det [ \{ DF \circ (D F)^*\}^{-1} ]  \in \bigcap_{1<p <\infty} L^p (\mu).
     \]
Here,  $\|\cdot\|_{\text{op},y}$ stands for 
the operator norm of 
a linear mapping of $T_y {\cal N}$ into itself with respect to the Riemannian
metric.
Since we  work on a compact manifold, 
the above condition is 
independent of the choice of a Riemannian metric.

%
By the way, there is another justification 
of Watanabe's composition $T \circ F$ on manifold.
It does not use (\ref{ibp.mfd}), but a localized version of Watanabe's 
composition theorem in (an open subset of) ${\mathbb R}^n$
(not Theorem \ref{pr.comp1} in this paper, but its 
slightly different variant in Yoshida \cite{yo}).
Here, we only give a sketch.
If a distribution $T$ has a very small support (like the delta functions),
then its support is contained in one coordinate chart.
Then we can find a suitable real-valued
 cutoff function on ${\cal W}$ and view that $F$,
restricted to the support of the cutoff function, 
takes values in an open subset of ${\mathbb R}^n$. 
Then, we can use the localized version of Watanabe's 
composition theorem to justify $T \circ F$.
Finally, recall that a general distribution on ${\cal N}$ can be written 
as a sum of distributions with small support 
thanks to a partition of unity on ${\cal N}$.
%


Let us get back to our original setting and 
assume {\bf (B1)}.
For $h \in {\cal K}_a^{min}$,
the Lagrange multiplier 
$
q(h)
=\sigma [\psi_1] (h)^{-1} D\psi_1 (h) \la h\ra
$ 
should be regarded as an element of 
$T^*_a {\cal N}$. 
Hence, 
$q(h)\circ  (\Pi_*)_{ \phi_1 (h)} \in T^*_{\phi_1 (h)} {\cal M}$
and
$q(h)\circ  (\Pi_*)_{ \phi_1 (h)}  \circ J_1(h) \in T^*_{x} {\cal M}$.
Recall the latter is the initial value $p_0$ of the Hamiltonian ODE
associated with Hamiltonian 
$$H(x,p) 
= \frac12 \sum_{i=1}^r  \bigl\la p, V_i (x)\bigr\ra^2
\qquad
\quad
(x,p) \in T^* {\cal M}= \cup_{ z \in {\cal M}} T^*_z {\cal M},
$$
where the pairing is between 
$T^*_{x} {\cal M}$ and  $T_{x} {\cal M}$.
Clearly, $H$ is a smooth function on the cotangent bundle $T^* {\cal M}$.

Then, if Assumption {\bf (A1)} is replaced by {\bf (A1)'},
Proposition \ref{pr.hml},
Remark \ref{re.rglr},
Proposition \ref{pr.smth},
Proposition \ref{pr.top=} and
Lemma \ref{lm.tw_62}
also hold 
in our manifold setting with trivial modifications of the statements.
(For instance, 
$\Pi_{{\cal V}}$ should be  replaced by $ (\Pi_*)_{ \phi_1 (h)}$, etc.)
A rough sketch of proof is as follows:
If we take a local coordinate chart on ${\cal M}$,
 ODEs for $\phi(h)$ and $J(h)$ can be written down in the coordinates.
Then, the computation for the Hamiltonian ODE is essentially the same 
as in Section \ref{sec.skltn}.
To prove differentiability or (rough path) continuity
of $q(h)$ or $p_0 = p_0(h)$, etc. in $h$, the extrinsic expression of
 the skeleton ODE on ${\mathbb R}^N$ is better.


As we will see in the next lemma,
in a neighborhood of the path $t \mapsto  \phi_t (h)$, 
$h \in {\cal K}_a^{min}$,
the 
RDE on ${\cal M}$ can be transferred to  an RDE on ${\mathbb R}^d$.
Hence, for  local analysis around $h$ of the RDE on manifold,
it is sufficient to deal with the corresponding RDE on the Euclidean space. 
This is one of  main advantages of  using rough path theory
since in the classical theory of It\^o maps
this kind of  ``local" operation is not easy and sometimes impossible.
\begin{lm}\label{lm.path_nbh}
 {\rm (i)}~ For any $h \in {\cal K}_a^{min}$, 
 we can find the following {\rm (1)}--{\rm (4)}:
\\
 {\rm (1)}~
An open neighborhood ${\cal U} = {\cal U}_h$ 
of ${\rm Im} \phi (h) :=\{\phi_t (h) \mid t \in [0,1] \}$ in ${\cal M}$.
\\
 {\rm (2)}~
A bounded,  connected, open subset ${\cal U}^{\prime}={\cal U}^{\prime}_h$
 in ${\mathbb R}^d$.
\\
 {\rm (3)}~
A diffeomorphism $\theta= \theta_h \colon {\cal U} \to {\cal U}^{\prime}$.
\\
 {\rm (4)}~Vector fields $V_i^{(h)}$ on ${\mathbb R}^d$ of class $C^{\infty}_b$
 such that $(\theta_*)_x V_i (x) = V_i^{(h)} (\theta (x))$ ($0 \le i \le r$).
 \\
 {\rm (ii)}~Consider 
 skeleton ODE (\ref{ode.def}) 
 and  RDE (\ref{rde_x.def}) on ${\mathbb R}^d$, 
but with the coefficients $V_i^{(h)} ~(0 \le i \le r)$ and the initial point $\theta (x)$. 
   If $k \in {\cal H}$ is sufficiently close to $h$ in ${\cal H}$-topology,
   then $t \mapsto \theta (\phi_t (k))$ solves the skeleton ODE on ${\mathbb R}^d$
   driven by $k$.  
      Likewise, if ${\bf w}$ and ${\bf h}={\cal L} (h)$ are sufficiently close in 
      $G\Omega_{\alpha, 4m}^{B}({\mathbb R}^d)$ 
      and $\| \lambda\|_{1 -H}$ is sufficiently small, 
      then $t \mapsto \theta (x + {\bf x}^1_{0,t} ) = \theta (x + \Phi ({\bf w}, \lambda)^1_{0,t} )$
      is the first level path of the solution of the RDE on ${\mathbb R}^d$
      driven by 
      $({\bf w},\lambda)$.
              \end{lm}

\Proof
First, we prove {\rm (i)}.
Clearly, 
the image of the path $\phi (h)$ is compact in ${\cal M}$
and there is  no self-intersection due to the energy-minimizing condition.
Moreover, by Remark \ref{re.rglr},  $\phi_t^{\prime} (h)$ never vanishes.
Therefore, we can use the standard technique 
of the tubular neighborhood  to prove such ${\cal U}, {\cal U}^{\prime}, \theta$ exist.
Then, existence of such $V_i^{(h)}$ is almost obvious.
Since the (Lyons-)It\^o maps associated with the skeleton ODE and the RDE
are continuous, 
the assertion {\rm (ii)} is easy.
\QED

\begin{re}\normalfont
To keep our notations simple, we will often identify  ${\cal U}$ and  ${\cal U}^{\prime}$
through $\theta$.
In such cases we will compute ODE/RDE/SDE on ${\mathbb R}^d$
with the coefficients $V_i^{(h)}$ and the initial point $x \in {\mathbb R}^d$.
The solution of SDE (\ref{sdeX.def}) associated with $V_i^{(h)}$
is denoted by $X^{\ve, (h)}$.
(This should not be confused with $X^{\ve, h}_t := X^{\ve}(t,x,w +h/\ve)$.)
Likewise,  $\phi (k), \Phi ({\bf w}, \lambda)$
associated with $V_i^{(h)}$
are denoted by 
$\phi^{(h)} (k), \Phi^{(h)} ({\bf w}, \lambda)$, respectively.
The ordinary terms and the remainder term
in the expansion of $X^{\ve, (h)}(t,x,w +k/\ve)$
in (\ref{asympX.1}) and (\ref{tylr_rp.eq})
are denoted by 
$f_j^{(h)} (w;k)$, $Q^{\ve, (h)}_{ j+1} (w;k)$
and 
$\hat{f}_j^{(h)} ({\bf w};k)$, $\hat{Q}^{\ve, (h)}_{ j+1} ({\bf w};k)$, respectively.
\end{re}


Next, we take the normal coordinate around $a \in {\cal N}$.
Then, there exists a diffeomorphism 
from a neighborhood of $a$ to a neighborhood of $0$ in $T_a {\cal N}$.
Note that $T_a {\cal N}$ is a vector space with an inner product.
Through this diffeomorphism, the delta function at $a \in {\cal N}$
corresponds to 
the delta function at $0 \in T_a {\cal N}$.
Moreover, since the tangent map of this diffeomorphism at $a$ 
is isometric and so is its transpose,
the determinant and the eigenvalues of 
the deterministic Malliavin covariance matrix of $\phi_1 (h)$ for $h \in {\cal K}_a$
remains the same 
even if they are calculated in the normal coordinate.

Let $h \in {\cal K}_a^{min}$
and let ${\cal U}, {\cal U}^{\prime}$ and $\theta$ be as 
in Lemma \ref{lm.path_nbh}.
(They all depend on $h$.)
If we take a small neighborhood of $\theta (\phi_1 (h) ) \in {\cal U}^{\prime}$,
then $\Pi$ (viewed as a map from ${\cal U}^{\prime}$) 
maps it 
into the normal coordinate of 
$a \in {\cal N}$.
In a standard way we can find a $C^{\infty}_b$-map $\Pi^{(h)}$
from ${\mathbb R}^d$ to $T_a {\cal N}$
which agrees with $\Pi$
when restricted to  the neighborhood of $\theta (\phi_1 (h) )$.
Consequently, for $\ve (h) >0$ and $\gamma (h) >0$
 small enough,  it holds that
\[
\Pi^{(h)} [x + \Phi^{(h)} (\ve {\bf W}, \lambda^{\ve})^1_{0,1} ]
= Y^{\ve}(1,x,w) \in T_a {\cal N}
 \quad
 \mbox{
 if $\ve \in [0, \ve (h)]$
and
$\ve {\bf W} \in U_{h, \gamma (h)} $.}
\]
(Precisely, this means 
that 
$Y^{\ve}(t,x,w)$ takes values in the normal coordinate of $a$ 
if $\ve \in [0, \ve (h)]$
and
$\ve {\bf W} \in U_{h, \gamma (h)} $
and the value is equal to the left hand side.)

Asymptotic expansion of 
$$
\Pi^{(h)} [x + \Phi^{(h)} (\tau_k (\ve {\bf W} ), \lambda^{\ve})^1_{0,1} ] 
=
\Pi^{(h)} [ X^{\ve, (h)} (1, x, w + \frac{k}{\ve})]
$$
was already done in Lemma \ref{lm.tnki_F}.
We  write $g_j^{(h)} (w;k)$ and $R_{j+1}^{\ve, (h)}  (w;k)$
for 
$g_j^{F} (w;k)$ and $R_{j+1}^{\ve, F}  (w;k)$ in Lemma \ref{lm.tnki_F}, 
respectively (with $F = \Pi^{(h)}$).
Thus, local analysis of the solution of RDE on manifold ${\cal M}$
was reduced to that in the linear setting.

\begin{re}\normalfont
The definition of ${\cal K}_a^{min}$ depends on the coefficient 
vector fields of the skeleton ODE.
If they are replaced by new vector fields $V_i^{(h)}~(0 \le i \le r)$,
${\cal K}_a^{min}$ may change and  
 $k \in {\cal K}_a^{min}$ that is close enough to $h$
may not belong to ${\cal K}_a^{min, (h)}$ anymore.
However, we need not worry about this for the following reasons:
Such  $k$ is still a local minimum of 
$k \mapsto \| k\|_{{\cal H}}^2 /2$ subject to 
$\psi^{(h)}_1 (k) =\Pi^{(h)} ( \phi^{(h)}_1 (k) )=0 \in T_a {\cal N}$
and
all  the computations for the skeleton ODE, but one,  require 
just the local minimum property 
via the Lagrange multiplier method
(So $k$ need not be a global minimizer of this conditional minimal problem
in the proof of the asymptotic expansion.)
The only exception where the global minimum property 
is really used is
the proof of the large deviation upper bound,
in which we will compute original $\psi_1$ and ${\cal K}_a^{min}$, 
not $\psi^{(h)}_1$ and  ${\cal K}_a^{min, (h)}$.
\end{re}


Now we compute the second order terms.
Let 
$h \in {\cal K}_a^{min}$ and take ${\cal U}$ and 
$V_i^{(h)}~(0 \le i \le r)$ as in Lemma \ref{lm.path_nbh}.
Suppose that $k$ is sufficiently close to $h$.
The Lagrange multiplier $q (k) \in T^*_a {\cal N}$, 
which is  
independent of $h$, is now an $n$-dimensional row vector 
and (\ref{Lag.def}) holds
with ${\cal V}$  replaced by ${\mathbb R}^n \cong T_a {\cal N}$.

We can compute $ D \phi^{(h)}_t (k)\la l  \ra$
and $D^2 \phi^{(h)}_t (k)\la l, \hat{l}  \ra$
in the same way as in
(\ref{D1phi.eq}) and (\ref{D2phi.eq}), respectively
(if $V_i$ are replaced by  $V_i^{(h)}$).

Since $k$ and $h$ are close, 
$\Pi (\phi_1 (k)) =\Pi^{(h)} (\phi^{(h)}_1 (k))$
and 
we do not distinguish them in this paragraph.
We denote by $\nabla$ the usual gradient for maps from 
an open subset of ${\mathbb R}^d$ to $T_a {\cal N} \cong {\mathbb R}^n$.
Then, we have
\[
D \psi_1^{(h)} (k)\la l\ra
= (\nabla \Pi)(\phi_1^{(h)} (k))  \la D_l \phi_1^{(h)}  (k)\ra
=
(\Pi_*)_{y} D_l \phi_1^{(h)} (k)
\qquad
(\mbox{with $y = \phi_1^{(h)}(k)$})
\]
and
\begin{equation}\label{2df_mnf.eq}
D^2 \psi_1^{(h)} (k) \la l, l' \ra
=
(\nabla \Pi)(\phi_1^{(h)} (k))  \la D^2_{l, l'} \phi_1^{(h)} (k)\ra
+
(\nabla^2 \Pi)(\phi_1^{(h)} (k)) 
 \la D_{l} \phi_1^{(h)} (k), D_{l'} \phi_1^{(h)} (k) \ra.
\end{equation}
Here, we wrote $D_l \phi_1^{(h)}  (k) = D \phi_1^{(h)}  (k) \la l \ra$
and 
$D^2_{l, l'} \phi_1^{(h)} (k) = D^2 \phi_1^{(h)} (k) \la l, l'\ra$
for simplicity.

Note that the same computations
as in (\ref{hess.cmp1}) and  (\ref{hess.cmp2})
 still hold  if $k$ is close to $h$.
 Hence, if $O_{h}$ is a sufficiently small neighborhood 
of $h$ in  ${\cal K}_a^{min}$,
then $\sup \{ C^{(h)} (k) \mid k \in O_{h}\} <1/2$ instead of 
 (\ref{cmax1/2.ineq}).


For the same reason as before, if $k$ is sufficiently close to $h$,
\begin{equation}\nn
\langle q  (k), g_2^{(h)}   (\pi^k w;k)\rangle  
- {\mathbb E} [\langle q(k), g_2^{(h)}  (\pi^k w;k)\rangle ]
\end{equation}
belongs to ${\cal C}_2$,
the second order homogeneous Wiener chaos,
and corresponds to the Hilbert-Schmidt bilinear form
\begin{equation}\nn
\frac12 \bigl\la q(k), D^2\psi_1^{(h)}(k) \la \pi^k \bullet , \pi^k \star \ra \bigr\ra.
\end{equation}
To check the correspondence, 
just apply $(1/2) D^2$ to $g_2^F (\pi^k w;k)$ with $F =\Pi^{(h)}$
 in Lemma \ref{lm.tnki_F}, {\rm (i)} and compare it with (\ref{2df_mnf.eq}).

From this and  $\sup \{ C^{(h)} (k) \mid k \in O_{h}\} <1/2$,
we have the following
under {\bf (A1)'} and {\bf (B1)}--{\bf (B3)}:
For any $h \in {\cal K}^{min}_a$,
there exits $c_1 =c_1 (h) >1$
and a neighborhood $O_{h}$ of $h$ in ${\cal K}^{min}_a$
such that
\begin{equation}\label{expint1.mfd}
\sup_{ k \in O_{h} }
{\mathbb E} \Bigl[
\exp ( c_1 \langle q(k), g_2^{(h)} (w;k)\rangle )
\delta_0 ( g_1^{(h)} (w;k),  {\bf i}_k \la w\ra)
\Bigr]
<\infty.
\end{equation}
Here, $\delta_0$ is
 the Dirac delta function on ${\mathbb R}^n \times T_k {\cal K}^{min}_a$.
(The proof is essentially the same as in Lemma \ref{lm.tw65}.)
Moreover, from (\ref{expint1.mfd}) we obtain the following:
For any $c_2 \in (1,c_1(h))$, there exists a constant 
$\gamma_1 = \gamma_1 (h) >0$ 
which is independent of $k \in  O_h$ and 
satisfies that
\begin{equation}\label{expint2.mfd}
\sup_{0 <\ve \le 1}\sup_{ k \in  O_{h}  }
{\mathbb E}
\Bigl[
\exp \Bigl( c_2  \bigl\la  q (k),   \frac{ R^{\ve, (h)}_{ 2}( w ;k)   }{ \ve^2}  \bigr\ra \Bigr)  
{\bf 1}_{U_{\gamma_1} }(\ve {\bf W})
~;~
| {\bf i}_k \la w \ra|^2 + \Bigl| \frac{ R^{\ve, (h)}_{1} ( w ;k) }{\ve} \Bigr|^2 
\le \kappa^2
\Bigr] 
<\infty
\end{equation}
for any $\kappa >0$.
(The proof is essentially the same as in Lemma \ref{lm.tw65}.)

%

In a similar way, Lemma \ref{lm.tw63} can be modified.
If $O_h$ is sufficiently small, then
\begin{equation}\label{lamb.mnf}
\lambda (h) 
:=\inf_{k \in  O_h   } 
\inf_{z :  |z|=1} z^* \sigma [\psi_1^{(h)} ](k) z >0.
\end{equation}
Taking $O_h$ smaller if necessary, 
we have the following:
There exist $\gamma_0= \gamma_0 (h), 
\ve_0=\ve_0 (h)  \in (0,1]$ such that
the smallest eigenvalue of 
$\sigma [(R^{\ve, (h)}_{ 1}( \,\cdot\,;k)/\ve, {\bf i}_k ) ] (w)$ is greater than 
$(\lambda (h) \wedge 1)/2$
if  $ \ve {\bf W} \in U_{\gamma_0}$, 
$k \in O_h$ and $0 <\ve \le \ve_0$.
(Here, $\psi_1^{(h)}$ and $R^{\ve, (h)}_{ 1} ( \,\cdot\,;k) /\ve$ are 
 regarded as $T^*_a {\cal N} \cong {\mathbb R}^n$-valued maps.)

%

In the same way as in Lemma \ref{lm.prm}, we have the following lemma.
Note that in this case the constants 
$\lambda, \gamma_0, \gamma_1, \ve_0$ all depend on $h \in {\cal K}_a^{min}$
(which makes no difference, however).
The definition of $\sigma=\sigma_h$ is in (\ref{csigma.eq}).
\begin{lm}\label{lm.prm_mfd}
For any $h \in {\cal K}^{min}_a$, we can find $O_h$ 
(or $\rho = \rho_h:= \sup_{k \in O_h} \| k-h \|_{{\cal H}} >0$),
$\gamma=\gamma_h >0$ 
and $\ve_0 =\ve_{0}(h) \in (0,1]$
such that (\ref{223_1.cond}),   (\ref{223_2.cond}), 
 (\ref{223_5.cond}) and the following three conditions hold.
The first one is:
 \begin{equation}\nn
 \inf_{k \in O_h}
 \inf_{z : \|z\|=1} 
z^*    \sigma [(R^{\ve, (h)}_{ 1}(\,\cdot\, ;k) /\ve, {\bf i}_k ) ](w)  z 
\ge \frac{\lambda (h) \wedge 1}{2}
\qquad
\mbox{if $\ve \in (0,\ve_0]$ and $\ve {\bf W} \in U_{ h-k, 2^{1/4m}\sigma\gamma}$}.
\end{equation}
The second one is:
\begin{equation}
\sup_{0 <\ve \le 1}\sup_{ k \in  O_h}
{\mathbb E}
\Bigl[
\exp \Bigl( \frac{c_2}{\ve^2}  \la  q (k), R^{\ve, (h)}_{ 2}(w ;k)\ra  \Bigr)  
{\bf 1}_{U_{ h-k, \gamma} } (\ve {\bf W})
~;~
| {\bf i}_k |^2 + \Bigl| \frac{ R^{\ve, (h)}_{1}(w ;k)}{\ve} \Bigr|^2 \le \kappa^2
\Bigr] 
<\infty.
\nn
\end{equation}
The last one is:
\begin{eqnarray*}
\{   x + \Phi (  {\bf w} ,  \lambda^{\ve})^1_{0,t} 
\mid  t \in [0,1], \ve \in [0, \ve_0],    {\bf w}  \in U_{h, \gamma} \}
\\
\cup \{ \phi (k)_t \mid t \in [0,1], k \in O_h\}
\subset 
{\cal U}_h \subset {\cal M},
\end{eqnarray*}
where ${\cal U}_h \subset {\cal M}$ is the tubular neighborhood of 
${\rm Im} \phi (h)$ introduced in Lemma \ref{lm.path_nbh}.
\end{lm}

\Proof
The proof is  the same as in Lemma \ref{lm.prm}.
\QED


In what follows let $\gamma_h$ be the constant in Lemma \ref{lm.prm}.
Since $\gamma_h >0$ for all $h \in {\cal K}_a^{min}$,
we can choose finitely many $(h_1, \gamma_1), \ldots, (h_N, \gamma_N)$
and construct asymptotic partition of unity 
$\chi^{\ve}$ and $\tilde{\chi}^{\ve}_{\nu}$
exactly in the same way as in 
(\ref{apu1.def})--(\ref{apu2.def}).
(Here, $ \gamma_{\nu}$ stands for $\gamma_{h}$ with $h = h_{\nu}$.)
Note that the definitions of 
$\chi^{\ve}$ and $\tilde{\chi}^{\ve}_{\nu}$
depend only on $ (h_{\nu}, \gamma_{\nu})_{1\le \nu \le N}$.
Below we set 
$\ve'_0 := \min_{1\le \nu \le N} \ve_0(h_{\nu}) \in (0,1]$
and 
$\lambda' := \min_{1\le \nu \le N} \lambda(h_{\nu}) >0$.


As in the Euclidean case, 
contribution of the rough paths away from ${\cal K}_a^{min}$ is negligible.
To prove it, 
we need to modify Theorem \ref{tm.ub_ldp}, 
the large deviation upper estimate, 
to the manifold case.
(Once this is done, a manifold version of Theorem \ref{lm.ldp_up}
can be obtained in the same way.)

The keys of the proof of Theorem \ref{tm.ub_ldp}  
are Kusuoka-Stroock's moment estimate (\ref{ks_bnd.ineq})  
of $\det \sigma[ Y^{\ve}_1]^{-1}$
and
the integration by parts formula for Watanabe distributions.
Both of them hold in the manifold case (see (\ref{ibp.mfd})). 
Therefore, we can prove Theorem \ref{tm.ub_ldp} 
in the manifold case, too.
Consequently, $I_1 (\ve)$ does not contribute to 
the asymptotic expansion.


\noindent
{\it Proof of Theorem \ref{tm.main_mfd}. }
It is sufficient to compute 
$$
I_2 (\ve) = \sum_{\nu=1}^N
{\mathbb E}[\tilde\chi^{\ve}_{\nu} G^{\ve}  \delta_a (Y^{\ve}_1) ]
$$ 
as $\ve \in (0, \ve'_0]$ tends to zero.
Here, $\delta_a$ is the Dirac measure at $a \in {\cal N}$
and $Y_1^{\ve} = \Pi ( X_1^{\ve})$
with $(X_t^{\ve})$ being the solution of the ${\cal M}$-valued SDE (\ref{sde.mfd}).

Let $\eta \in C_0^{\infty} (T_a{\cal N}, {\mathbb R})$
be as in (\ref{sonyu.eq}) in Section 7.
Then, defined by $\eta_l (y) = l^n \eta (l y)$,
$\{ \eta_l \}_{l \ge 1}$ 
approximates $\delta_0$,
where $\delta_0$ is the delta function at $0$ on 
$T_a {\cal N} \cong {\mathbb R}^n$ with respect to the Lesbegue measure
associated with the inner product on $T_a {\cal N}$.
Denote by $\hat{\eta}_l \in C_0^{\infty} ({\cal N}, {\mathbb R})$
the  function which corresponds to $\eta_l$ 
via the exponential map at $a$.
Then, $\{ \hat\eta_l \}_{l\ge 1}$
approximates the delta function at $a$ on ${\cal N}$.
Therefore,
\begin{eqnarray}
{\mathbb E}[\tilde\chi^{\ve}_{\nu} G^{\ve}  \delta_a (Y^{\ve}_1) ]
 =
 \lim_{l \to \infty}
{\mathbb E}[\tilde\chi^{\ve}_{\nu} G^{\ve}  \hat\eta_l (Y^{\ve}_1) ]
=
\lim_{l \to \infty}
{\mathbb E}[\tilde\chi^{\ve}_{\nu} G^{\ve}  \eta_l (Y^{\ve, (\nu)}_1) ]
\label{onnaji.eq}
\end{eqnarray}
since $\tilde\chi^{\ve}_{\nu} =0$ if $\ve{\bf W} \notin U_{h_{\nu}, \gamma_{\nu}}$.
Here $Y^{\ve, (\nu)}_1$ is a shorthand for 
\[
Y^{\ve, (h)}_1
= 
\Pi^{(h)} [ X^{\ve, (h)} (1, x, w)]
\qquad
\mbox{with $h =h_{\nu}$}.
\]
%
%
%
%
%
%
Note that the right hand side of (\ref{onnaji.eq})
is essentially the same as the left hand side of (\ref{sonyu.eq}).
Thus,
the asymptotic problem has been reduced to the one on the Euclidian space.
Note also that we need not prove non-degeneracy of $\sigma [Y_1^{\ve, (h)} ]$
because we use Propositions \ref{pr.comp1} and \ref{pr.comp2}.

Define ${\cal E}_j^{(\nu)} (w;k) \in \tilde{\mathbb D}_{- \infty}$
in the same way as in (\ref{calE.def}).
(This time it depends on $\nu$ (or $h_{\nu}$), however.)
By the same computation as in (\ref{asFin.eq}), we have
the following asymptotics: for any $L \in {\mathbb N}$,
\begin{eqnarray}
\lefteqn{
{\mathbb E}[\tilde\chi^{\ve}_{\nu} G^{\ve}  \delta_a (Y^{\ve}_1) ]
}
\nn\\
&=&
\int_{ O_{\nu}} 
{\mathbb E}
\Bigl[
\exp (  \la  q(k), g_2^{(\nu)} (w;k)  \ra  )
\Bigl(
\sum_{j=0}^{L} \ve^{j} {\cal E}_{j}^{(\nu)} (w;k) + O (\ve^{L+1})
\Bigr)
\tilde\chi^{\ve}_{\nu} (w +\frac{k}{\ve})
\Bigr]
\omega(dk)
\nn\\
&=&
\int_{ O_{\nu}} 
{\mathbb E}
\Bigl[
\exp (  \la  q(k), g_2^{(\nu)} (w;k)  \ra  )
\Bigl(
\sum_{j=0}^{L} \ve^{j} {\cal E}_{j}^{(\nu)} (w;k) 
\Bigr)
\Bigl(
\sum_{l=0}^{L} \ve^{l} {\cal D}_{l}^{(\nu)} (w;k) 
\Bigr)
\Bigr]
\omega(dk)
+O (\ve^{L+1})
\nn
\end{eqnarray}
as $\ve \searrow 0$.
Note that
the asymptotic expansion inside $\int _{ O_{\nu}}\omega (dk)$
is uniform in $k \in {O_{\nu}}$.
$
$
Vanishing of $c_{2j+1}$ under {\bf (C2)} can be shown in the same way 
as in the Euclidean case.
This completes the proof of Theorem \ref{tm.main_mfd}.
\toy


\vspace{15mm}

\begin{flushleft}
\begin{tabular}{ll}
Yuzuru \textsc{Inahama}
\\
Graduate School of Mathematics,   Kyushu University,
\\
Motooka 744, Nishi-ku, Fukuoka, 819-0395, JAPAN.
\\
Email: {\tt inahama@math.kyushu-u.ac.jp}
\end{tabular}
\end{flushleft}


\begin{flushleft}
\begin{tabular}{ll}
Setsuo \textsc{Taniguchi}
\\
Faculty of Arts and Science,   Kyushu University,
\\
Motooka 744, Nishi-ku, Fukuoka, 819-0395, JAPAN.
\\
Email: {\tt se2otngc@artsci.kyushu-u.ac.jp}
\end{tabular}
\end{flushleft}

\end{document}